\documentclass[11pt,
twoside]{amsart}
\usepackage{amsmath}
\usepackage{amssymb} 
\usepackage{amsfonts} 
\usepackage{latexsym, amsthm, mathrsfs}
\usepackage[latin1]{inputenc}
\usepackage{fancybox}
\usepackage{amscd, upgreek}
\input{xy}
\usepackage[all]{xy}
\numberwithin{equation}{section}
\newtheorem{theorem}{Theorem}[section]
\newtheorem{lemma}[theorem]{Lemma}

\theoremstyle{definition} %fettes label, aber nich kursiv

\theoremstyle{plain} %fettes label. text kursiv

\newtheorem{proposition}[theorem]{Proposition}
\newtheorem{corollary}[theorem]{Corollary}

\newtheorem*{maintheorem*}{Main Theorem}
\newtheorem*{conjecture*}{Conjecture}
\newtheorem*{theorem*}{Theorem}
\newtheorem*{proposition*}{Proposition}
\newtheorem{definition}[theorem]{Definition}

\newtheorem{notation/observation}[theorem]{Notation/Observation}

\theoremstyle{remark}  %kursives label. text roman
\newtheorem{remark}[theorem]{Remark}
\newtheorem*{remarks*}{Remarks}
\newtheorem*{remark*}{Remark}

\newtheorem*{claim*}{Claim}
\newcommand{\nc}{\newcommand}
\nc{\nothing}[1]{}
\nc{\comment}[1]{#1}
\nothing{   %dies kann fuer alte papers wieder aufleben!
\nc{\dom}{{\rm dom}} \nc{\card}{{\rm card}} \nc{\lh}{{\rm lh}}
\nc{\lgg}{{\rm lg}} \nc{\rge}{\mbox{\rm range}} \nc{\cf}{{\rm cf}}
nc{\nex}{\mbox{\rm next}} \nc{\uhr}{\restriction} \nc{\supt}{{\rm
supt}} \nc{\supp}{{\rm supp}} \nc{\Lim}{{\rm Lim}} \nc{\Leb}{{\rm
Leb}} \nc{\modd}{{\rm mod}} \nc{\RO}{{\rm RO}} \nc{\prob}{{\rm
Prob}} }

\nc{\On}{{\rm On}} \nc{\Ord}{{\rm On}}
\nc{\nco}{\DeclareMathOperator} \nco{\rk}{rk} \nco{\order}{o}
\nco{\ppower}{pp}
\nco{\pcf}{pcf} %possible cofinality
\nco{\tcf}{tcf} %true cofinality
\nco{\tlim}{tlim} %true lim
\nco{\limtext}{lim} %true lim
\nco{\prodt}{{\textstyle \prod}}
\nco{\symdiff}{\triangle} \nco{\dom}{dom} \nco{\card}{card}
\nco{\lh}{lh} \nco{\lt}{lt} \nco{\lgg}{lg} \nco{\hgt}{ht}
\nco{\rge}{range} \nco{\otp}{otp} \nco{\trunk}{tr} \nco{\cf}{cf}

\nco{\cfstar}{\cf^\ast}

\nco{\nex}{next} \nc{\uhr}{\restriction} \nco{\reduction}{red}
\nco{\supt}{supt} \nco{\supp}{supp} \nco{\Lim}{Lim} \nco{\Leb}{Leb}
\nco{\modd}{mod} \nco{\invariant}{inv} \nco{\id}{id} \nco{\RO}{RO}
\nco{\poss}{pos}
\nco{\Inc}{Inc} %for the canonical forms
\nco{\Ge}{Ge}
\nco{\hdrop}{\hat{drop}} \nco{\Gen}{Gen} \nco{\Property}{Pr}
\nco{\GT}{GT} \nco{\refl}{refl} \nco{\nmg}{nmg} \nco{\en}{en}
\nco{\Tc}{Tc}
\nc{\potom}{\ensuremath{{\cal P}(\omega)}}
\nc{\potinf}{\ensuremath{[\omega]^\omega}}
\nc{\pfin}{\ensuremath{{\cal P}(\omega)/{\rm fin}}}

\nc{\potfin}{\ensuremath{[\omega]^{<\omega}}}
\nc{\inn}{\ensuremath{{\omega^{\uparrow \omega}}}}
\nc{\baire}{{}^\omega \omega} \nc{\bair}{{}^\omega \omega}
\nc{\hoch}{^{<\omega}} \nc{\hocho}{^{\omega}} \nco{\tree}{tree}

\nc{\prooff}[1]{{\bf Proof} of #1:}
\nc{\proofend}{\makebox{} \hfill ${\bf \dashv}$ \\}
\nc{\proofendof}[1]{\makebox{} \hfill $\boldmath{\dashv}_{\rm #1}$
\\} \nc{\beq}{\begin{eqnarray*}} \nc{\eeq}{\end{eqnarray*}}
\nc{\bde}{\begin{list}} \nc{\ede}{\end{list}}

\newenvironment{myrules}
{\begin{list}{} {
 \setlength{\leftmargin}{0.8cm}
 \setlength{\labelwidth}{0.8cm}
 \setlength{\labelsep}{0.2cm}
 \setlength{\parsep}{0.5ex plus 0.2ex minus 0.1 ex}
 \setlength{\itemsep}{0.3ex plus 0.2 ex minus 0ex}
}}{\end{list}}

\newenvironment{myrules1}
{\begin{list}{} {
 \setlength{\leftmargin}{1.0cm}
 \setlength{\labelwidth}{0.8cm}
 \setlength{\labelsep}{0.3cm}
 \setlength{\parsep}{0.5ex plus 0.2ex minus 0.1 ex}
 \setlength{\itemsep}{0.5ex plus 0.2 ex minus 0ex}
}}{\end{list}}

{\end{list}}

\newcounter{subalph}
{\end{list}}

\newcommand{\greek}[1]{\ifthenelse{\value{#1}=1}{\mbox{$\alpha$}}%
  {\ifthenelse{\value{#1}=2}{\mbox{$\beta$}}{%
   \ifthenelse{\value{#1}=3}{\mbox{$\gamma$}}{%
   \ifthenelse{\value{#1}=4}{\mbox{$\delta$}}{%
   \ifthenelse{\value{#1}=5}{\mbox{$\varepsilon$}}{%
   \ifthenelse{\value{#1}=6}{\mbox{$\zeta$}}{%
   \ifthenelse{\value{#1}=7}{\mbox{$\eta$}}{%
   \ifthenelse{\value{#1}=8}{\mbox{$\theta$}}{%
   \ifthenelse{\value{#1}=9}{\mbox{$\iota$}}{%
   \ifthenelse{\value{#1}=10}{\mbox{$\kappa$}}{%
   \ifthenelse{\value{#1}=11}{\mbox{$\lambda$}}{%
   \ifthenelse{\value{#1}=12}{\mbox{$\mu$}}{%
   \ifthenelse{\value{#1}=13}{\mbox{$\nu$}}{%
   \ifthenelse{\value{#1}=14}{\mbox{$\xi$}}{%
   \ifthenelse{\value{#1}=15}{\mbox{$\rm o$}}{%
   \ifthenelse{\value{#1}=16}{\mbox{$\pi$}}{%
   \ifthenelse{\value{#1}=17}{\mbox{$\varrho$}}{%
   \ifthenelse{\value{#1}=18}{\mbox{$\sigma$}}{%
   \ifthenelse{\value{#1}=19}{\mbox{$\tau$}}{%
   \ifthenelse{\value{#1}=20}{\mbox{$\upsilon$}}{%
   \ifthenelse{\value{#1}=21}{\mbox{$\varphi$}}{%
   \ifthenelse{\value{#1}=22}{\mbox{$\chi$}}{%
   \ifthenelse{\value{#1}=23}{\mbox{$\psi$}}{\mbox{$\omega$}%
  }}}}}}}}}}}}}}}}}}}}}}}}%24 Zuklammern

\newcounter{subgreek}
{\end{list}}

\newcounter{subarabic}
{\end{list}}

\newcounter{subroman}
{\end{list}}

\newcount\skewfactor

\def\mathunderaccent#1#2 {\let\theaccent#1\skewfactor#2
\mathpalette\putaccentunder}
\def\putaccentunder#1#2{\oalign{$#1#2$\crcr\hidewidth
\vbox
to.2ex{\hbox{$#1\skew\skewfactor\theaccent{}$}\vss}\hidewidth}}
\def\name{\mathunderaccent\tilde-3 }

\nc{\nname}{\name}

\nc{\even}{\ensuremath{\rm Even}} \nc{\odd}{\ensuremath{\rm Odd}}

\nc{\al}{$\alpha$\  } \nc{\om}{\omega}
\nc{\omm}{\ensuremath{\omega_1}} \nc{\ep}{\varepsilon}
\nc{\tk}{\tilde{K}}
\nc{\concat}{{}^\frown}   %im math mode
\nc{\force}{\Vdash} \nc{\fb}{f_{\overline{M}}} \nc{\such}{\, : \,}
\nc{\la}{\langle}
\nc{\ra}{\rangle}
\nc{\meager}{\ensuremath{{\cal M}}} 
\nc{\lebesgue}{\ensuremath{{\cal N}}} 
\nc{\nulll}{\ensuremath{{\cal N}}}
\nc{\ksigma}{\ensuremath{{\bf K}_\sigma}}
\nc{\ideal}{\ensuremath{{\cal I}}} \nc{\ga}{\ensuremath{\frak a}}
\nc{\AAA}{{\cal A}}   %im math mode
\nc{\gc}{\ensuremath{\frak c}} \nc{\gs}{\ensuremath{\frak s}}
\nc{\gh}{\ensuremath{\frak h}} \nc{\gd}{\ensuremath{\frak d}}
\nc{\gb}{\mathfrak{\lowercase{b}}} \nc{\gro}{\ensuremath{\frak g}}
\nc{\gu}{\ensuremath{\frak u}} \nc{\gr}{\ensuremath{\frak r}}
\nc{\gt}{\ensuremath{\frak t}} \nc{\fff}{\ensuremath{\frak f}}
\nc{\gm}{\ensuremath{\mathfrak{mcf}}}
\nc{\gge}{\ensuremath{\mathfrak e}}
\nc{\cfupro}{\ensuremath{\cf(\upro)}}
\nc{\cfvpro}{\ensuremath{\cf(\vpro)}} \nc{\gp}{\ensuremath{\frak p}}
\nc{\gk}{\ensuremath{\frak k}}

\nc{\add}{\mbox{\ensuremath{{\rm add}}}}
\nc{\cov}[1]{\mbox{\ensuremath{{\rm cov}(#1)}}}
\nc{\unif}[1]{\mbox{\ensuremath{{\rm unif}(#1)}}}
\nc{\cof}[1]{{\mbox{\ensuremath{\rm cof}(#1)}}}

\nc{\addd}[2]{\mbox{\ensuremath{{\rm add}^{#1}(#2)}}}   %fuer ch4.
\nc{\covv}[2]{\mbox{\ensuremath{{\rm cov}^{#1}(#2)}}}   %unter #1 steht
\nc{\uniff}[2]{\mbox{\ensuremath{{\rm unif}^{#1}(#2)}}} %das Modell
\nc{\coff}[2]{{\mbox{\ensuremath{\rm cof}^{#1}(#2)}}}

\nc{\cd}{Cicho\'n's Diagram}

\nc{\MA}{\mbox{\sf MA}} \nc{\PFA}{\mbox{\sf PFA}}
\nc{\OCA}{\mbox{\sf OCA}} \nc{\GCH}{\mbox{\sf GCH}}
\nc{\CH}{\mbox{\sf CH}} 
\nc{\zfc}{\mbox{\sf ZFC}}
\nc{\sch}{\mbox{\sf SCH}} \nc{\ZF}{\mbox{\sf ZF}}
\nc{\zc}{\mbox{\sf ZC}}
\nc{\NCF}{\mbox{\sf NCF}} %! nicht mehr variabel programmiert,
\nc{\FD}{\mbox{\sf FD}}   %in den \indices ohne makro, damit sie
\nc{\SFT}{\mbox{\sf SFT}} \nc{\fourG}{\mbox{\sf 4G}}
\nc{\fourI}{\mbox{\sf 4I}} \nc{\past}{\ ;{\rm past}\;}
\nc{\Borelhood}{Borel measurability} %according to Kechris
\nc{\Pieinseins}{\mbox{${\bf \Pi}^1_1$}}
  \nc{\seinseins}{\mbox{${\bf\Sigma}^1_1$}}
\nc{\seinszwei}{\mbox{${\bf\Sigma}^1_2$}}
\nc{\seinsdrei}{\mbox{${\bf\Sigma}^1_3$}}
\nc{\Deleinszwei}{\mbox{${\bf\Delta}^1_2$}}

\nc{\up}{\ensuremath{{\cal U}\mbox{\ensuremath{\rm -prod}}\,\omega}}
\nc{\upp}{\ensuremath{{\cal U}'\mbox{\ensuremath{\rm
-prod}}\,\omega}} \nc{\upro}{\ensuremath{{\cal
U}\mbox{\ensuremath{\rm -prod}}\,\om}}
\nc{\fupro}{\ensuremath{f({\cal U})\mbox{\ensuremath{\rm
-prod}}\,\om}} \nc{\vpro}{\ensuremath{{\cal V}\mbox{\ensuremath{\rm
-prod}}\,\om}} \nc{\fpro}{\ensuremath{{\cal F}\mbox{\ensuremath{\rm
-prod}}\,\om}}

\nc{\cff}[1]{{\text{cf}\,(#1)}}           %cf mit richtiger parameter
\nc{\cu}{\ensuremath{\cal U}}             %Vorsicht, dies gibt im
\nc{\ai}{\ensuremath{\forall^\infty}}     %geaendert ensuremath dazu
\nc{\ei}{\ensuremath{\exists^\infty}}     %auch geaendert
\nc{\ww}{\ensuremath{\omega^\omega}}      %auch geaendert
\nc{\gw}{groupwise dense}

\nc{\kk}{car\-dinal cha\-rac\-teris\-tic} \nc{\joker}{\ast}
\nc{\gtc}{Galois-Tukey connection} %Vorsicht: heisst generalized

\nc{\av}[1]{{\rm Av}_{#1}} \nc{\eps}{\varepsilon}
\renewcommand{\epsilon}{\varepsilon}

\nc{\n}{{\bf n}}                 %for the blueprints
\nc{\m}{{\bf m}}

\nc{\marginparr}[1]{\marginpar{#1}}
\nc{\footnoteee}{} % I think the first footnote is no more a question
\nc{\footnotee}{}  % for some ad libitum notes

\newcommand{\cal}{\mathcal}

\nc{\divs}{{c_0 \setminus \ell^1}} \nc{\divser}{(\divs,
\leq^*)/\thickapproy} \nc{\bfin}{\RO(\pfin
\setminus\{0\},\subseteq^*)} \nc{\bdivser}{\RO(\divser)}
\nc{\inc}{{\rm INC}} \nc{\com}{{\rm COM}}
\nc{\thickapproy}{\makebox{}\!\!\thickapprox}
\nc{\approy}{\makebox{}\!\!\approx} \nc{\lessi}{\leqslant}
\nc{\gessi}{\geqslant} \nc{\interior}[1]{{\rm int}(#1)}
\nc{\closure}[1]{{\rm cl}(#1)} \nc{\Vo}{Vojt\'a\v{s}}

\nc{\precedeseq}{\leq^*} %%%take your favorite partial order!!!!
\nc{\precedes}{\prec} \nc{\stronger}{\leqslant_{\bf P}}
\nc{\underlline}[1]{\hat{#1}} \nc{\PO}{{\bf P}}
\nc{\charak}{\text{ch}} \nc{\symom}{{\rm{Sym}(\omega)}}
\nc{\needed}{needed\ } \nc{\neededc}{needed} \nc{\Needed}{Needed\ }
\nc{\wneeded}{weakly needed\ } \nc{\Wneeded}{Weakly needed\ }
\nc{\wneededc}{weakly needed}

\nco{\dcl}{dcl}
 \nco{\WCR}{WCR}
 \nco{\TCR}{TCR}
 \nc{\mup}{m_{\rm up}} \nc{\mdn}{m_{\rm dn}}
\nco{\may}{may}
\nco{\aver}{av} % neu fuer F427
\nco{\norm}{nor} % neu fuer F427  SCHOENER ohne fett
\nco{\val}{val} % neu fuer F427
\nco{\dis}{dis} % neu fuer F427-->MdSh:731
\nco{\basis}{basis}
\nco{\pos}{pos}
\nco{\spec}{spec}
\nc{\err}{\mbox{err}}
\nc{\eee}{\mbox{e}}
\nco{\Expect}{Exp}
\nco{\rt}{rt}
\nco{\pr}{pr}
\nco{\suc}{succ}
\nco{\splitt}{split}
\nco{\halv}{h}
\nco{\Add}{Add}
\nco{\Cov}{Cov}
\nco{\Unif}{Unif}
\nco{\Cof}{Cof}
\nco{\htt}{ht}

\nco{\tr}{tr}

\nco{\Levy}{Levy}

\nc{\bbforcing}{\mathbb A} \nc{\itername}{\mathfrak q}
\nc{\iterp}{\mathfrak p} \nc{\iterq}{\mathfrak q} \nc{\invcm}{\rm
inv_{cm}} \nc{\invcf} {\rm inv_{cf}} \nc{\invgm}{\rm inv_{gm}}

\nc{\Q}{\mathbb Q}
\nc{\subsetsim}{\underset{\raise0.6em\hbox{$\sim$}}{\subset}}
\newcommand{\subsim}{\underset{\raise20pt\hbox{$\rightarrow$}}{\rightarrow}}
\newcommand{\ssim}{\overset{\raise-40pt\hbox{$\leftarrow$}}{\subsim}}
\nc{\rest}{\restriction}

\nc{\sE}{\mathscr E} \nc{\cE}{\mathscr E}

\nc{\F}{\mathbb F} \nc{\bF}{\mathbb F} \nc{\bR}{\mathbb R}
\nc{\M}{\mathbb M} \nc{\bQ}{\mathbb Q} \nc{\bP}{\mathbb P}

\nc{\cP}{\mathscr P} \nc{\cU}{\mathscr U} \nc{\cV}{\mathscr V}
\nc{\cC}{\mathscr C} \nc{\cD}{\mathscr D} \nc{\cG}{\mathscr G}
\nc{\cF}{\mathscr F} \nc{\cA}{{\mathscr A}} \nc{\cB}{\mathscr B}
\nc{\cR}{\mathscr R} \nc{\cT}{\mathcal T}
\nc{\cK}{\mathscr K}\nc{\cO}{\mathscr O}

\nc{\cI}{\mathcal I} \nc{\cJ}{\mathcal J} \nc{\cH}{\mathcal
H}\nc{\cS}{\mathcal S}

\nc{\cW}{\mathscr W}

\nc{\god}{\mathfrak{od}} % od
\nc{\mcf}{\mathfrak{\lowercase{mcf}}} \nc{\roth}{[\omega]^{\omega}}
\nco{\last}{last} % fuer F767
\nco{\filter}{fil} \nco{\semifilter}{semi} \nco{\sfil}{sfil}
\nco{\ssfil}{ssfil}

\nco{\upwards}{up} \nco{\inter}{inter} \nco{\FU}{FU} \nco{\FP}{FP}
\nco{\OP}{OP} \nco{\set}{set} \nco{\seq}{seq}

\nco{\acc}{acc}

\nc{\forks}{\underset{\raise0.6em\hbox{$\smile$}}{\mid}}
\renewcommand{\setminus}{\smallsetminus}
\nc{\bV}{{\bf V}} \nco{\fil}{fil} \nco{\CFF}{CFF} \nco{\rfl}{rfl}

\nc{\fin}{\emptyset}
\nc{\ba}{\bar{a}}
\nc{\bb}{\bar{b}}
\nc{\bc}{{\bf c}} %%%Vorsicht

\nc{\fc}{{\bf c}} \nc{\fd}{{\bf d}}

 \nc{\bd}{\bar{d}}
  \nc{\be}{\bar{e}}
\nc{\fto}{\ensuremath{\omega^{\omega,\text{fto}}}}

\nc{\bx}{{\bf x}}
\nc{\by}{{\bf y}}
 \nc{\bw}{\bar{w}}

\nc{\bg}{{\bf g}}

\nc{\bK}{{\bf K}} \nc{\bq}{\bar{q}}

\nc{\bp}{\bar{p}} %Vorsicht, in manchen papers durch \bf p interpretieren lassen

 \nc{\bG}{{\bf G}}

\nco{\bH}{H}

\nco{\bB}{{\bf B}}
 \nc{\sk}{\smallskip} \nc{\mk}{\medskip}
\nc{\reri}{\Rsh} %Harpoon zuerst aufwaerts dann nach rechts fuer
 \renewcommand{\phi}{\varphi}
\nc{\ord}{{\rm ord}}
\begin{document}

%-----------Title-------------------
\title[Preserving Souslin Trees]{Many countable support iterations of
proper forcings preserve Souslin trees}

\author{Heike Mildenberger and Saharon Shelah}
\address{Einstein Institute of Mathematics,
The Hebrew University, Edmond Safra Campus Givat Ram, Jerusalem
91904, Israel}
\address{Abteilung f\"ur Mathematische Logik, Mathematisches Institut,
 Universit\"at Freiburg, Eckerstr.~1, 79104 Freiburg im Breisgau, Germany}

\email{heike.mildenberger@math.uni-freiburg.de, shelah@math.huji.ac.il}

\thanks{
\\
The first author acknowledges Marie Curie grant PIEF-2008-219292
of the European Union. The second author's research was partially
supported by the United States-Israel Binational Science Foundation
(Grant no. 2002323). This is the second author's publication
no.~973. Corresponding author: H.~Mildenberger.}

\begin{abstract}
We show that many countable support iterations of proper forcings preserve
Souslin trees. We establish sufficient conditions in terms of games and 
we draw connections to other preservation properties. We present a proof of preservation
properties in countable support iterations in the so-called Case A
that does not need a division into forcings that add reals and those who do not.
\end{abstract}

\date{November 16, 2010, Revisions: August 2011, March 2013, August 2013}
\maketitle

 \setcounter{section}{-1}
\section{Introduction}\label{S0}

This work is related to Juh\'asz' question \cite{Miller1993}: ``Does
Ostaszewski's club principle imply the existence of a Souslin
tree?'' We recall the club principle (also written $\clubsuit$):
There is a sequence
$\langle A_\alpha \such \alpha \mbox{ a limit ordinal} < \omega_1\rangle$  with the following properties: 
For every countable limit ordinal $\alpha$,
$A_\alpha$ is cofinal in $\alpha$ and for any uncountable $X  \subseteq \omega_1$
there are stationarily many $\alpha$ with $A_\alpha \subseteq X$.
Such a sequence is called a $\clubsuit$-sequence. 
The club principle was introduced in \cite{ostaszewski}.

Partial positive answers are known: Let $\meager$ denote the ideal of meagre sets.
In every model of the club principle and $\cov{\meager}>\aleph_1$
by Miyamoto \cite[Section~4]{brendle:baumgartnerband} there are
Souslin trees. Brendle showed
\cite[Theorem~6]{brendle:baumgartnerband}: In every model of the club
principle and $\cof{\meager} = \aleph_1$ there are Souslin trees.
In this paper we give examples of models 
satisfying the club principle, the existence of Souslin trees, 
$\cov{\meager}=\aleph_1$ and
$\cof{\meager}=\aleph_2$ (i.e., neither of the sufficient conditions mentioned above holds).

Assume that we start with a ground model satisfying 
$\diamondsuit_{\omega_1}$ and that we force with a proper
 countable support iteration 
$\langle \bP_\alpha,
\name{\bQ}_\beta \such \alpha\leq\kappa,\beta<\omega_2\rangle$
of length $\omega_2$. 
For this scenario in \cite{Mi:clubdistr} we showed: 
If the single step forcings  are suitable forcings from
\cite{RoSh:470} (with finite or countable ${\bH}(n)$, see Section
2.1), then the final model will satisfy the club principle. 
Note that the assumption of
the diamond in the ground model is actually not necessary, since
after $\omega_1$ iteration steps of any forcing with two incompatible conditions
with countable support $\diamondsuit_{\omega_1}$
holds anyway \cite[Ch.~7, Theorem 8.3]{Kunen} and the length of our iterations is  $\omega_2$.

Let us look at the countable support iteration of length
$\omega_2$ of Miller forcing:
According to the mentioned result, after $\omega_1$ many steps we get 
$\diamondsuit_{\omega_1}$ 
 and 
therefore a Souslin tree in the intermediate extension. Lemma~2.1. together with the results in Section~\ref{S4} show
that any countable support iteration of Miller forcing preserves Souslin trees.
Hence after $\omega_2$
many iteration steps there is a Souslin tree.
Moreover by \cite{Mi:clubdistr} the club principle holds.
 It is known
that in the Miller model $\gd= \aleph_2$ (and hence $\cof{\meager}=\aleph_2$)
and $\cov{\meager}=
\aleph_1$.  
A countable support iteration of
length $\omega_2$ of Blass--Shelah forcing gives
another model of $\gd= \aleph_2$ 
and $\cov{\meager}=
\aleph_1$ and the club principle. Blass--Shelah forcing 
is not $\omega$-Cohen preserving (see Def.~\ref{3.1}) and
increases the splitting number (see \cite[Prop.~3.1]{BsSh:242}). 
Besides these two particular examples,
the main technical work in this paper is a study of
the preservation of Souslin trees.

We refer the reader to
\cite{BJ} for the definitions of cardinal characteristics, and to
\cite{Mi:clubdistr} for reading about the the club principle. For background about
properness we refer the reader to \cite{Sh:f} and the more detailed 
introductions in \cite{Goldstern93, abraham:handbook}.
In  forcing notions, $q> p$ means that
$q$ is stronger than $p$. The paper is
organised as follows: 

In Section~1 we give some
conditions on a forcing in terms of  games that imply that the
forcing is $(T,Y,\cS)$-preserving. A special case of
$(T,Y,\cS)$-preserving is preserving the Souslinity of an
$\omega_1$-tree.

In Section~2
we show that for some tree-creature forcings from \cite{RoSh:470}
the player COM has a winning strategy in one of the games from Section~1.
Hence these forcings preserve Souslin trees.
Without the games, we show that some linear creature forcings from
\cite{RoSh:470} are $(T,Y,\cS)$-preserving. There are non-Cohen
preserving examples.

For the wider class of non-elementary proper forcings we
show in Section~3 that $\omega$-Cohen preserving for certain candidates implies
$(T,Y,\cS)$-preserving.

In Section~4 we give a less general but hopefully more easily readable
presentation of a result from \cite[Chapter~18, \S 3]{Sh:f}: If all
iterands in a countable support iteration are proper and
$(T,Y,\cS)$-preserving, then also the iteration is
$(T,Y,\cS)$-preserving. This is  a presentation of the so-called
Case A in which a division in forcings that add reals and those who
do not is not needed.

\section{A sufficient condition for $(T,Y,\cS)$-preserving}
\label{S1}
We introduce two games $\Game^\iota(\bP,p)$, $\iota = 1,2$,
that are  games about the completeness of the notion of
forcing $\bP$ above $p$. Similar games appear in \cite{RoSh:942,
RoSh:890, RoSh:888}. We let $\name{\bG_{\bP}} = \{ (\check{p},p) \such p \in
\bP\}$ be the standard name for a $\bP$-generic filter. If it is
clear which $\bP$ is meant we write just $\name{\bG}$.

\begin{definition}\label{1.1}
Let $\bP$ be a notion of forcing and $p \in \bP$. We define the
games $\Game^{\iota}(\bP,p)$, $\iota = 1,2$. The moves look
the same for both games, and only in the 
winning conditions they are different.

\begin{myrules1}
\item[(1)] The game $\Game^1(\bP,p)$ is played in $\omega$ rounds. In round
$n$, player COM chooses an $\ell_n\in \omega\setminus \{0 \}$ and a
sequence $\langle p_{n,\ell} \such \ell < \ell_n \rangle$ of
conditions $p_{n,\ell} \in \bP$ and then player INC plays $\langle
q_{n,\ell} \such n < \ell_n \rangle$ such that $p_{n,\ell} \leq
q_{n,\ell}$. After $\omega$ rounds, COM wins the game iff there is
$q \geq p$ such that for each $n$,
$$\{q_{n,\ell} \such \ell < \ell_n\} \mbox{  is predense above }q.$$

\item[(2)]The game  $\Game^2(\bP,p)$  is played in $\omega$ rounds that look exactly
like the rounds in $\Game^1(\bP,p)$. After $\omega$ rounds, COM wins
the game iff for every infinite $u \subseteq \omega$ there is $q_u
\geq p$ such that
$$q_u \Vdash (\exists^\infty n\in u)(\exists \ell < \ell_n)
(q_{n,\ell} \in \name{\bG}).$$
\end{myrules1}
\end{definition}

\begin{definition}\label{1.2}
For $\iota=1,2$, we say $\bP$ has property $\Pr^\iota$ and write $\Pr^\iota(\bP)$ iff
for every $p \in \bP$, in the game $\Game^\iota(\bP,p)$ the player
COM has a winning strategy.
\end{definition}

We fix a suffiently large regular cardinal $\chi$.
We write $H(\chi)$ for the set of sets of hereditary cardinality
less than $\chi$, and let $\cH(\chi) = (H(\chi),\in,<_\chi  ^*)$ with
a well-order $<_\chi^*$ on $H(\chi)$.

\begin{definition}\label{1.3}
Let $\alpha(\ast)$ be an uncountable ordinal. Let $\cS \subseteq
[\alpha(\ast)]^{\omega}$ be stationary, let $\iota=1,2    $, and
let $\bP$ be a forcing.  Then $\Pr^\iota_{\cS}(\bP)$ denotes the
following property: For every sufficiently large $\chi$ and every
countable $N \prec \cH(\chi)$ with $\bP \in N$,
and $N \cap \alpha(\ast) \in \cS$, for every 
  $p \in \bP\cap N$ player COM has a winning strategy
in the game $\Game^\iota(N,\bP,p)$. The game $\Game^\iota(N,\bP,p)$
is defined like the $\Game^\iota(\bP,p)$ except that we require that
every initial segment of a play is in $N$.
\end{definition}

$\Pr^\iota(\bP)$ implies $\Pr_{\cS}^\iota(\bP)$ for any $\cS$,
and $\Pr^1(\bP)$ implies $\Pr^2(\bP)$.

\begin{lemma}\label{1.4}
In all the games any winning strategy for COM can be modified by 
playing at each stage a larger number 
$\ell_n$ and stronger conditions, that is, the resulting function 
is a winning strategy for COM as well.
\end{lemma}

\begin{definition}\label{1.5}
 Let $\cS \subseteq
[\alpha(\ast)]^{\omega}$ be stationary. $\bP$ is $\cS$-proper if for any $ N \prec
\cH(\chi))$ such that $N \cap \alpha(\ast) \in \cS$, for any $p \in
\bP\cap N$ there is $q \geq p$ that is $(N, \bP)$-generic.
$q$ is $(N,\bP)$-generic means: 
For any $D \in N$, if $D$ is dense in $\bP$ then  $q \Vdash \name{\bG}
\cap D \neq \emptyset$.

\end{definition}

\begin{definition}\label{1.6}
\begin{myrules}
\item[(1)] A forcing $\bP$ is \emph{${}^\omega \omega$-bounding} if
for every $\bP$-name $\name{f}$  for a function from $\omega$ to $\omega$ and for any $p$, there are
$g \in {}^\omega \omega$ and $q\geq p$, $q\in \bP$, such that 
$q\Vdash \forall n \name{f}(n) \leq g(n)$.
\item[(2)] A forcing $\bP$ is  \emph{almost ${}^\omega \omega$-bounding} 
if for every name
$\name{f}$ for a function from $\omega$ to $\omega$, for any $A\subseteq\omega$ and any 
$p \in \bP$ there are $q\geq p$ and 
$g \in {}^\omega \omega$ such that $q \Vdash (\exists^\infty n \in A)
( \name{f(n)} \leq g(n))$.
\end{myrules}
\end{definition}

\begin{lemma}\label{1.7}
\begin{myrules1}
\item[(1)] If $\Pr_\cS^1(\bP)$, then $\bP$ is
$\cS$-proper, and $\bP$ is
${}^\omega \omega$-bounding.
\item[(2)] If $\Pr_\cS^2(\bP)$, then $\bP$ is
$\cS$-proper, and $\bP$ is
almost ${}^\omega \omega$-bounding.
\end{myrules1}
\end{lemma}

\begin{remark}\label{1.8}
The reverse implications do not hold: The NNR forcing from \cite[Ch.~IV]{Sh:f} 
is a counterexample to both, as Theorem~\ref{main} will
show.
\end{remark}

\proof We prove (2). Item (1) is proved similarly.
 Let $\name{f}$ be a $\bP$-name for a function from $\omega$
to $\omega$. Fix a winning strategy ${\bf st}$ for COM in
$\Game^2(\bP,p)$. Let $\bP, \name{f}, {\bf st} \in N \prec
\cH(\chi)$, $N \cap \alpha(\ast) \in \cS$, $p \in \bP\cap N$. Let
$\langle \name{\tau}_k \such k < \omega \rangle$ be a list of the
$\bP$-names in $N$ of ordinals. In  round  $n$, INC
plays such that for 
for every $\ell < \ell_n$, $q_{n,\ell}$ forces a
value to $\name{f}(i)$ for $i \leq n$ and a value to $\name{\tau}_i$
for $i <n$. 
Let $g(n)$ be the maximum of the values forced to $\name{f}(n)$ by
$q_{n,\ell}$, $\ell<\ell_n$.
Fix an infinite $u \subseteq \omega$ and let $q_u$ witness that COM wins. 
Then $q_u$ is $N$-generic: Let $\tau_k$ be a $\bP$-name in $N$ for an ordinal.
 Then $q_u$ forces that here are infinitely many $k' > k$, $k'\in u$ and 
$\ell\in \omega$ that that $q_{k',\ell} \in G$. This $q_{k',\ell}$ decides 
$(\tau_m)_{ m < k'}$ in $N$ and forces $f(k')$ to be some value less than $g(k')$.
\proofend

\begin{remark}\label{1.9}
Sacks forcing satisfies ${\rm Pr}^1$.
COM can fix
$\ell_n=2^n \cdot \ell_{n-1}$ and play the restrictions of 
$q_{n-1,i}$, $i <\ell_{n-1}$, 
to the members
of its $n$-th splitting
front as $p_{n,i}$, $i <\ell_n$.
\end{remark}

The following  versions of the games that work for all starting
points in a countable model simultaneously  are interesting for
themselves. However, \ref{1.10} and \ref{1.11} will not be
used in the sequel so that a reader who is mainly interested in
preserving Souslin trees can skip them.

\begin{definition}\label{1.10}
 Let $N\prec \cH(\chi)$.
We define a game $\Game^\iota(N,\bP)$: The moves are as in
$\Game^\iota(N,\bP,p)$. The winning conditions read
for $\iota =1$: For
every $p \in \bP\cap N$ there is a $q\geq p$
such that for all but finitely many $n$,
$\{q_{n,\ell}\such \ell < \ell_n\}$ is predense above $q$.
For $\iota =2$: For
every $p \in \bP\cap N$ and infinite $u$ there is a $q_u\geq p$
as in $\Game^2(\bP,p)$.
\end{definition}

\begin{lemma}\label{1.11} If $\Pr_\cS^\iota(\bP)$ and $N \cap \alpha(\ast) \in \cS$,
then COM has a winning strategy in $\Game^\iota(N,\bP)$.
\end{lemma}

\proof Let $N \cap \bP = \{ p_j \such j < \omega\}$. Let ${\bf
st}_j$ be a strategy for COM in $\Game^\iota(N,\bP,p_j)$.

\smallskip

Let in the $n$-th move strategy ${\bf
st}_j$ tell COM to choose $\overline{p_{j,n}} = \langle p_{j,n,\ell}
\such \ell < \ell_{j,n}\rangle$. Then COM moves in
$\Game^\iota(N,\bP)$ by letting $\ell_n = \sum_{j \leq n}
\ell_{j,n-j}$ and $\overline{p_n} = \overline{p_{0,n}} \concat \dots \concat
\overline{p_{j,n-j}} \concat \dots \concat \overline{p_{n,n-n}}$.
\proofend

Now we describe Souslin trees.

\begin{definition}\label{1.12}%(See \cite[Ch.~XVIII, 3.9]{Sh:f})
\begin{myrules}
\item[(1)] An $\omega_1$-tree is a tree of size
$\omega_1$ with at most countable levels and height $\omega_1$.
\item[(2)] Let $(T,<_T)$ be a tree. We let $T_{<_T s}= \{ t \in T\such t < s\}$,
and let $T_{\leq_T s}$, $T_{>_T s}$ be defined analogously.
\item[(3)] Let $(T,<_T)$ be a tree. Then we write 
$T_\alpha$ for $\{s\in T \such  (T_{<_T s},<_T) \cong \alpha\}$ and call
$T_\alpha$ the $\alpha$-th level of $T$.
\item[(4)]  
 Moreover, we require that the trees are normal
i.e., for every node $t$ on level $\alpha<\omega_1$ for 
every $\omega_1 >\beta> \alpha$
there are  $t'' \neq t' >_T t$ on level $\beta$.
\end{myrules}
\end{definition}

\begin{definition}\label{1.13} A Souslin tree is an $\omega_1$ tree 
that has no uncountable chains and no uncountable antichains.
\end{definition}

A notion of forcing $\bP$ preserves any Souslin tree if it preserves any
normal Souslin tree. This is seen as follows: Let $T$ be a Souslin tree.
We let $A = \{ t \in T \such T_{\geq_T t }$ is at most countable and 
$t$ is minimal with the property$\}$. Since $T$ is a Souslin tree, $A$ is at most countable. We let $T' = T \setminus A$. $T'$ is a Souslin tree
in $\bV^\bP$ iff $T$ is a Souslin tree in $\bV^\bP$.

Let $T$ be a normal $\omega_1$-tree. 
Let $b$ be a cofinal branch.
By normality, there are 
cofinally many $\alpha<\omega_1$ such that there are $t_\alpha \in b\cap T_\alpha$
and $t'_\alpha >_T t_\alpha$, $t'_\alpha  \not\in b$.
Then these $t'_\alpha$ form an antichain.
So $T$ is Souslin iff it does not have any uncountable
antichain.

\begin{definition}\label{1.14}
We conceive a normal $\omega_1$-tree $(T,<_T)$ 
without cofinal branches as a forcing notion. A stronger condition is
higher up in the tree.
For $\delta \in \omega_1$, we let $Y(\delta) \subseteq T_\delta$. 
Let $Y= \bigcup \{Y(\delta)\such \delta \in \omega_1\}$ and
reversely, given $Y \subseteq T$ we let $Y(\delta) = \{ t \in Y
\such t \in T_\delta\}$. Some of the $Y(\delta)$ may be empty. Let
$\cS\subseteq [\omega_1]^{\omega}$ be stationary.

 We say $T$ is {\em $(Y,\cS)$-proper} iff $Y
\subseteq T$ and  for every
sufficiently large $\chi$ for every countable $N \prec \cH(\chi)$
with $\{T,\cS\} \subset N$ and $N\cap \omega_1\in \cS$, $\delta =
N\cap \omega_1$, for $t \in Y(\delta)$, $T_{<_T t} :=\{s
\such s <_T t\}$ is $(N,T)$-generic.
\end{definition}

The definition gets stronger the larger $Y$ is. 
For characterising Souslin trees $Y$ is taken to be a union 
of stationarily many level of the tree.
Another application is gotten by taking $Y(\delta)\neq\emptyset$ for
stationarily many $\delta$.
The following known lemma,  which is \cite[Claim 3.9 B]{Sh:f},
characterises normal Souslin trees.

\begin{lemma}\label{1.15}
Let $(T,<_T)$ be a normal $\omega_1$-tree.
The following are equivalent:
\begin{myrules1}
\item[(1)]
 $T$ is Souslin.
\item[(2)]
  $T$
is $(Y,\cS)$-proper for every stationary $\cS \subseteq
[\omega_1]^{\omega}$ and for every $Y$ of the form $\bigcup_{\delta
\in W} T_\delta$, such that $W\subseteq \{ \sup(a) \such a \in
\cS\}$ stationary.
\item[(3)]
  $T$
is $(Y,\cS)$-proper for some stationary $\cS \subseteq
[\omega_1]^{\omega}$ and for some $Y$ of the form $\bigcup_{\delta
\in W} T_\delta$, such that $W\subseteq \{ \sup(a) \such a \in
\cS\}$ stationary.
\end{myrules1}
\end{lemma}

\proof (1) implies (2): Let $T$ be a Souslin tree and let $N \prec
\cH(\chi)$ with $T \in N$, $N \in \cS$. Let $\delta = N \cap
\omega_1 \in W$. We show that every node $t$ on level $\delta$ is
$(N,T)$-generic: Let $I \in N$ be dense in $T$. Now let in $N$, $I'
\subset I$ be a maximal antichain in $T$. $N \models \mbox{``}I'
\mbox{ is countable''}$, so $I' \subseteq N$.
 Now
$\{ s  \in T \such (\exists r \in I')(r \leq_T s)\}  \cap \{ s \in T
\such s <_T t\} \neq \emptyset$, since otherwise $I' \cup \{t\}$ is
an antichain, in contradiction to the fact that by $N \prec
\cH(\chi)$ the set $I'\in N$ is also a maximal antichain in $T$ in
the sense of $\cH(\chi)$ and in the sense of $\bV$.

\smallskip

(3) implies (1). We fix $\cS$ and $Y$ as in (3). We consider the
case that $A \subseteq T$ is an uncountable maximal antichain and
take $N \prec \cH(\chi)$ with $T, A \in N$, $N\in \cS$, $\delta =
N \cap \omega_1 \in W= \{ \delta \in \omega_1 \such Y(\delta)=
T_\delta\}$. Then $A$ is dense in $T$ in $N$. However, since $A$ is
uncountable, there is $t' \in A \setminus N$. Let $t = t' \rest
\delta \in T_\delta$. The node $t$ is incompatible with every $a \in
A \cap N$, so $t$ cannot lie above an $a \in A \cap N$, so $T_{<_T
t}$ is not $(N,T)$-generic. \proofend

\begin{definition}\label{1.16}
We say $\bP$ is {\em $(T,Y,\cS)$-preserving} iff the following
holds: Let $\cS \subseteq \omega_1$ be stationary. There is  $x \in
H(\chi)$, for every countable $N \prec \cH(\chi)$ with $\{x,Y,T,\bP,\cS\}
\subseteq N$ and $p \in \bP \cap N$: if $N \cap \omega_1 =
\delta$, $N \cap \omega_1 \in \cS$, and for every $t \in Y(\delta)
$, $\{ s \such s <_T t\}$ is
$(N,T)$-generic,  then there is $q \geq_{\bP} p$ such that $q$
is $(N,\bP)$-generic \nothing{Heike added $q$ is $(N,\bP)$-generic.
It saves us some work with $N_1 \in N_2$, $N_1 \prec N_2$ and so
on.}
 and
\[q \Vdash_{\bP} (\forall t \in Y(\delta)) (\{ s \such s <_T
t\} \mbox{ is } (N[\name{\bG_{\bP}}],T)\mbox{-generic}).
\]
\end{definition}

We remark that the quantifier 
``for every countable $N \prec \cH(\chi)$ with $\{x,Y,T,\bP,\cS\}
\subseteq N$ and $p \in \bP \cap N$''
can be weakened and that the particular
choice of $x \in H(\chi)$ is not essential, see \cite[Theorem~2.13]{abraham:handbook}

In Section~\ref{S4} we show that ``$T$ is $(Y,\cS)$-proper'' is
preserved by countable support iterations of proper iterands if each
iterand preserves it. Since we are mainly interested in countable support iterations
(because of the club principle), we can focus onto the question:
Which iterands preserve ``$T$ is $(Y,\cS)$-proper''?

A sufficient criterion
is given by $\Pr^2_\cS(\bP)$.

\begin{theorem}\label{main}%1.20
Assume $\alpha(\ast)=\omega_1$ and $\cS \subseteq \omega_1$ is
stationary.  Let $T$ be an $\omega_1$-tree and $Y \subseteq T$. If
$\Pr^2_\cS(\bP)$, then $\bP$ is $(T,Y,\cS)$-preserving.
\end{theorem}

\proof  Assume $N \prec H(\chi)$, $N \cap
\omega_1 \in \cS$,  $N \cap \omega_1 = \delta$, and 
$\bP \in N$, $p \in N\cap \bP$, and assume
 for every $t\in Y(\delta)$, $\{ s \such s <_T t\}$ is
$(N,\bP,p)$-generic. Let $x = {\bf st}$ for a winning strategy ${\bf st}$ for
player COM in $\Game^2(N,\bP,p)$.
We show that there is a $q$ as required
in the previous definition.

\smallskip

Let $Y = \{ t^\delta_k \such k < \gamma_\delta, \delta \in
W\}$ for suitable $\gamma_\delta\leq\omega$. Let $\{ \name{\cI}_n \such n \in \omega \}$ list the
$\bP$-names of open dense sets in the forcing $T$ that are in $N$
and let $\{\cJ_n
\such n \in \omega\}$ list the open dense sets in $\bP$ in $N$.  Now
we take a play $\langle
\langle \overline{p_n},\overline{q_n} \such n \in \omega
\rangle$ in which COM plays
according to ${\bf st}$. INC plays in every round $n$ in every 
$i < \ell_n$ the condition $q_{n,i}$ so strong that
$q_{n,i} \in \bigcap_{r<n}\cJ_r$ and such that for every $k < n$ there is
$t_{i,n,k} <_T t^\delta_k$ such that
\[
q_{n,i} \Vdash_{\bP} t_{i,n,k} \in \bigcap_{k'< n}\name{\cI}_{k'}.
\]
Why can INC play like this? Given $i <\ell_n$ and
a starting point $q'$, for $k<n$ 
he can strengthen $q_{n,i}\geq q'$ so that 
$q_{n,i} \Vdash_{\bP} t_{i,n,k} \in \bigcap_{k'< n}\name{\cI}_{k'}$
for a suitable $t_{i,n,k} <_T t^\delta_k$. Since $\{s
\such s<_T t_k^\delta\}$ is $(N,T)$-generic, there is such a
$t_{i,n,k} <_T t^\delta_k$, $t_{i,n,k} \in \cJ$.
Now he repeats this for each $k < n$.
Since $ \bigcap_{k'< n}\name{\cI}_{k'}$ is (forced by the
weakest conditions to be) open dense in the forcing $T$, the set
$\cJ = \{ s \in T \cap N \such q\not\Vdash_{\bP} s \not\in
\bigcap_{k'< n}\name{\cI}_{k'}\}$ is dense in $T$ in the ground
model (before forcing with $\bP$).

\smallskip

COM wins the play because he played according to the strategy. So
for every $u$, in particular for $u = \omega$,  there is $q_u \geq
p$ such that
\begin{equation}\label{1}
q_u \Vdash (\exists^\infty n \in u )(\exists \ell < \ell_n)
(q_{n,\ell} \in \name{\bG}_\bP).
\end{equation}
Let $k \in \omega$ and $q'\geq q_u$ be given.
Then there is $q{''}\geq q'$ and $n \geq k$ such that $q{''} \Vdash
n\in u$. 
So there is $i <\ell_n$, $q{''}\Vdash q_{n,i}\in \name{\bG}_\bP$ and hence
\begin{equation}\label{2}
q{''} \Vdash_{\bP}
t_{i,n,k} \in \bigcap_{k'< n}\name{\cI}_{k'}\wedge t_{i,n,k}\leq_T t^\delta_k.
\end{equation}
Now we unfreeze $k$ and combine the equations
\eqref{1} and \eqref{2} and thus get
\[
q_u \Vdash (\forall k < \omega)(T_{<_T t^\delta_k} \mbox{
is } (N[\name{\bG}_\bP],T)\mbox{-generic.})
\]
\relax From $q_{n,i}\in \bigcap_{r<n} \cJ_r$ we also get
that $q_u$ is $(N,\bP)$-generic. \proofend

\begin{corollary}\label{maincor}%1.20
If
$T$ is a Souslin tree, $\cS$ is stationary, and $\bP$ is a notion of
forcing with 
$\Pr^2_\cS(\bP)$, then $T$ is Souslin in $\bV^\bP$ .
\end{corollary}

\proof We let $Y=
\bigcup \{T_\delta \such \delta \in \cS\}$. By Lemma~\ref{1.15}
we have: $T$ is $(Y,\cS)$-proper iff it is a Souslin tree.
Now Theorem~\ref{main} shows the preservation of 
``$T$ is $(Y,\cS)$-proper''.
\proofend 

Historical remarks:
Our notion of $\cS$-properness this is called 
$(\{\cS\},\emptyset,\emptyset)$-properness in \cite[Def.~IV, 2.2.]{Sh:f}.
The notions of ${}^\omega \omega$-bounding and almost ${}^\omega \omega$- bounding
appeared in \cite{Sh:207}. A general study of preservation of these 
and related properties in iterations is in 
% (see \cite[Def.~VI, 2.8 A]{Sh:f}) (see \cite[Def.~VI, 3.5]{Sh:f})
\cite[Ch.~VI]{Sh:f}. A even more extensive study of preservation properties is
carried out in Chapter XVIII of \cite{Sh:f}. In \cite[XVIII 3.9 D]{Sh:f}
a variant of our definition of $(T,Y,\cS)$-preserving is mentioned.

\section{Many creature forcings $\bP$ 
preserve Souslin trees}\label{S2}

We begin this section with a proof that Miller forcing
hat $\Pr^2$. Then we look at other creature forcings.
We give a short self-contained introduction to creatures in general.
In the third subsection we consider tree creatures and give 
sufficient conditions for
$\Pr^1$ and
$\Pr^2$. In the fourth subsection we are concerned with linear creatures.
For these forcings we have not found strategies
in our games. However, some forcings of this kind preserve Souslin trees for
other reasons.

\subsection{A game on the Miller forcing}\label{S2.1}
Conditions in the Miller order are superperfect
trees $p \subseteq \omega^{<\omega}$ .
% also
%called rational perfect forcing (for a definition see
%\cite[Def.~7.3.43]{BJ}), has $\Pr^2$. 
A tree is called superperfect iff for any node $\eta\in p$ there is $
\varrho \trianglerighteq \eta$ that has infinitely many immediate
 successors in $p$. Here $\trianglerighteq$ denotes the end extension of
finite sequences.
Stronger conditions are perfect subtrees.

Miller forcing answers our question
about the consistency of the club principle together with the
existence of a Souslin tree and $\cov{\meager}=\aleph_1$ and
$\cof{\meager}=\aleph_2$, since it is well-known \cite{ncf3,
Vojtas87, Rothberger39} that in the Miller model $\cov{\meager} \leq
\gu = \aleph_1$ and $\cof{\meager} \geq \gd = \aleph_2$. The
``Miller model'' means any countable support iteration of Miller
forcing over a ground model of $\CH$.

The Miller conditions such that each splitting node has
infinitely many immediate successors are dense in the Miller 
order. From now on we work only with such conditions.
For $r \in \bP$ let $\rt(r)$ be the trunk, that is the
shortest $\eta$ such that $\suc_r(\eta)= \{\eta\concat n \such \eta
\concat n \in r \}$ is infinite.

\begin{theorem}\label{2.1}
Miller $\bP$ forcing has $\Pr^2(\bP)$.
\end{theorem}

\proof  We assume that
 all moves of both players in the game $\Game^2(\bP,p)$ below have infinite
splitting in each splitting node. 
Let $v \subseteq {}^{\omega>}\omega$. We let $\dcl(v) = \{ \eta
\rest k \such \eta \in v, k < \lg(\eta) \}$ be the downwards
closure of $v$. A set $v$ is a tree iff $v = \dcl(v)$.

 We describe a strategy ${\bf st}$
for COM in $\Game^2(\bP,p)$. On the side after the $n$-th
move COM chooses a finite set of nodes $v_{n}$ that are among the splitting nodes of INC's previously chosen conditions.
 COM play so that  the sequence
$\langle \overline{p}_n, \overline{q}_n, v_{n} \such n \in \omega\rangle$ has
the following properties:
\begin{myrules1}
\item[(0)] $\ell_0 = 1$, $p_{0,0} = p$, $q_{0,0} \geq p_{0,0}$,
$v_0 = \{ \tr(q_{0,0})\}$.
\nothing{\item[(1)] $v_{n} $ is  a finite subset of the splitting
 nodes of $p$, more precisely
  $v_{n} \setminus v_{n-1}$ is a subset of the union of the sets of splitting nodes
 of $q_{n,\ell}$, $\ell < \ell_n$.}
\item[(1)] For $n \geq 1$, given $v_{n-1}$, COM chooses $\ell_{n} =
|v_{n-1}|$ and for $\eta \in v_{n-1}$, $\eta = \rt(q_{n',\ell})$ for
some $n'<n$, $\ell < \ell_{n'}$ he lets
$$
m(\eta,n) = \min\{ k \such \eta \concat k \in q_{n',\ell} \setminus
\dcl(v_{n-1})\}.
$$ Since the $q_{n',\ell}$ is a Miller condition
 and each $\eta \in v_{n-1}$ is a splitting node of
$q_{n',\ell}$ for some $n'<n$ and $\ell < \ell_{n'}$ and $v_{n-1}$
is finite, for each $\eta \in v_{n-1}$, $m(\eta,n)$ is defined. Let
$\{\eta^n_\ell \such \ell < \ell_n\}$ enumerate $v_{n-1}$ and let
$\eta^n_\ell = \rt(q_{n',\ell'})$. Now COM chooses $p_{n,\ell}
=q_{n',\ell'}^{[\eta^n_\ell \concat m(\eta_\ell,n)]}$.
\item[(2)] INC plays $q_{n,\ell} \geq p_{n,\ell}$.
\item[(3)] Now COM chooses his new helper:
$v_{n} = v_{n-1} \cup\{ \rt(q_{n,\ell}) \such \ell < \ell_n\}$, and
the round is finished. Indeed $\ell_{n+1} = 2 \ell_n$ and $\ell_0 =
1$, but this is not important.
% Since $q_{n,\ell} \geq q_{n,\ell} \geq$
% $q_{n-1,\ell'}$, $\eta_\ell$ is a split node also of earlier conditions.
\end{myrules1}

\smallskip

The strategy ${\bf st}$ is a winning strategy for COM: Let $u
\subseteq \omega$ be infinite. By induction on $n \in u$ we choose
$s_n \subseteq v_{n} \setminus v_{n-1}$ . If $n= \min(u)$, then $s_n
\subseteq v_{n} \setminus v_{n-1}$ is a singleton. For $n
> \min(u)$, let 
\begin{multline*}
 s_n = s_{\max(u\cap n)} \cup\\
 \{\eta \in v_{n}
\setminus v_{n-1} \such \nu= \triangleleft\mbox{-}\max\{ \varrho \in
v_{n} \such \varrho \trianglelefteq \eta \} \in s_{\max(u\cap n)}\}.
\end{multline*}

 Lastly we let
\[
q_u = \{ \varrho \such (\exists n \in u)( \exists \eta \in s_n)
(\varrho \trianglelefteq \eta) \}.
\]
By definition, $q_u$ is a tree.
It is a Miller tree, since for every $n \in u$, for every $\eta \in
s_n$, $\eta$ is a splitting node in $q_u$ since  $\eta=
\rt(q_{n,\ell}) \in v_{n} \setminus v_{n-1}$. We show that for this
pair $(n,\ell)$, an infinite subset of $\suc_{q_{n,\ell}}(\eta)$  is
a subset of $\suc_{q_u}(\eta)$: For any $k
> n$ and there is  $\nu \trianglerighteq \eta \concat m(\eta,k)$, such that
$\nu  \in v_{k}
\setminus v_{k-1}$, by the choice of $\la v_k \such k \in \omega\ra$.
 For such a $\nu$ we
have $\max_{\triangleleft}\{\varrho \in v_k \such \varrho
\trianglelefteq \nu\} = \eta \in s_n \subseteq s_{\max(u \cap k)}$.
If $k \in u$, then $\nu\in s_k$.

\smallskip

 Moreover $q_u \Vdash_\bP (\exists^\infty n \in u)(
\exists \ell < \ell_n) (q_{n,\ell} \in \name{{\bf G}_{\bP}})$:
Suppose that not. Let $r \geq_{\bP} q_u$ be a Miller condition such
that $r\Vdash \forall n \in u (n \geq k \rightarrow (\forall \ell <
\ell_k)( q_{k,\ell} \not\in \name{{\bf G}_{\bP}})$. By strengthening
$r$, we may assume that $\rt(r) \in s_n$ for some $n \geq k$, so $r
\geq q_{n,\ell}$ for some $\ell < \ell_n$.  This is a contradiction.
\proofend

The properties $\Pr^1$ and $\Pr^2$ hold for
tree-creature forcings with the lim norm or the lim-sup norm. We
explain this now.

\subsection{Forcings with (tree) creatures}\label{S2.2}
Now we give more examples. In order to describe the relevant
properties we give a brief review of the definitions to forcings
with creatures. This concept is explained in the book
\cite{RoSh:470}, and it is divided into two main streams: one kind of
creature forcing is forcing with creatures such that the conditions
are written in an $\omega$-sequence like Blass--Shelah forcing
\cite{BsSh:242}. Another example of a forcing with
$\omega$-sequences of creatures is the (historically first) creature
forcing in \cite{Sh:207} that forces $\gb < \gs$. The other stream
is forcing with tree creatures. For historical reasons the first
kind is often just called ``creature forcing'' and the second kind
is called ``tree creature forcing''. In this subsection we give a
very short introduction to the main concepts.

%%BIS HIER
Let $\bH \colon \omega \to \cH(\omega_1)$ be a function such that
$(\forall n)(| \bH(n)| \geq 2 \wedge 0 \in \bH(n))$.

\begin{definition}\label{2.2}%(\cite[Def.~1.1.1]{RoSh:470})
Let $\chi$ be a regular cardinal. A triple $t = (\norm[t],
\val[t],\dis[t])$ is a {\em weak creature for $\bH$, $\chi$} if the
following hold:
\begin{myrules1}
\item[(a)] $\norm[t] \in \bR^{\geq 0} \cup \{ \infty\}$,
\item[(b)] $\val[t]$ is a non-empty subset of $$\Bigl\{(x,y)
\in \bigcup_{m_0<m_1<\omega} \prod_{i<m_0} \bH(i) \times
\prod_{i<m_1} \bH(i) \such x \triangleleft y \Bigr\},$$
\item[(c)] $\dis[t] \in \cH(\chi)$.
\end{myrules1}
The family of weak creatures for $\bH$ and $\chi$ is denoted by
$\WCR[\bH]$.
\end{definition}

We omit the parameter $\chi$ since in the following
 $\dis[t]$ is constant or empty.

\begin{definition}\label{2.3}%(\cite[Def.~1.1.3]{RoSh:470})
We say $\bH$ is {\em finitary} if $\bH(n)$ is finite for each $n$,
we say $\bH$ is of {\em countable character} if $\bH(n)$ is at most
countable for every $n$. We say $K \subseteq \WCR[\bH]$ is {\em
finitary} if $\bH$ is finitary and for every $t \in K$, $\val[t]$ is
finite.
\end{definition}

By our choice of  $\bH \colon \omega \to \cH(\omega_1)$ 
all the creatures in the following will be of countable character.

\begin{definition}\label{2.4}
%(\cite[Def.~1.1.4]{RoSh:470}) 
Let $K \subseteq \WCR[\bH]$.
\begin{myrules1}
\item[(1)] A function $\Sigma \colon [K]^{\leq \omega}
\to {\mathcal P}(K)$ is called a {\em sub-com\-position operation on
$K$} if the following holds:
\begin{myrules1}
\item[(a)]
(Transitivity:) if $\cS \in [K]^{\leq \omega}$ and for each $s \in
\cS$ we have $s \in \Sigma(\cS_s)$ for some $\cS_s\in \dom(\Sigma)$, 
then $\Sigma(\cS) \subseteq
\Sigma(\bigcup_{s \in \cS} \cS_s)$.
\item[(b)]  We write $\Sigma(r) $ for $\Sigma(\{r\})$.
$r \in \Sigma(r)$ for each $r \in K$ and $\Sigma(\emptyset)=
\emptyset$.
\end{myrules1}
\item[(2)] In the situation described above $(K,\Sigma)$ is called a
{\em weak creating pair}.
\end{myrules1}
\end{definition}

\begin{definition}\label{2.5}%(\cite[Def.~1.1.6]{RoSh:470})
Let $(K,\Sigma)$ be a weak creating pair for $\bH$.
\begin{myrules1}
\item[(1)] For a weak creature $t \in K$ we define its {\em basis} with
respect to $(K,\Sigma)$ as
$$\basis(t) = \Bigl\{ w \in \bigcup_{m<\omega} \prod_{i<m} \bH(i) \such
(\exists s \in \Sigma(t))(\exists u) (\langle w,u\rangle \in
\val[s])\Bigr\}.$$
\item[(2)]
For $w \in \bigcup_{m<\omega} \prod_{i<m} \bH(i)$ and $\cS \in
[K]^{\leq \omega}$ we define the set $\pos(w,\cS)$ of {\em possible
extensions of $w$ from the point of view of $\cS$} as
$$\pos^*(w,\cS)= \{ u \such \exists s \in \Sigma(\cS))(\langle
w,u\rangle \in \val[s])\},$$
\begin{equation*}
\begin{split}
\pos(w,\cS)= &\{ u \such \mbox{ there are $m\in \omega$ and
disjoint sets $\cS_i$ for
} i< m,
\bigcup_{i<m} \cS_i = \cS,\\
& \mbox{and a sequence } 0<\ell_1< \dots < \ell_{m-1} < \lg(u) \mbox{such that}\\
&u \rest \ell_1 \in \pos^*(w,\cS_0) \mbox{and} \\
&u\rest \ell_2 \in \pos^*(u\rest\ell_1,\cS_1), \dots, u \in
\pos^*(u\rest\ell_{m-1}, \cS_{m-1})\}.
\end{split}
\end{equation*}
\end{myrules1}
\end{definition}

\subsection{Tree creatures}
\relax From now on we specialise on tree creatures.
They have the special property that $\val[t]$ has just one root.

\begin{definition}\label{2.6}%(\cite[Def.~1.3.1]{RoSh:470})
\begin{myrules1}
\item[(1)] A {\em quasi tree}
 $(T,\triangleleft_T)$ is a set
of finite sequences, ordered by initial segment, and there is
a $\triangleleft_T$ smallest element $\rt(T)$, called the root of
$T$.
\item[(2)]
A quasi tree is called a {\em tree} if it is closed under initial
segments. If $T$ is a quasi tree we denote its closure under initial
segments by $\dcl(T)$. (This is the smallest tree containing $T$.)

\item[(3)]
We define the {\em set of immediate successors of $\eta$ in $T$},
the {\em restriction of $T$ to $\eta$}, the {\em splitting points of
$T$} and the {\em maximal points of $T$} by
$$
\suc_T(\eta) = \{ \nu \in T \such \eta \triangleleft_T \nu \wedge
\neg (\exists \rho \in T) (\eta \triangleleft_T \rho \triangleleft_T
\nu) \},$$
$$T^{[\eta]} = \{ \nu \in T \such \eta \trianglelefteq_T \nu \},$$
$$\splitt(T)= \{ \eta \in T \such |\suc_T(\eta)| \geq 2 \},$$
$$\max(T)= \{ \nu \in T \such
\neg (\exists \rho \in T)(\nu \triangleleft_T \rho\}.$$

\item[(4)]
The {\em $n$-th level} of $T$ is
$$T_n = \{ \eta \in T \such \eta \mbox{ has $n$
$\triangleleft_T$-predecessors}\}.$$ 

\item[(5)]
A {\em branch of $T$} is a maximal subset of $T$ that is linearly ordered by
$\triangleleft_T$. The {\em set of infinite
branches through $T$} is
\begin{equation*}
\lim(T)= \{ \eta \such \eta \mbox{ is an $\omega$-sequence and }
\wedge (\exists^\infty k)( \eta \rest k \in T)\}
\end{equation*}

A quasi tree is {\em well-founded} if there are no  infinite
branches through it.

\item[(6)]
A subset $F$ of a quasi tree $T$ is called a {\em front of $T$} if
every infinite branch of $T$ and every finite branch of $T$ passes
through this set, and the set consists of
$\triangleleft_T$-incomparable elements.
\end{myrules1}
\end{definition}

\begin{definition}\label{2.7}%(Compare to a part of \cite[Def.~1.3.3]{RoSh:470})
\begin{myrules1}
\item[(1)]
A weak creature $t \in \WCR[\bH]$ is a {\em tree creature} if
$\dom(\val[t])$ is a singleton $\{\eta\}$ and no two distinct
elements of $\rge(\val[t])$ are $\triangleleft$-comparable (so also
not compatible as finite partial functions
 since every $\eta \in \rge(\val[t])$ has as a domain
some $n \in \omega$). $\TCR[\bH]$ is the family of all tree
creatures for $\bH$.

\item[(2)] $\TCR_\eta[\bH]=\{ t \in \TCR[\bH] \such \dom(\val[t]) =
\{ \eta\}\}$.

\item[(3)]
A sub-composition operation $\Sigma$ on $K \subseteq \TCR[\bH]$ is a
{\em tree-composition} (and then $(K,\Sigma)$ is called a
tree-creating pair for $\bH$) if the following holds:

\begin{myrules1}
\item[($\ast$)] If $\cS \in [K]^{\leq \omega}$ and $\Sigma(\cS) \neq \emptyset$,
and $\cS = \{ s_\nu \such \nu \in \hat{T}\}$ for some well-founded
quasi tree $\hat{T} \subseteq \bigcup_{n<\omega} \prod_{i<n} \bH(i)$
and if for each finite sequence $\nu \in \hat{T}$,
$ s_\nu \in \TCR_\nu[\bH]$ and  for $\nu\in \hat{T} \setminus \max(\hat{T})$,
$ \rge(\val[s_\nu]) = \suc_{\hat{T}}(\nu)$
and if $t \in \Sigma(\{s_\nu \such \nu \in \hat{T}\})$ {\em then} $t \in
\TCR_{\rt(\hat{T})}[\bH]$ and
\nothing{Equation~%formerly \tag{$\spadesuit$}
needed to be changed in comparison to the book. Otherwise the part
on tree creatures in the book \cite{RoSh:470} all the time gives
away one layer of the tree in the subcomposition of a tree creature
condition.}
\begin{equation*}\label{change}
%formerly \tag{$\spadesuit$}
\rge(\val[t])\subseteq \bigcup\{ \rge(\val[s_\nu]) \such \nu \in
\max(\hat{T})\}.
\end{equation*}
\end{myrules1}
We write $\Sigma(s_\nu \such \nu \in \hat{T})$ instead of
$\Sigma(\{s_\nu \such \nu \in \hat{T} \})$.
 If $\hat{T}= \{
\rt(\hat{T}) \}$, $t = s_{\rt(\hat{T})}  \in
\TCR_{\rt(\hat{T})}[\bH]$ 
%and $\rge(\val[t]) = \max(\hat{T})$ 
then
we will write $\Sigma(t)$ instead of $\Sigma(s_\nu\such \nu \in
\hat{T})$.
\nothing{\item[(4)] A tree composition $\Sigma$ on $K$ is
{\em bounded} if for each $t \in \Sigma(s_\nu \such \nu \in
\hat{T})$ we have
$$\norm[t] \leq \max\{\norm[s_\nu] \such  (\exists \eta \in
\rge(\val[t]))(\nu \triangleleft \eta)\}.$$ }
\end{myrules1}
\end{definition}

So for a tree creating pair, 
if $t \in \TCR_\eta[\bH]$, then  $\basis(t) = \{\eta\}$ and
$\pos^\ast(\eta,\{t\})=\suc_t(\eta) = \rge(\val[t])$. We write only
$\pos^\ast(t)$ for $\pos^\ast(\eta,\{t\})$.

The next definition introduces requirements on the norms of
the creatures in a condition. We focus on the limsup condition and the lim condition. To speak in a uniform way about both variants, we introduce a
parameter $e$, and $e=0$ stands for the limsup case, and $e=1$ stands for 
the lim case.
\begin{definition}\label{2.8}%(Compare to \cite[Def.~1.3.3]{RoSh:470})
Let $(K,\Sigma)$ be a tree-creating pair for $\bH$, such that
there are
$t \in \TCR_\eta[\bH] \cap K$ of arbitrary high norm.
\begin{myrules1}
\item[(1)]
We define the forcing notion $\bQ^{\rm tree}_e(K,\Sigma)$ for $e =
0,1$  by letting $p = \langle t_\eta \such \eta \in  T \rangle\in
\bQ^{\rm tree}_e(K,\Sigma)$:
\begin{myrules1}
\item[(a)] $T\subseteq\bigcup_{n < \omega} \prod_{i<n} \bH(i)$
 is a non-empty quasi tree with $\max(T) = \emptyset$, and
\item[(b)]
$t_\eta \in \TCR_\eta[\bH] \cap K$ and $\pos^\ast(t_\eta) =
\suc_T(\eta)$, and
\item[(c)] in the  lim case ($e = 1$) we require for $\eta \in \lim(T)$,
$$
\lim \langle \norm[t_{\eta\rest k}] \such k < \omega, \eta\rest k \in T
\rangle=\infty.
$$ 
In the limsup case ($e = 0$) we require for $\eta \in
\lim(T)$ the sequence
$$
\limsup\langle \norm[t_{\eta\rest k}] \such
k < \omega, \eta\rest k \in T \rangle=\infty.
$$
 \end{myrules1}
 We define the
forcing order $\leq = \leq_{\bQ^{\rm tree}_e}$ by $\langle t^1_\eta
\such \eta \in T^1 \rangle \leq \langle t^2_\eta \such \eta \in
T^2\rangle$ iff $T^2 \subseteq T^1$ and for each $\eta \in T^2$
there is a quasi tree $\hat{T}_{0,\eta} \subseteq
(T^1)^{[\eta]}$ such that $\dcl(\hat{T}_{0,\eta})$ is well-founded and
 $t^2_\eta \in \Sigma(\{t^1_\nu \such \nu
\in \hat{T}_{0,\eta}\})$. 
If $t = \langle t_\eta \such \eta \in T
\rangle$ then we write $\rt(p) = \rt(T)$ and $T^p = T$ and $t^p_\eta
= t_\eta$ etc.
\item[(2)] If $p \in
\bQ_e^{\rm tree}(K,\Sigma)$ then we let $p^{[\eta]} = \langle
t^p_\nu \such \nu \in (T^p)^{[\eta]}\rangle$ for $\eta \in T^p$.
\end{myrules1}
\end{definition}

We write $\suc_p(\eta)$  for $\suc_{T^p}(\eta)$.
  
 The prominent real added by $\bQ_e^{\rm tree}(K,\Sigma)$, $ e = 0,1$, is
$\name{W}$ with 
$$
\Vdash_{\bQ_e^{\rm tree}(K,\Sigma)} \name{W} =
\bigcup\{\rt(p) \such p \in \name{\bG}_{\bQ_e^{\rm
tree}(K,\Sigma)}\}.
$$
Usually the conditions on the norm imply that 
$\name{W}$ is forced to be not in $\bV$.

\begin{definition}\label{2.9}%(\cite[Def.~1.3.7]{RoSh:470})
Let $(K,\Sigma)$ be a tree-creating pair $p \in \bQ^{\rm
tree}_e(K,\Sigma)$, $e=0,1$. A set $A \subseteq T^p$ is called {\em
an $e$-thick antichain} if it is an antichain in
$(T^p,\triangleleft)$ and for every condition $q \geq p$ the
intersection $A \cap \dcl(T^q)$ is not empty.
\end{definition}
\begin{proposition}\label{2.10}
%(Compare to part of \cite[Prop.~1.3.7]{RoSh:470})
\begin{myrules1}
\item[(1)] Let $e=0,1$. $\bQ_e^{\rm tree}$ is a partial order.
 Each $e$-thick antichain $A$ in $T^p$ gives  a maximal antichain
$\{p^{[\eta]} \such \eta\in A\}$ in  $\bQ_e^{\rm tree}$ above $p$.
Every front of $T^p$ is an $e$-thick antichain in $T^p$.
\item[(2)]
We define
\begin{equation*}
\begin{split} F_n^m(p) = \{ \eta \in T^p
\such &\norm[t^p_\eta]
> n \mbox{ and}\\
 & |\{\eta'\in T^p \such \eta' \triangleleft \eta \wedge
 \norm[t^p_{\eta'}]> n\}| = m \}.
\end{split}
\end{equation*}
Each $F^m_n(p)$ is a front of $T^p$ and an $e$-thick antichain of
$T^p$ for $e = 0,1$.
\item[(3)] If $K$ is finitary and   $\dcl(T^p)$ is well-founded, then
every front of
$T^p$  is finite.
\item[(4)] $p \leq p^{[\eta]} \in \bQ^{\rm tree}_e$ and $\rt(p^{[\eta]}) = \eta$.
\end{myrules1}
\end{proposition}

\proof  See  \cite[Prop.~1.3.7]{RoSh:470}.

Now the forcings notions with the normed trees let us define strengthenings of 
the forcing order $\leq$ that are natural candidates for Axiom A
(a definition can be found, e.g., in \cite[Def.~7.1.1]{BJ}).                         

\begin{definition}\label{2.11}%(Part of \cite[Def.~1.3.10]{RoSh:470})
Let $(K,\Sigma)$ be a tree-creating pair and let $p,q \in \bQ^{\rm
tree}_0(K,\Sigma)$.
\begin{myrules1}
\item[(1)] For the limsup case, we define $\leq_n^0$ for $n < \omega$ by\\
$p \leq_0^0 q$ if $p \leq q$ and $\rt(p)= \rt(q)$,\\
 $p\leq^0_{n+1}
q$ if $p \leq_0^0 q$ and if $\eta \in F_n^0(p)$ and $\nu \in T^p$
and $\nu \trianglelefteq \eta$ then $\nu \in T^q$ and
$t^q_\nu=t^p_\nu$.

\item[(2)] For the lim case we define $\leq_n^1$ for $n < \omega$ by\\
$p \leq_0^1 q$ if $p \leq q$ and $\rt(p)= \rt(q)$,\\
$p\leq^1_{n+1} q$ if $p \leq_0^1 q$ and if $\eta \in F_n^0(p)$ and
$\nu \in T^p$ and $\nu \trianglelefteq \eta$ then $\nu \in T^q$ and
$t^q_\nu=t^p_\nu$ and
$$
\{ (\eta,t_\eta^q) \such \eta \in T^q \wedge \norm[t_\eta^q] \leq n \}
\subseteq \{(\eta,t_\eta^p) \such \eta \in T^p\}.
$$
\end{myrules1}
\end{definition}

Note that $t_\nu^p=t_\nu^q$ means also that the immediate successors of 
$\nu$ in $p$ coincide with the immediate successors of $\nu$  in $q$. 

\begin{proposition}\label{2.12}
%(Part of \cite[Prop.~1.3.11]{RoSh:470}) 
Suppose that $e = 0,1$, $p_n
\in \bQ^{\rm tree}_e(K,\Sigma)$ and for $n \in \omega$, $p_n
\leq_{n+1}^e p_{n+1}$. Then the limit condition $p =\lim_{n\to
\omega} p_n$ is defined by $T^p = \bigcap_{n<\omega}T^{p_n}$  and
for $\eta \in \bigcap_{n<\omega} T^{p_{n}}$ we take the creature
$t^{p}_\eta=\bigcap_{n\in\omega} t_\eta^{p_n}$ (note that this is actually a 
finite intersection since the descending sequence
 $t_\eta^{p_n}$, $n\in \omega$ eventually becomes constant) into $p$. 
Then $p \in \bQ^{\rm tree}_e(K,\Sigma)$ and
$p \geq^e_{n+1} p_n$ for each $n$.
\end{proposition}

\begin{proposition}\label{2.13}%2.13
%(Compare to \cite[Prop.~2.3.1]{RoSh:470}) 
Let $p \in \bQ^{\rm
tree}_e(K,\Sigma)$, $n < \omega$, and let $A \subseteq T^p$ be an
antichain in $T^p$ such that $(\forall \eta \in A)( \exists \nu \in
F^0_n(p))(\nu\triangleleft \eta)$.
Assume that for each $\eta\in A$ we have a condition $q_\eta \in
\bQ^{\rm tree}_e(K,\Sigma)$ such that $p^{[\eta]} \leq_0^e q_\eta$
and
$$\mbox{if  $e=1$ then } (\forall \eta\in A)(\forall
\nu \in T^{q_\eta})
((\nu \trianglerighteq \eta \wedge \norm(t^{q_\eta}_\nu)\leq n) \rightarrow
t^{q_\eta}_\nu = t^p_\nu).
$$
 Then there exists $q \in \bQ^{\rm tree}_e(K,\Sigma)$ such that
$p \leq_{n+1}^e q$, $A \subseteq T^q$, $q^{[\eta]} = q_\eta$ for
$\eta \in A$ and if $\nu \in T^p$ is such that there is no $\eta\in
A$ with $\eta \trianglelefteq \nu$ then $\nu \in T^q$ and
$t^q_\nu=t^p_\nu$.
\end{proposition}

Since we repeatedly use the construction from
Proposition~\ref{2.13} in a re-ordered setting for the lim case,
we name it:

\begin{definition}\label{2.14}
We call the $q$ constructed from $p$, $A$ and $q_\eta$, $\eta \in
A$, as in the previous proposition: $$q = p \rest \{\nu \in T^p \such
\forall \eta\in A \nu \not\trianglerighteq_p \eta\} \concat
\sum_{\eta \in A} q_\eta.$$
When we use this expression we assume
that the conditions on $p$, $A$, $q_\eta,\eta\in A$, as given in the
proposition are fulfilled.
\end{definition}

\begin{definition}\label{2.15}%(\cite[Def.~2.3.4]{RoSh:470})
A tree-creating pair $(K,\Sigma)$ is {\em t-omittory} if for each
system $\langle s_\nu \such \nu \in \hat{T} \rangle$ such that
$\dcl(\hat{T})$ is a well-founded  tree and $\rt(s_\nu) =\nu$ and
$\pos^\ast(s_\nu) = \suc_{\hat{T}}(s_\nu)$ for $\nu \in \hat{T}
\setminus\max(\hat{T})$ and
for every $\nu_0 \in \hat{T}$ such that $\pos^\ast(s_{\nu_0}
)\subseteq \bigcup\{ \rge(\val[s_\nu]) \such \nu \in
\max(\hat{T})\}$ \nothing{parallel change to the change in the
definition of a tree composition in $(\spadesuit$) in
Def.~\ref{2.7}.} there is $s \in \Sigma(s_\nu \such \nu \in
\hat{T})$ such that
$$\norm[s] \geq \norm[s_{\nu_0}] -1 \mbox{ and }\pos^\ast(s) \subseteq
\pos^\ast(s_{\nu_0}).$$
\end{definition}

Note that t-omittoriness implies that the domain of $\Sigma$ 
contains $(s_\nu \such \nu \in
\hat{T})$ for  all well-founded subtrees $\hat{T}$.
A suitable equivalent formulation of Miller forcing is t-omittory.

Now there is an important construction we want to recall and use 
in the proofs of Lemma~\ref{2.17} and of Prop.~\ref{2.18} and
of Theorem~\ref{2.19}.

\begin{lemma}\label{2.16}
%(See \cite[Rem.~2.3.5]{RoSh:470}) 
Let $e = 0,1$  and let
$(K,\Sigma)$ be t-omittory. If $p \leq q$ then there is $r \in
\bQ^{\rm tree}_e(K,\Sigma)$ such that $p \leq^e_0 r$ and $\dcl(T^r)
\subseteq \dcl(T^q)$ and $t_\nu^r = t^q_\nu$ for  $\nu \in T^r
\setminus\{\rt(T^r)\}$ and $\rt(q) \in \dcl(T^r)$ and
$\norm(t^r_{\rt(r)}) \geq \norm(t^q_{\rt(q)})-1$.
\end{lemma}

\proof We let $\eta = \rt(q)$ and let $T^\ast$ be a well founded
quasi tree such that $(\forall \nu \in T^\ast) (\suc_{T^\ast}(\nu) =
\pos^*(t_\nu^p))$ and $\rt(T^\ast) = \eta$ and $t^q_\eta \in
\Sigma(t^p_\nu \such \nu \in T^\ast)$.
We let $T^- = \{\rt(p)\}
\cup \{\nu \in T^p \such \nu \triangleleft \eta\} \cup \{\eta\}$.
\nothing{{\bf Here was 
$\pos^*(t^q_\eta)$ in  \cite[Rem.~2.3.5]{RoSh:470}).}}
$T^-$ is a well-founded quasi tree and
we may apply t-omitting to $\langle t^p_\nu\such \nu
\triangleleft\eta \such \nu \in T^p\rangle \concat \langle
t^q_\eta\rangle$ and $\eta$. Thus we get $t^r_{\rt(p)} \in
\Sigma(\{t^p_\nu \such \nu \triangleleft \eta , \nu \in T^p\} \cup\{
t^q_\eta\})$ such that $\pos^*(t^r_{\rt(p)}) \subseteq
\pos^*(t^q_\eta)$ and $\norm(t^r_{\rt(r)}) \geq
\norm(t^q_{\rt(q)})-1$. Note that by transitivity of $\Sigma$,
$t^r_{\rt(p)} \in \Sigma(t^p_\nu \such \nu \in T^- \cup T^\ast)$.
For $\nu \in T^q$ such that $(\exists \nu' \in \pos^*(t_{\rt(p)}^r))
(\nu \trianglerighteq \nu')$ let $t^r_\nu = t^q_\nu$. \proofend

\begin{lemma}\label{2.17} %added July 26, 2010
Suppose that $(K,\Sigma)$ is a finitary t-omittory tree-creating
pair. Then $\bQ_1^{\rm tree}(K,\Sigma)$ is dense in $\bQ_0^{\rm
tree}(K,\Sigma)$.
\end{lemma}
\proof Given a $p \in \bQ^{\rm tree}_0(K,\Sigma)$, we repeatedly use
Lemma~\ref{2.16} to change it into a stronger condition in
$\bQ^{\rm tree}_1(K,\Sigma)$. \proofend

\begin{proposition}\label{2.18}
\begin{myrules1}
\item[(1)]
Suppose that $(K,\Sigma)$ is a finitary t-omittory tree-creating
pair. Then player COM has a winning strategy in 
$\Game^1(\bQ_1^{\rm tree}(K,\Sigma),p)$.
\item[(2)]
Suppose that $(K,\Sigma)$ is a finitary creating pair that is
t-omittory. Then player COM has as winning strategy in $\Game^1(\bQ_0^{\rm
tree}(K,\Sigma),p)$.
\end{myrules1}
\end{proposition}
\proof By Lemma~\ref{2.17} we need to prove only (2).
 We describe a strategy ${\bf st}$ for COM in
$\Game^1(\bP,p)$. The play will be $\langle \overline{p}_n, \overline{q}_n,
\such n \in \omega\rangle$. This time we let
$v_{n} = \{\rt(q_{n,\ell}) \such \ell < \ell_n\}$.
COM plays so that the play has the following properties:
\begin{myrules1}
\item[(0)] $\ell_0 = 1$, $p_{0,0} = p$, $q_{0,0} \geq p_{0,0}$,
$v_0 = \{ \rt(q_{0,0})\}$.
\item[(1)]
  $v_{n} \setminus v_{n-1}$ is a subset of some nodes $\eta_{n,\ell}$
 of $q_{n,\ell}$, $\ell < \ell_n$. In contrast to the 
proof of Theorem~\ref{2.1}, now $v_m \cap v_n =
 \emptyset$ and we only need to look at the $q_{n,\ell}$ from the
 previous round.
\item[(2)]  COM  lets for $\eta \in v_{n-1}$, $\eta =\eta^n_\ell=
\rt(q_{n-1,\ell})$ for some $\ell <\ell_{n-1}$,
$$
F(n,\eta) = \mbox{ a front in }\{\zeta \triangleright
\rt(q_{n-1,\ell}) \such \zeta \in T^{q_{n-1,\ell}} \setminus
\dcl(v_{n-1}), \norm(t^{q_{n-1,\ell}}_\zeta)>n\}.
$$
Since the $q_{n-1,\ell}$, $\ell < \ell_{n-1}$, are tree conditions
with pairwise incomparable roots and each $\eta \in v_{n}$ is a node
of $q'_{n-1,\ell}$ for some $\ell < \ell_{n-1}$ and $v_{n-1}$ is
finite, for each $\eta \in v_{n-1}$, $(\eta,n)$ is defined. Now COM
chooses for each $\eta \in v_{n-1}$, with $\eta= \rt(q_{n-1,\ell})$,
for each $\zeta \in F(\eta,n)$, and for each  $\rho\in
\rge(\val[t^{q_{n-1,{\ell}}}_\zeta])$,
$$p_{n,\eta,\rho} =(q_{n-1,{\ell'}})^{[\rho]}.$$
COM lets $\ell$ be the sum of these finite cardinalities of the
fronts $F(n,\eta)$ for all $\eta \in v_{n-1}$, of all  $\rho$'s for
all $\zeta \in F(n,\eta)$ and rearranges his move as
\begin{multline*}
\{p_{n,\ell} \such \ell <\ell_n\}= \{p_{n,\eta,\rho} \such \eta
=\eta^n_\ell\in v_{n-1}, \ell < \ell_{n-1},\zeta \in F(\eta,n),\\
 \rho \in
\rge(\val[t^{q_{n-1,{\ell}}}_\zeta])\}.
\end{multline*}
\item[(3)] INC plays $q_{n,\ell} \geq p_{n,\ell}$.
\end{myrules1}

Now we prove $\Pr^1(\bP)$.
We let
\[
q = \langle \varrho \such \varrho \in \rge(\val[t^{q_{n-1,\ell}}_\zeta])
\such  n\in \omega, \ell< \ell_{n-1}, \zeta \in F(\eta_\ell,n)\rangle.
\]
 By definition, $q$ is a quasi tree.
 It is a condition since we have
 by Prop.~\ref{2.13} that there is a sequence 
$\langle q_n\such n < \omega\rangle$
 such that
 $q_0 = q$,
 $q_{n+1} = q_n \rest \{\eta \such \forall \nu \triangleright \eta 
\nu\not\in \bigcup_{\ell < \ell_n}
 F_{n,\eta^n_\ell}\}\concat
 \sum_{\eta \in \bigcup_{\ell<\ell_n}F(n,\eta_\ell^n)} (q)_{n,\ell}^{[\zeta]}$.
We consider $q' \geq q$ and assume $q'$ forces that
$\{q_{n,\ell}\such \ell <\ell_n\}$ is not dense.  
By the tree omittoriness, we can find $q_n'\geq q', q_n$ as in Lemma~\ref{2.16}
such that  $\pos(t_\eta^{q'_n}) \subseteq
\pos(t_\eta^{q_n})$ for each $\eta\in T^{q'_n}$.
Then
 $q'_n \Vdash_\bP ( \exists \ell < \ell_n)
(q_{n,\ell} \in \name{{\bf G}_{\bP}})$ holds since 
 a subset of 
$\{\rt(q_{n,\ell}) \such \ell < \ell_n\}$ is a front in $T^{q'_n}$.
 \proofend

\begin{theorem}\label{2.19} 
Suppose that $\bigcup_{m<\omega} \bH(m)$ is countable
and $(K,\Sigma)$ is a tree creating pair for $\bH$ that is
t-omittory. Then the forcing
COM has a winning strategy in $\Game^2(\bQ_0^{\rm tree}(K,\Sigma),p)$.
\end{theorem}

\proof We show that there is a strategy for COM in $\Game^2(\bP,p)$.
This time he plays a tuple as a side move that is more complex than 
the one in Theorem~\ref{2.1}. 
We let 
 $F(n,\eta)$ and 
$\val[t^{q_{n-1,\ell}}_\zeta]$ be as in the proof of 
Proposition~\ref{2.18}. Now both of them are infinite. So in order to visit
them all, we organise the induction so that we at stage $n$ we visit
all the finitely many previous stages again and add just one node
$\zeta(n,\eta)$ of each previous $F(n',\eta)$, $n'< n$,
$\eta=\eta^n_\ell = \rt(q_{n',\ell'})\in v_n$ for a specific $n'<n$
and $\ell' < \ell_{n'}$,  and we add  one node $\rho(n,\eta)$ of
$\val[t^{q_{n',\ell'}}_{\zeta(n,\eta)}]$. Again $\ell_n = 2^n$. Both
kinds of tasks, the $F(n',\eta)$ and the
$\val[t^{q_{n,\ell'}}_{\zeta(n,\eta)}]$ will now appear as the
fronts $I_{{\bf x}, \eta}$ in the construction below.
\medskip

 We describe a strategy ${\bf st}$ for COM in
$\Game^2(\bP,p)$.

First, let $\langle \varrho^\ast_k \such k< \omega \rangle$ list
$\dcl(T^p)$.

We say ${\bf x}=(v_{\bf x}, \overline{p}_{\bf x}, \overline{I}_{\bf x})$ is an
{\em expanded state in $\Game^2(\bP,p)$} if ${\bf x}$ consists of
\begin{myrules1}
\item[(a)] $v = v_{\bf x}$ a finite, non-empty  set of splitting nodes of
$p$ with sufficiently high norm that has a root, $\rt(v_{\bf x})$, 
such that $\varrho \in v_{\bf
x} \rightarrow \rt(v_{\bf x}) \trianglelefteq \varrho$,
\item[(b)] a tuple of conditions
$\overline{p}_{\bf x} =\langle p_{{\bf x},\eta} \such \eta \in
v_{\bf x} \rangle$ such that $p \leq p_{{\bf x},\eta}$, $\eta =
\rt(p_{\bx,\eta})$ and $\norm(t^{p_{\bx,\eta}}_\eta) > 1$,
\item[(c)] a tuple $\overline{I}= \overline{I}_{\bf x}= \langle I_{{\bf x},\eta} \such \eta
\in v_{\bf x}\rangle$ of fronts such that $I_{{\bf x}, \eta}
\subseteq \dcl(T^{p_{{\bf x},\eta}})\setminus \{\eta\}$ is a front in 
$p_{{\bx},\eta}$,
it can be taken the direct successors of $\eta$ in $\dcl(T^{p_{{\bf x},\eta}})$,
\item[(d)] if $\eta \in v_{\bf x}$ and $\eta \triangleleft \nu \in I_{{\bf x},\eta}$
and $\eta \triangleleft \varrho \triangleleft \nu$ then $\varrho
\not\in v_{\bf x}$.
\end{myrules1}

For two expanded states ${\bf x}$, ${\bf y}$, we say ${\bf y} \in
\suc({ \bf x})$ if
\begin{myrules1}
\item[($\alpha$)] $v_{\bf x} \subseteq v_{\bf y}$ and $\rt(v_{\bf x})=
\rt(v_{\bf y})$, $(v_\by \setminus v_\bx)\cap \dcl(v_\bx) = \emptyset$,
 \item[($\beta$)] $\overline{p}_{\bf x} = \overline{p}_{\bf y} \rest
v_{\bf x}$,
\item[($\gamma$)] $\overline{I}_{\bf x} = \overline{I}_{\bf y} \rest
v_{\bf x}$,
\item[($\delta$)] if $\eta \in v_{\bf x} \setminus v_{\bf y}$ then
$\norm[t_\eta^{p_{{\bf y},\eta}}] \geq |v_{\bf x}|$,
%we write
%$t_\eta[q_{{\bf y},\eta}]$ instead of $t_\eta^{q_{{\bf y},\eta}}$,
\item[($\eps$)]
if $\eta \in v_{\bf x}$ and $k < \omega$ is minimal such that
$\varrho^\ast_k \in I_{{\bf x},\eta}$ and $(\neg \exists \varrho)
(\varrho^\ast_k \trianglelefteq \varrho \in v_{\bf x})$ then
$(\exists \varrho)( \varrho^\ast_k \trianglelefteq \varrho \in
v_{\bf y})$.
\end{myrules1}

COM chooses  on the side after the $n$-th move ${\bf x}_{n} \geq 
{\bf x}_{n-1}$ such the play
 $\langle \overline{p}_n, \overline{q}_n,
{\bf x}_{n} \such n \in \omega\rangle$ has the following properties:
\begin{myrules1}
\item[(0)]  $\ell_0 = 1$, $p_{0,0}= p$, $q_{0,0} \geq p$,
$\rt(q_{0,0}) = \eta$, now COM chooses $\nu \in T^{q_{0,0}}$ such
that $\norm[t^{q_{0,0}}_\nu] > 1$ and $v_{{\bf x}_0} = \{\nu\}$
 $p_{{\bf x}_0,\nu} =
q_{0,0}^{[\nu]}$, $I_{{\bf x}_0,\eta} =
\suc_{\dcl(T^{q_{0,0}})}(\nu)$,
\item[(1)]
In the $n$-th move COM first lets for $\eta \in v_{{\bf x}_{n-1}}$,
\[
k_{n,\eta} = \min\{ k \in \omega \such \varrho^\ast_k \in I_{{\bf
x}_{n-1},\eta} \wedge (\neg \exists \varrho )(\varrho^\ast_k
\triangleleft \rho \in v_{{\bf x}_{n-1}}\}.\]
COM makes the
move $\overline{p}_n= \langle p_{{\bf x}_{n-1},\eta}^{[\varrho^\ast_{k_{n,\eta}}]}
 \such \eta \in v_{{\bf
x}_{n-1}} \rangle$. Then INC makes moves $\langle q_\eta^{\ast,n}
\such \eta \in v_{{\bf x}_{n-1}}\rangle$ so $p_{{\bf x}_{n-1},\eta}
\leq q_\eta^{\ast,n}$ and $\varrho^\ast_{k_{n,\eta}} \trianglelefteq
\rt(q_\eta^{\ast,n})$.
\item[(2)] Now on the side COM chooses $\langle \nu_\eta^n
 \such \eta \in v_{{\bf x}_{n-1}}\rangle$ such that
($\nu_\eta^n \in T^{q^{\ast,n}_\eta}$ and
$\norm(t_{\nu_\eta^n}[q^{\ast,n}_\eta]) > |v_{{\bf
x}_{n-1}}|=2^{n-1}$ and $\varrho^{\ast,n}_{k_{n,\eta}}
\trianglelefteq v_\eta^n$).
\item[(3)]
COM defines ${\bf x}_n$ with the following properties: 
\[v_{{\bf
x}_n} = v_{{\bf x}_{n-1}} \cup \{ \nu_\eta^n \such \eta \in v_{{\bf
x}_{n-1}} \},\]
\[p_{{\bf x}_n,\nu} = \left\{
\begin{array}{ll}
p_{{\bf x}_{n-1},\nu} & \mbox{ if } \nu \in v_{{\bf
x}_{n-1}};\\
(q^{\ast,n}_\eta)^{[\nu]} & \mbox{ if } \eta \in v_{{\bf x}_{n-1}}
\wedge \nu = \nu^n_\eta,
\end{array} \right.
\]
\[I_{{\bf x}_n,\nu} = \left\{
\begin{array}{ll}
I_{{\bf x}_{n-1},\nu} & \mbox{ if } \nu \in v_{{\bf
x}_{n-1}};\\
\suc_{q^{\ast,n}_\eta}(\nu_\eta^n) & \mbox{ if } \eta \in v_{{\bf
x}_{n-1}} \wedge \nu = \nu^n_\eta.
\end{array} \right.
\]
Now the round is finished.
\end{myrules1}
\smallskip

Now we prove $\Pr^2(\bP)$. Let $\langle \overline{p}_n, \overline{q}_n, {\bf
x}_n \such n < \omega\rangle$ be a play in which COM uses ${\bf st}$
and let $u \subseteq \omega$ be infinite. We show how to define
$q_u$.

For $m_1 < m_2 < \omega$ we define a function $f_{m_1,m_2}$ as
follows
$$\dom(f_{m_1,m_2})= \{ \nu \such(\exists \eta \in v_{{\bf
x}_{m_1}})(\nu \in I_{{\bf x}_{m_1}, \eta} \wedge (\exists \rho \in
v_{{\bf x}_{m_2}})( \nu \trianglelefteq \rho)) \}$$ and for $\nu \in
\dom(f_{m_1,m_2})$ we let $f_{m_1,m_2}(\nu) \in v_{{\bf x}_{m_2}}$
be such that $\nu \triangleleft f(\nu)$  and $(\neg \exists
\rho)(f_{m_1,m_2}(\nu) \triangleleft \varrho \in v_{{\bf
x}_{m_2}})$, that is, $f_{m_1,m_2}(\nu)$ is $\triangleleft$-maximal.

Next we choose $w_{u,n}$ by induction on $n \in u$ such that
$w_{u,n} \subseteq v_{{\bf x}_{n}}$.

Case 1:  $n = \min(u)$. We let $\eta \in v_{{\bf x}_{n}}$ be
$\triangleleft$-maximal and let $w_{u,n}=\{\eta\}$.

 Case 2: $n \in u$, $n > \min(u)$. $m = \max(u \cap
n)$, $$ w_{u,n}= w_{u,m} \cup \{ f_{m,n}(\nu) \such (\exists \eta
\in w_{u,m})( \nu \in I_{{\bf x}_m,\eta} \wedge (\neg \exists
\varrho)( \nu \trianglelefteq \varrho \in v_{{\bf x}_m}))\}.$$
We 
define $q_u \in \bP$ by induction on $n$ such that 
 by $T^{q_u} \subseteq \dcl(\bigcup_{m \in u}
w_{u,m})$. For $\zeta\in T^{q_u}$ we let
$f_{m,n}(\nu) = \zeta$ and set  $\pos(t^{q_u}_\zeta)\subseteq
 \pos(t^{ p_{{\bf x}_n, \zeta}}_\zeta)$ and 
$t^{q_u}_\zeta=
t^{p_{{\bf x}_n, \zeta}}_\zeta$. 
Note that $\dcl(\bigcup_{m \in u} w_{u,m})$ is a tree
without maxima since $(\forall m < n)(m,n \in u \rightarrow(\forall
\eta \in w_{u,m})(\exists \nu \in I_{{\bf x}_m, \eta})(\exists
\varrho \in w_{u,n}) (\rho = f(\nu)))$. 
We show:
\begin{equation} \tag{$\odot$}q_u
\Vdash (\exists^\infty n \in u) \{ q^{\ast,n}_\eta \such \eta \in
w_{u,n}\} \mbox{ is predense}.\end{equation}

Since $w_{u,n} \subseteq v_{{\bf x}_n}$ from ($\odot$) we get $q_u
\Vdash_\bP (\exists^\infty n \in u)( \exists \ell < \ell_n)
(q_{n,\ell} \in \name{{\bf G}_{\bP}})$. Suppose that ($\odot$) is
false. Let $r \geq_{\bP} q_u$ be a condition such that $r\Vdash
(\forall n \in u )(n \geq k \rightarrow (\forall \eta \in w_{u,k})(
q_{\eta}^{\ast,k} \not\in \name{{\bf G}_{\bP}}))$. 
Now we use t-omittoriness. By strengthening
$r$ according to Lemma~\ref{2.16}, 
we may assume that $\rt(r) \in w_{u,n}$ for some $n \geq k$,
so, by our construction of $q_u$, $r \geq q^{\ast,n}_\eta$ for some
$\eta \in w_{u,n}$. This is a contradiction. \proofend

\begin{definition}\label{2.20}
%\cite[Def.~3.2.3]{RoSh:470}
Let $(K,\Sigma)$ be a tree-creating pair for $\bH$, $k < \omega$.
\begin{myrules1}
\item[(1)] A tree creature $t \in K$ is called {\em $k$-big} if $\norm[t]
> 1$ and for every function $h \colon\pos(t) \to k $ there is 
$s \in \Sigma(t)$ such that $h \rest
\pos(s)$ is constant an $\norm[s] \geq \norm[t] -1$.
\item[(2)] We say $(K, \Sigma$) is {\em $k$-big} if every $t \in K$
with $\norm[t]>1$ is $k$-big.
\end{myrules1}
\end{definition}

Although t-omittoriness does not literally imply bigness, it gives
an analogue of bigness if the function $h$ from 
Definition~\ref{2.20}(1) colours the second but highest level of a tree of
possibilities, since this level corresponds to the top level of the tree
$\hat{T}$ in the definition of t-omittory. 
So every tree-creature forcing construction performed with
bigness can also be done with t-omittoriness. 
This sheds some light on the the
conditions in \cite[Lemma 2.3.6 and Theorem 2.3.7]{RoSh:470}.

Remarks: The definitions are
taken from \cite[Chapters~1--3]{RoSh:470}.
 Lemma~\ref{2.19} is an analogue to the result
that for linear creatures  various kinds of limits of norms
give the same notion of forcing under a suitable condition on finitariness 
and omittoriness \cite[Prop.~2.1.3]{RoSh:470}.
Proposition~\ref{2.18} adds an intermediate step to the implication that
for finitary t-omittory pairs $(K,\Sigma)$ the forcing notion $\bQ_0^{\rm
tree}(K,\Sigma)$ is ${}^\omega\omega$-bounding
(see \cite[Conclusion~3.1.1]{RoSh:470}).
\nothing{ Now we can adapt \cite[Lemma~2.3.6 (2) for $e=1$]{RoSh:470} to add an
intermediate step to . That
conclusion   says that these forcings are
${}^\omega\omega$-bounding.
}
Theorem~\ref{2.19} is a strengthening of 
the implication:
If  $\bigcup_{m<\omega} \bH(m)$ is countable
and $(K,\Sigma)$ is a tree creating pair for $\bH$ that is
t-omittory, then the forcing
$\bQ_0^{\rm tree}(K,\Sigma)$ is almost ${}^\omega \omega$-bounding
(which is \cite[Theorem 4.3.9]{RoSh:470}).

\subsection{The case of linear creature forcings}\label{S2.4}

In this section we look at a second main kind of creature forcings,
namely forcings with $\omega$-sequences of creatures. This kind is
sometimes just called creature forcing, for historic reasons. The
best-known examples are Blass--Shelah forcing \cite{BsSh:242} and the
forcing from \cite{Sh:207}.

Blass--Shelah forcing \cite{BsSh:242} fulfils
the conditions of the next theorem. 
The assumptions resemble the assumptions made in 
Prop.~\ref{2.18} and Theorem~\ref{2.19}. 
Under these conditions, the various limit conditions on the divergence of 
norms coincide:
$\bQ^\ast_{{\rm w},\infty}(K,\Sigma)$ and
$\bQ^\ast_{\infty}(K,\Sigma)$ and $\bQ^\ast_{{\rm
s},\infty}(K,\Sigma)$ are equivalent forcings by
\cite[Prop.~2.1.3]{RoSh:470}.
We first recall these families of notions of forcing:

\begin{definition}%\cite[Def.~1.1.7]{Rosh:470}.
Suppose that $(K,\Sigma)$ is a weak creating pair for
$\bf H$ and ${\mathcal C}(\norm)$ is a property of $\omega$-sequences
of weak creatures from $K$ (i.e., ${\mathcal C}(\norm)$ is a subset of
$K^\omega$). We define the forcing notion $\bQ_{{\mathcal C}(\norm)}(K,\Sigma)$.
Conditions are pairs $(w,T)$ such that for some $k_0 < \omega$,
\begin{myrules}
\item[(a)] $w \in \prod_{i<k_0} {\bf H}(i)$.
\item[(b)] $T = \langle t_i \such i<\omega\rangle$ where
\begin{myrules}
\item[(i)] $t_i \in  K$
 for each $i$.
\item[(ii)] $w \in \basis(t_i)$ for some $i < \omega$, and for each
 finite set $I_0 \subseteq \omega$ and
$u \in \pos(w, \{t_i \such i \in I_0\})$  there
is $i \in \omega\setminus (\max(I_0)+1)$ such that $u \in \basis(t_i)$.
\end{myrules}
\item[(c)]
the sequence $\langle t_i \such i<\omega\rangle$ satisfies the conditions
${\mathcal C}(\norm)$.
\end{myrules}
The order is given by $(w_1,T^1)\leq (w_2, T^2)$ if and only if
for some disjoint sets ${\mathcal S}_0$, ${\mathcal S}_1$, $\ldots
\subseteq\omega$ we have
$w_2 \in \pos( w_1, \{t^1_\ell \such \ell \in {\mathcal S}_0\})$ and
 $t_i^2 \in \Sigma(\{t^1_\ell \such \ell \in {\mathcal S}_{i+1} \})$ for
each $i <\omega$
where $T^i = \langle t^i_i \such i<\omega\rangle$.

If $p = (w,T)$ we let $w^p = w$ and $T^p = T$ and if $T^p = \langle t_i \such
i< \omega \rangle$ then we let $t_i^p = t_i$. We may write $(
w,t_0,t_1\dots )$ instead of $(w,T)$ when $T = \langle t_i \such i
<\omega \rangle$.
\end{definition}

If $(K,\Sigma)$ is a weak creating pair and $\mathcal C(\norm)$
is a property of sequences of elements of $K$ then
$\bQ_{\mathcal C}(\norm)$ is a forcing notion.
Now we explain what properties $\mathcal C(\norm)$ are meant in (c)
of the previous definition.

\begin{definition} %\cite[1.1.10] 470
For a weak creature $t$ let us denote
$$m_{\rm dn}^t = \min\{\lg(u) \such u \in \dom(\val[t])\}
$$
We introduce the following basic properties of sequences of weak
creatures which may serve as ${\mathcal C}(\norm)$
\begin{myrules}
\item[($s\infty$)] A sequence $\la t_i \such i<\omega\ra$ satisfies
${\mathcal C}^{s\infty}(\norm)$ if and only if
$$\forall i<  \omega) (\norm(t_i) > \max\{i,m_{\rm dn}(t_i)\}).$$
\item[($\infty$)] A sequence $\la t_i \such i<\omega\ra$ satisfies
${\mathcal C}^{\infty}(\norm)$ if and only if
$$\lim_{i\to\infty}\norm(t_i) =\infty.$$
\item[($w\infty$)] A sequence $\la t_i \such i<\omega\ra$ satisfies
${\mathcal C}^{w\infty}(\norm)$ if and only if
$$\limsup_{i\to\infty}\norm(t_i) =\infty.$$

\end{myrules}
The forcing notions corresponding to the above properties for a
weak creating pair $(K,\Sigma)$ will be denoted by $\bQ_{s\infty}(K,\Sigma)$,
$\bQ_{\infty}(K,\Sigma)$, $\bQ_{w\infty}(K,\Sigma)$.
\end{definition}

Adding more properties
to a weak creature gives an ${\bf H}$-creature:

\begin{definition} %\cite[2.1.2}{RoSh:470}
Let $t$ be a weak creature for ${\bf H}$.
\begin{myrules}
\item[1.] If there is $m<\omega$ such that
$ \forall \la u , v\ra \in \val[t]$, $\lg(u) =m$, then this unique
$m$ is called $m_{\rm dn}^t$.
\item[2.] If there is $m<\omega$ such that
$ \forall \la u , v\ra \in \val[t]$, $\lg(v) =m$, then this unique
$m$ is called $m_{\rm up}^t$.
\item[3.] If both $m_{\rm dn}^t$ and $m_{\rm up}^t$ are defined then
$t$ is called an
$(m_{\rm dn}^t,m_{\rm up}^t)$-creature of just a creature.
\item[4.] ${\rm CR}_{m_{\rm dn}^t,m_{\rm up}^t}[{\bf H}] =
\{ t \in {\rm WCR}[{\bf H}]\such m_{\rm dn}^t = m_{\rm dn}, m_{\rm up}^t =
m_{\rm up}\}$. The set ${\rm CR}[{\bf H}] = \bigcup_{m_{\rm dn}<m_{\rm up} < \omega}
 \rm CR_{m_{\rm dn}^t,m_{\rm up}^t}[{\bf H}]$ is called the set of ${\bf H}$-creatures.
\end{myrules}
\end{definition}

\begin{definition}\label{2.24}%\cite[1.2.2]{RoSh:470} 
Suppose that
$K \subseteq CR[\bf H]$ and $\Sigma$ is a subcomposition operation on $K$.
 We say that $\Sigma$ is a \emph{composition on $K$} and we say 
$(K,\Sigma)$ is
a \emph{creating pair for $\bf H$} if
\begin{myrules}
\item[(1)]if ${\mathcal S} \in [K]^{\leq \omega}$ and
$\Sigma({\mathcal S}) \neq \emptyset$ then ${\mathcal S}$ is finite and
for some enumeration $S=\{t_0,\dots, t_{m-1}\}$ we have
$m_{\rm up}^{t_i} = m_{\rm dn}^{t_{i+1}}$
for $i<m-1$, and
\item[(2)] for each $s \in \Sigma(t_0,\dots, t_{m-1})$ we have $m_{\rm dn}^s=  m_{\rm dn}^{t_0}$
and $m_{\rm up}^s = m_{\rm up}^{t_{m-1}}$.
\end{myrules}
\end{definition}

\begin{definition}\label{2.25}%\cite[1.2.6]{RoSh:470}
Let $(K,\Sigma)$ be a creating pair and ${\mathcal C}(\norm)$ be a property
of $\omega$-sequences of creatures. The forcing notion
$\bQ^\ast_{{\mathcal C}(\norm)}(K,\Sigma)$ is a suborder of
$\bQ_{{\mathcal  C}(\norm)}(K,\Sigma)$ consisting of these conditions
$(w,t_0,t_1,\dots)$ for which additionally
$\forall i \in \omega$, $m_{\rm up}^{t_i} = m_{\rm dn}^{t_{i+1}}$.
\end{definition}

\begin{definition}\label{2.26}%\cite[Def.~2.1.1]{RoSh:470})
Let $(K,\Sigma)$ be a weak creating pair for ${\bf H} $.
\begin{myrules}
\item[1.] For $t \in K$, $m_0 \leq m_{\rm dn}^t$, $m_{\rm up}^t \leq m_1$
we define the creature $s = t \reri [m_0,m_1)$ by
\begin{eqnarray*}
\norm[s] &=& \norm[t],\\
\val[s] &=& \{\langle w,u\rangle \in \prod_{i<m_0} {\bf H}(i) \times
\prod_{i<m_1}{\bf H}(i) \such \langle v \rest m_{\rm dn}^t,
u \rest m_{\rm up}^t\rangle \in \val[t] \wedge\\
&& w \triangleleft u \wedge (\forall i \in [m_0,m_{\rm dn}^t) \cup
[m_{\rm up}^t,m_1)) (u(i)=0)\}.
\end{eqnarray*}
Note that $t \reri [m_0,m_1)$ is well-defined only if
$ \val[s] \neq\emptyset$  and then $m_{\rm dn}^s = m_0$ and $m_{\rm up}^s = m_1$.

\item[2.] The creating pair $(K,\Sigma)$ is {\em omittory} if it has the following properties:
\begin{myrules}
\item[(o$_1$)] If $t \in K$ and $u \in \basis(t)$ then $u
\concat 0_{[m_{\rm dn}^t,m_{\rm up}^t)} \in \pos(u,t)$ but there is $v
\in \pos(u,t)$ such that
$v\rest [m_{\rm dn}^t,m_{\rm up}^t) \neq   0_{[m_{\rm dn}^t,m_{\rm up}^t)}$.
\item[(o$_2$)] For every $(t_0, \dots, t_{n-1})$
sequence of $(K,\Sigma)$-creatures, if for every $i <n$,
 $m_{\rm dn}^{t_{i+1}}= m_{\rm up}^{t_i}$ then for 
every $i < n$,
$t_i \reri [m_{\rm dn}^{t_0}, m_{\rm up}^{t_{n-1}}) \in \Sigma(t_0,\dots ,t_{n-1})$.
\item[(o$_3$)] If $t, t\reri [m_0,m_1) \in K$ then for every $u \in
\basis(t \reri [m_0,m_1))$ and $v\in \pos(u, t\reri [m_0,m_1))$  we have
$$v(n) \neq 0 \wedge n \in [\lg(u),\lg(v)) \rightarrow
n \in [m_{\rm dn}^t,m_{\rm up}^t).$$
Note that (o$_1$) implies that in the cases relevant for (o$_2$) 
the creature $t
\reri [m_{\rm dn}^{t_0}, m_{\rm up}^{t_{n-1}})$ is well defined.
\end{myrules}
\end{myrules}
\end{definition}

\begin{definition}
%\cite[Def.~2.2.5]{ 470}
An omittory creating pair $(K,\Sigma)$ is {\em omittory-big} if for every
$ k<\omega$ there is $m<\omega$ such that if $t\in K$, $\norm(t) >m$,
$u \in \basis (t)$, $c\colon \pos(u,t) \to \{0,1\}$ then there is
$s \in \Sigma(t)$
such that $\norm(s) \geq k$ and $c\rest \pos(u,s) \setminus
\{0_{[m_{\rm dn}^t,m_{\rm up}^t)}\}$ is constant. We call $m$ an \emph{omittory
bigness witness for $k$.}
\end{definition}

\begin{definition}
%\cite[Def.~1.1.3 (2)]{RoSh:470}
$(K,\Sigma)$ is finitary, that means every $\fc_i^p$ has a
finite range and $\Sigma(\cS) \neq \emptyset$ only for finite
subsets $\cS\subseteq K$ and also $\Sigma(\fc_0, \dots, \fc_{n-1})$
is finite.
\end{definition}

If $p = (\eta^p, \bc_0^p,bc_1^p\dots)$ is a condition and $n \in \omega$,
and $\nu \in \pos(\eta^p,\bc_0^p,\dots,c_{n-1}^p)$ then we let 
$p^{[\nu]} = (\nu, \bc_n^p,\bc_{n+1}^p,\dots)$.

\begin{theorem}\label{2.29}
Assume $\bP= \bQ^\ast_{{\rm w}\infty}(K,\Sigma)$
% \cite[Def.~1.1.6,Def.~ 1.1.7, Def.~1.1.10]{RoSh:470}
 is finitary
 and omittory
%\cite[Def.~2.1.1]{RoSh:470}
and is omittory-big.  %\cite[Def.~2.1.5]{RoSh:470}
%\nothing{I changed big into omittory-big}
Then $\bQ$ is
$(T,Y,\cS)$-preserving.
\end{theorem}

\proof Assume that $\chi \geq 2^{2^\omega}$ and $N\prec \cH(\chi)$
is countable  and that $\bP \in N$, $p \in N \cap \bP$, $T, \cS, Y
\in N$. Let $\delta = N \cap \omega_1$ and $N \cap \omega_1
\in \cS$. Let $T$ be $(Y,\cS)$-proper. Assume that $p =
(\eta^p,\bc_0^p,\bc_1^p, \dots)$.

 We show that there is
$q \geq p$ that is $(N,\bP)$-generic and such that for every $t \in
Y(\delta)$, $$q \Vdash T_{<_T t} \mbox{ is
$(N[\name{\bG_\bP}],T)$-generic}.$$

Now we use the Axiom A structure:
We enumerate all the  $\bP$-names $\name{\cI} \in N$
for dense sets in $T$  as $\{\name{\cI}_n \such n < \omega \}$, all
the $\cJ \in N$ that are dense in $\bP$ as $\{\cJ_n \such n <
\omega\}$ and all the $t \in Y(\delta)$ a $\{t_n \such n <
\omega\}$.
 We
choose $p_n$ by induction on $n \in [n_*,\omega)$ such that
\begin{myrules1}
\item[(a)] $p_n \in \bP \cap N$,
\item[(b)] $p_n \leq p_{n+1}$,
\item[(c)] $p_{n^*} = p$,
\item[(d)] for some countable $\cJ^\ast_n \subseteq \cJ_n$
$\cJ^\ast_n$ is predense above $p_{n+1}$,
\item[(e)] if $k > n \geq n_\ast$ then $\norm[\bc_k^{p_n}] \geq n$,
\item[(f)] $(\eta^{p_n},\bc_0^{p_n},\bc_1^{p_n},\dots,
\bc_{n-1}^{p_n})
=(\eta^{p_{n+1}},\bc_0^{p_{n+1}},\bc_1^{p_{n+1}},\dots,
\bc_{n-1}^{p_{n+1}})$,
\item[(g)] if $\nu \in \pos(\eta^{p_n},\bc_0^{p_n},\bc_1^{p_n},\dots,
\bc_{n-1}^{p_n})$ and there are $q \geq p_n$ and  $s \in T$
satisfying $(\ast)_{n,\nu,p_n,q,s}$ below, then $(\ast)_{n,\nu,p_n,
p_{n+1}^{[\nu]}, s}$.
\end{myrules1}
Here we use
\begin{equation}\tag*{$(\ast)_{n,\nu,p,q,s}$} \begin{split}
& \nu \in
\pos(\eta^{p},\bc_0^{p},\bc_1^{p},\dots, \bc_{n-1}^{p}) \wedge\\
& \eta^q = \nu \wedge \\
& q \geq (\nu,\bc_n^{p}, \dots) \geq p \wedge \\
& s <_T t_n \wedge \\
& q \Vdash s \in \bigcap_{k<n}\name{\cI_k}.
\end{split}
\end{equation}
There is no problem in carrying this induction as $\bP$ is finitary
and omittory.

In the end we let $$p_\omega = \lim_{n \to \omega} p_n = (\eta^p,
\bc_0^{p_{n_\ast}}, \dots,\bc_{n_\ast-1}^{p_{n_\ast}},
\bc_{n_\ast}^{p_{n_\ast+1}},\bc_{n_\ast+1}^{p_{n_\ast+2}},\dots).$$
By (f), $p_\omega$ is $(N,\bP)$-generic. Now we strengthen
$p_\omega$ once more to get a condition $p_{\omega+1} \geq p_\omega$
that forces that every $t \in Y(\delta)$ is
$(N[\bG_\bP],T)$-generic. This strengthening is carried out as follows:

Now for $n < \omega$ we let $C_n = \pos(\eta^{p_\omega},
\bc_0^{p_\omega}, \dots , \bc_{n-1}^{p_\omega})$. So $C =
\bigcup_{n<\omega} C_n$ is a tree. We colour this tree in two
colours: $\bc\colon C \to \{\mbox{yes}, \mbox{``no''}\}$ for $\nu \in
C$, $\bc(\nu) = \mbox{``yes''}$, iff for some $s<_T t$,
$(\ast)_{n, \nu,p_{\omega},p_{\omega}[\nu],s}$ and no otherwise.
If $n \leq n_1 < n_2$ and $\nu_i
\in C_{n_i}$ and $\nu_1 \triangleleft \nu_2$ and $\bc(\nu_1) =
\mbox{yes}$, then $\bc(\nu_2)= \mbox{yes}$, since $p_{\omega}^{[\nu_1]}
\leq p_{\omega}^{[\nu_2]}$.

Now by \cite[Theorem~2.2.6]{RoSh:470} we have the following
consequence of omittory-big: There is $p_{\omega+1} \geq_0 p_\omega$
such that \nothing{$(\eta^{p_\omega}, \bc_0^{p_{\omega}},
\dots,\bc_{n_\ast-1}^{p_{\omega}})= (\eta^{p_{\omega+1}},
\bc_0^{p_{\omega+1}}, \dots,\bc_{n_\ast-1}^{p_{\omega+1}})$} the
following holds: If $\nu_i \in C_{n_i}$, $n_1 < n_2$, and $\nu_i \in
\{\rt(q) \such p_{\omega+1} \leq q\}$ and $\nu_1 \triangleleft
\nu_2$ then $\bc(\nu_1)= \bc(\nu_2)$.

We check that the uniform colour is ``yes''. Suppose for a
contradiction that $(\forall \nu_\ast \triangleright
\eta^{p_{\omega+1}})( \bc(\nu_\ast) = \mbox{no})$. We let
$p_{\omega+2} = (\nu_\ast,
\bc_{m_\ast}^{p_{\omega+1}},\bc_{m_\ast+1}^{p_{\omega+1}},\dots ) \geq
p_{\omega+1}$   for a suitable $m_\ast$ with $m_ \ast \in C$. So
there are $s < t_{m_\ast}$ and $q \geq p_{\omega+2}$ with  $q \Vdash
s \in \bigcap_{k<m_\ast} \cI_k$.
 As $\cI_k$, $k < m_\ast$, are dense subsets of $(T,<_T)$
that have names in $N$ there is such a pair $(s,q)$. Now
$\bc(\rt(q)) = \mbox{yes}$. So the uniform colour cannot be ``no''.
\proofend

We recall Blass--Shelah forcing in order to see that it fulfils 
the conditions of the previous theorem.

\begin{definition}\label{2.30}%\cite[Def.~2.4.2]{RoSh:470}.
We define a depth function on $\{ A \subseteq
[\omega]^{<\omega}\such 2 \leq |A| <\omega\}$ as follows:
\begin{equation*}\begin{split}
{\rm dp}(A) \geq 0, & \mbox{ always,}\\
 {\rm dp}(A) \geq 1, & \mbox{ if $A\neq \emptyset$,}\\
{\rm dp(A)} \geq n+2, & \mbox{ if for every set $X \subseteq \omega $ one
of the following conditions holds}\\
& {\rm dp}(\{a \in A\such a \subseteq X\}) \geq n+1,\\
&\mbox{ or } {\rm dp}(\{ a \in A \such a\subseteq \omega \setminus X \})
\geq n+1.
\end{split}
\end{equation*}
\end{definition}

\begin{definition}\label{2.31}
Blass--Shelah forcing is $\bQ^*_{s \infty}(K,\Sigma)$ with the
following creating pair $(K,\Sigma)$: We let $H(m)=2$ for $m \in
\omega$. A creature $t \in {\rm CR}[H]$ is in $K$ if $m_{\rm dn}^t
+2 < m_{\rm up}^t$ and there is a sequence $\langle A_u^t \such u
\in \prod_{i< m_{\rm dn}^t} H(i) \rangle$ such that for every $u \in
\prod_{i<m_{\rm dn}^t} H(i)$ the following holds:
\begin{myrules}
\item[($\alpha$)] $A_u^t$ is a non-empty family of subsets of
$[m_{\rm dn}^t,m_{\rm up}^t)$ such that each member of $A_u^t$ has
at least 2 elements,
\item[($\beta$)] $\langle u,v\rangle \in \val[t]$ iff $u \triangleleft v$ and
$\{ i \in [m_{\rm dn}^t,m_{\rm up}^t ) \such v(i) = 1 \} \in A_u^t
\cup \{\emptyset\}$
\item[($\gamma$)] $\norm(t) = \min \{\log_2({\rm dp}(A_u^t)) \such u \in \prod_{i<
m_{\rm dn}^t} H(i)\}$.
\end{myrules}
Suppose $t_0, \dots, t_n$ in $K$ are such that $m_{\rm dn}^{t_{i+1}}
= m_{\rm up}^{t_i}$ for $i < n$. Then $s \in \Sigma(t_0, \dots, t_{n}
)$ iff $s \in K$ and $m_{\rm dn}^s = m_{\rm dn}^{t_0}$ and $m_{\rm
up}^s = m_{\rm up}^{t_{n}}$ and for every $\langle u,v \rangle \in
\val(s)$ for every $i \leq n$, $\langle v \rest m_{\rm dn}^{t_i}, v
\rest m_{\rm up}^{t_i}) \in \val[t_i]$.
\end{definition}

Blass-Shelah forcing is finitary, omittory and omittory-big.
So by Theorem~\ref{2.29} it is $(T,Y,\cS)$-preserving.

There is a parallel result without the property ``omittory'' but
with strong enough bigness and halving.
\begin{theorem}\label{2.32}
Assume that $\bP= \bQ_{\rm w}(K,\Sigma)$ is creature forcing with
the following properties:
\begin{myrules1}
\item[(a)] $p \in \bP$ has the form $(f,\fc_0,\fc_1, \dots) =
(w^p,\fc_0^p,\fc_1^p,\dots)$ with $\liminf \la \norm(\fc_n) \such n
< \omega \ra = \infty$.
\item[(b)] $(K,\Sigma)$ is finitary.
\item[(c)] For some sufficiently fast increasing sequence $\overline{k}=\la k_i
\such i < \omega \ra$ we have the following strong versions of bigness and
halving: First, we assume that there is a function $i \colon K \to
\omega$  such that
\begin{myrules1}
\item[--] $\fc \in\Sigma(\fc_0,\dots, \fc_{n-1}) \rightarrow
i(\fc)\leq \max\{i(\fc_j) \such j < n\}$,
\item[--] in every condition $p$,  $i(\fc^p_0) < i(\fc^p_1) < i(\fc^p_2)\dots$,
\item[--]
for every $\fc\in K$ and $n$ we have $ |\{(f,\fc^p_0, \dots,
\fc^p_{n-1}) \such p \in \bP \wedge \fc_{n-1}^p = \fc\}| \leq
k_{i(\fc)}$.
\end{myrules1}
Now for such a sequence $\overline{k}$ and function $i$ 
we require:
\begin{myrules1}
\item[($\alpha$)]
$\norm(\fc) \in \{\frac{m}{n} \such n \leq k_{i(\fc)}, m \leq
k_{i(\fc)}!!\}$,
\item[($\beta$)]
for every $p \in \bP$, $n \in \omega$,
 $|\pos(f^p,\fc_0^p, \dots,\fc^p_{n-1})| \ll
k_{i(\fc^p_n)}$,
\item[($\gamma$)]
(bigness) for every $p \in \bP$, $n \in \omega$, $d \colon
\pos(f^p,\fc^p_0, \dots, \fc^p_{n}) \to k_{i(\fc^p_n)}$  there is
$\fc\in \Sigma(\fc^p_n)$ such that $\norm(\fc) \geq \norm(\fc^p_n) -
\frac{1}{k_{i(\fc^p_n)}}$ and for every $g \in \pos(f,\fc^p_0,
\dots,\fc^p_{n-1})$, $d \rest \pos(g, \fc)$ is constant,
\nothing{\item[($\delta$)]
(halving in the no gluing case) for every $\fc$ there is $\fc' \in
\Sigma(\fc)$ such that $\norm(\fc') \geq \norm(\fc) -
\frac{1}{k_{i(\fc)}}$ and if $\fd' \in \Sigma(\fc')$ and
$\norm(\fd') \geq 1$ then there is $\fd \in \Sigma(\fc)$ such that
$\norm(\fd) \geq \norm(\fc) - \frac{1}{k_{i(\fc)}}$ and $\fd$ and
$\fd'$ are equivalent,
 that is $i(\fd) = i(\fd')$ and $\pos(g,\fd) = \pos(g,\fd')$ for every
 $g$
 (both maybe undefined),}
 \item[($\delta$)] (halving with gluing) if $p \in \bP$, $m(\ast) <
 \omega$
 then we can find $q \in \bP$ with the following properties
 \begin{myrules1}
 \item[--] $p \leq q$,
 \item[--] $f^p = f^q$,
 \item[--] $\fc^p_m = \fc^q_m$ for $m < m(\ast)$,
 \item[--] if $m \geq m(\ast)$ then $\norm(\fc^q_m) \geq
 \inf\{\norm(\fc^p_\ell) \such \ell \in [m(\ast), \infty)\} -
 \frac{1}{k_{i(\fc^p_{m(\ast)})}}$,
\item[--] if $q \leq r$, $f^r=f^q$, $\fc^q_m = \fc^r_m$ for $m < m(\ast)$ and
$\norm(\fc^r_m) \geq 1$ for $m \geq m(\ast)$ then there is $q_1$
such that
\begin{myrules1}
\item[($*$)] $p \leq q_1$,
\item[($*$)] $f^{q_1} = f^p$,
\item[($*$)] $\fc^{q_1}_m = \fc^p_m$ for $m < m(\ast)$,
\item[($*$)] if $m \geq m(\ast)$ then $\norm(\fc^q_m) \geq
 \inf\{\norm(\fc^p_\ell) \such \ell \in [m(\ast), \infty)\} -
 \frac{1}{k_{i(\fc^p_{m(\ast)})}}$,
\item[($*$)] $q_1$ and $r$ are equivalent in a strong sense for some
$n(\ast) \geq m(\ast)$ we have $m \geq n(\ast)
 \rightarrow \fc^{q_1}_m = \fc^r_m$ and $\pos(f^{q_1},\fc^{q_1}_0, \dots,
 \fc^{q_1}_{n(\ast)-1}) =\pos(f^{r},\fc^r_0, \dots,
 \fc^{r}_{n(\ast)-1})$.
 \end{myrules1} %ast
 \end{myrules1} %--
 \end{myrules1} %greek
\end{myrules1} %latin
Then $\bP$ is $(T, \cS, Y)$-preserving.
\end{theorem}

\proof Assume that $\chi \geq 2^{2^\omega}$ and $N\prec \cH(\chi)$
is countable  and that $\bP \in N$, $p \in N \cap \bP$, $T, \cS, Y
\in N$. Let $\delta = N \cap \omega_1$ and $N \cap \omega_1
\in \cS$. Let $T$ be $(Y,\cS)$-proper. Assume that $p =
(\eta^p,\bc_0^p,\bc_1^p, \dots)$.

 We show that there is
$q \geq p$ that is $(N,\bP)$-generic and such that for every $t \in
Y(\delta)$, $$q \Vdash T_{<_T t} \mbox{ is
$(N[\name{\bG_\bP}],T)$-generic}.$$

We enumerate all pairs $(\name{\cI},t)$ of  $\bP$-names $\name{\cI}
\in N$ for dense sets in $T$ and $t \in Y(\delta)$ as
$\{(\name{\cI}_n,t_n) \such n < \omega \}$, all the $\cJ \in N$ that
are dense in $\bP$ as $\{\cJ_n \such n < \omega\}$, each object
appearing infinitely often in each enumeration.

We choose $(p_n, m_n)$ by induction on $n \in \omega$ such that
\begin{myrules1}
\item[(a)] $p_n \in \bP \cap N$,
\item[(b)] $p_n \leq p_{n+1}$,
\item[(c)] $m_n < m_{n+1}<\omega$, $m_0=0$,
\item[(d)] $p_n \Vdash (\exists t < t_n)
(t \in \bigcap_{k\leq n}\name{\cI_k})$,
\item[(e)] $p_{0} = p$,
\item[(f)] $p_n \in \cJ_n$,
\item[(g)] $f^{p_n} = f^p$,
\item[(h)] $m_n \geq \min \{ m> m_{n-1} \such (\forall r \geq
m)(\norm(\fc^{p_{n-1}}_r) \geq n +1) \}$,
\item[(i)] if $m < m_n$, then $\fc_m^{p_n} = \gc_m^{p_{n-1}}$.
\end{myrules1}

If we succeed then we can take the fusion $$q=(f^{p_0} \fc^{p_0}_0,
\dots \fc^{p_0}_{m_0-1},\fc^{p_1}_{m_0}, \dots, \fc^{p_1}_{m_1-1},
\dots)$$ and by (h) and (i), $q$ fulfils the norm conditions and
hence $q \in \bP$, and obviously $q \geq p$.

So suppose that $p_n$ and $m_n$  have been defined we are to define
$p_{n+1}$.

Let $i_n = i(\fc^{p_n}_{m_n})$ and let  $\{g_\ell \such \ell <
\ell_n\}$ list $\pos(f^{p_n},\fc_0^{p_n}, \dots,
\fc^{p_n}_{m_n-1})$. By the conditions on $\bP$, we have $\ell_n
\leq k_i$.

Now we choose $p_{n,\ell}$ by induction on $\ell < \ell_n$ such that
\begin{myrules1}
\item[(a)] $p_{n,\ell} \in \bP\cap N$,
\item[(b)] $p_{n,0} = p_n$,
\item[(c)] $f^{p_n} = f^{p_{n,\ell}}$,
\item[(d)] if $m < m_n$ then $\fc^{p_n}_m = \fc^{p_{n,\ell}}_m$,
\item[(e)] if $m \geq m_n$ then $\norm(\fc^{p_{n,\ell}}_m) \geq n +1
- \frac{1}{k_{i_n}}$,
\item[(f)] if there is $q=(g_\ell,\fc^q_{m_n}, \fc^q_{m_n+1}, \dots ) \geq
p_{n,\ell}$ such that $q \in \cI_n$ then $\fc^q_j =
 \fc^{p_{n,\ell+\frac{1}{2}}}_j$ for $j \geq m_n$; otherwise we apply
 halving to $(g_\ell, \fc^{p_{n,\ell}}_{m_n},
 \fc^{p_{n,\ell}}_{m_n+1}, \dots)$ and get $q$ as in the halving
 with gluing, and let again $\fc^q_j =
 \fc^{p_{n,\ell+\frac{1}{2}}}_j$ for $j \geq m_n$,
\item[(g)] if there is $q=(g_\ell,\fc^q_{m_n}, \fc^q_{m_n+1}, \dots ) \geq
p_{n+ \frac{1}{2},\ell}$ such that
 $q \Vdash (\exists t \in T) (t< t_n \wedge t \in
 \cJ_n)$ then $\fc^q_j =
 \fc^{p_{n,\ell+1}}_j$ for $j \geq m_n$; otherwise we apply
 halving to $(g_\ell, \fc^{p_{n,\ell+\frac{1}{2}}}_{m_n},
 \fc^{p_{n,\ell+\frac{1}{2}}}_{m_n+1}, \dots)$ and get $q$ as in the halving
 with gluing, and let again $\fc^q_j =
 \fc^{p_{n,\ell+1}}_j$ for $j \geq m_n$.
\end{myrules1}
It is easy to carry on the induction. In the end we let $p_{n+1} =
p_{n,\ell_n}$. Now we have to show that for each $\cI$ and each $(t,
\name{\cJ})$ (that appear under infinitely many indices) after
finitely many of these $n$ where e.g. $\cI_n = \cI$, in items (f)
and (g) the first alternative will be applied. This is because of
the strong version of bigness. We colour $\pos(\fc^p,\fc^{p_n}_0
\dots \fc^{p_n}_{m_n-1}, \fc^{p_n}_{m_n})$ by $\{0,1\}$ assigning
$c(\hat{g}) =1$ if there is $q$ with $f^q = \hat{g}$ (no conditions
on the rest) and $q \in \cI$ (in the case of (g): and $q \Vdash
(\exists t')( t' < t\wedge t' \in \name{\cJ})$). For every $\hat{g} \in
\pos(\fc^p,\fc^{p_n}_0 \dots \fc^{p_n}_{m_n-1})$ there is a uniform
colour. Now we go one level back: For ``most'' of the  $\hat{g} \in
\pos(\fc^p,\fc^{p_n}_0 \dots \fc^{p_n}_{m_n-1})$, their uniform
colour is the same, and for most of the most of the next level and
so on. So we get back to the root. Its colour is at some time $n$
the colour 1, since otherwise we succeed in constructing the fusion
$p_\omega$ that has no extension in $\cI$ or  no extension in $\{ r
\in \bP \such r \Vdash (\exists t' < t)(t' \in \name{\cJ})\}$, so
$\cJ^* = \{ s \in T \such s \not\leq t \vee p_\omega \Vdash s
\not\in \name{\cJ}\}$ is dense in $T$ and witnessing that $T$ is not
$(\cS,Y)$ proper. In any case this is a contradiction. So we get
$\fc'_{m_n}$ with large norm and colour 1 and are done. There are
$\ell_n$ substeps and in each step we lose maximally $\frac{1}{k_n}$
of $\norm(\fc^{p_n}_{m_n})$ so in the end it is still large enough
for a fusion. \proofend

\section{A sufficient condition for $(T,Y,\cS)$-preserving for nep
forcings}\label{S3}

The property of preserving Cohen generic reals over countable elementary
submodels 
proved to be a useful property of forcings. Preserving Cohen reals
is slightly stronger than preserving non-meager sets 
(see \cite[Section 3.2]{RoSh:470}). Preserving Cohen reals
is preserved in countable support iterations \cite[Ch.~XVIII, 3.10]{Sh:f}.
 In this section we show that a 
relative of this property, namely ``$\bP$ preserves $\omega$ Cohen reals
over countable elementary submodels and over certain transitive models called
candidates'', 
guarantees that $\bP$ preserves Souslin trees.
The candidates will replace the elementary $N \prec \cH(\chi)$. 
When a forcing notion $\bP$ has also for these non-elementary 
countable models suitable generic conditions then $\bP$ is called
``nep'' -- non-elementary proper. There are many versions of this definition:
We can specify which candidates are considered and which conditions
are imposed on genericity.  A standard reference to non-elementary proper forcing is \cite{Sh:630}.

Let $N \prec \cH(\chi)$. $x \in \baire$ is called Cohen over $N$, if
for every comeager $G_\delta$-set $C \subseteq \baire$ with code in $N$, 
$x \in C$. For Borel codes see \cite[Section~25, p.504]{Jech3}.
We recall the the original definition for 
proper forcing with elementary submodels, \cite[Def.~3.2.1]{RoSh:470}.

\begin{definition}\label{3.1} 
\begin{myrules1}
\item[(1)] Let $\bP$ be a proper forcing notion. We say $\bP$ is {\em
$\omega$-Cohen preserving} iff the following holds: For every $N
\prec \cH(\chi)$ such that $\bP \in N$, for every $p \in \bP\cap N$
for every $\{ x_n \such n \in \omega\}$ such that every $x_n$ is a
Cohen real over $N$, there is an $(N,\bP)$-generic condition $q \geq
p$ such that
$$q \Vdash (\forall n \in \omega) (x_n \mbox{ is Cohen over }
N[\name{\bG_{\bP}}]).$$
\item[(2)]
$\bP$ is {\em Cohen preserving} iff the above holds for just one
Cohen real.
\end{myrules1}
\end{definition}

By \cite[Ch.~XVIII, 3.10]{Sh:f} also $\omega$-Cohen
preserving is preserved in countable support iterations.
Cohen forcing itself is Cohen preserving, whereas random forcing is 
not, since the ground model reals are a meager set in the extension.
For creature forcings \cite[Ch.~3]{RoSh:470} gives some structural
properties on the building blocks of the forcing that
imply Cohen preserving.

The notion ``nep'' --- non-elementary proper --- was introduced and
investigated in \cite{Sh:630} and it is actually a reach family of
notions with many parameters.
We give a short introduction to our instance of nep.
Our presentation is a compromise between at least covering all the
creature forcings from \cite{RoSh:470}  and many technicalities.
 Explanations and useful work
with nep forcings can also be found in \cite{KrSh:828}.

In one respect we introduce more technique than needed for the
creature forcings from \cite{RoSh:470}: We like to allow a parameter
$\mathfrak B$ with domain $|\mathfrak B|\subseteq
\cH(\omega_1)$ and countable signature.

 Why are we so interested in allowing definitions with
parameters ${\mathfrak B}$?
In the light of the theorem in this section, an interesting question
is to consider to which extent forcings specialising Aronszajn trees
(by finite approximations, by countable approximation as in
\cite{MdSh:778} or by uncountable conditions as in the NNR forcing
from \cite[Ch.~V, \S 6]{Sh:f}) are nep.
Here are some partial answers:

In all ground models in which then NNR does not add
reals it is $\omega$-Cohen preserving.  However,
as our Theorem~\ref{3.11} shows, NNR is not
$\omega$-Cohen preserving in other models $M[G]$, where $G$ collapses $\omega_1$ of $M$ or it is not nep in the strong sense required in the theorem. 
The NNR
 forcings
are defined with Aronszajn trees as parameters in the definition. An
Aronszajn with its tree order tree can be written as a subset of
$H(\omega_1)$ and so still is a parameter allowed in size in the
definitions of nep we give. 

The forcings from \cite{MdSh:778} add a real that makes the ground
model meager (this is not yet published work by Mildenberger and Shelah), and hence they are not Cohen preserving.

\begin{definition}\label{3.2}
\begin{myrules}
\item[(1)]
A {\em fragment} $\zfc^*$ is an ${\mathcal L}(\in)$-theory extending ${\rm ZC}^-$,
$\zfc$ without replacement and without the power set axiom.
\item[(2)]
Let $K$ be a class of notions of forcing.
We say $\zfc^*$ is {\em $K$-good},
if $\bP$ is a forcing notion in $K$ and $\beth_\omega(|\bP|)$
exists then the forcing $\bP$ preserves $\zfc^*$.
\end{myrules}
\end{definition}

Now let $T$ be a  fragment of \zfc\ and $K$ be a definable class of forcing
notions or a set of forcings. Then by using the definability of
$\Vdash_{\bP}$ for $\bP \in K$ \cite[Ch.~VII, \S 4]{Kunen}
and adding successively the requirements $\Vdash_{\bP} \phi$ for $\phi$ in the
previous stage
we get a fragment $T_1 \supseteq T$ that $M \models T_1$ ensures
$M[G_\bP] \models T$.
Now we iterate and take the union. This need not be a finite fragment anymore.
So in practice, in order to get consistency relative to \zfc\, we take
$\zfc^\ast = \zc^-$. Then for every uncountable regular $\kappa$,
$(H(\kappa),\in)\models \zfc^*$. Now if no forcing in $K$ collapses
$\kappa$ to $\omega$ or to a singular cardinal, $\zfc^*$ is $K$-good.

We fix $\lambda = (2^{|H(\omega_1)|})^+$, $\chi = |H(\lambda)|^{++}$ the set
$K
= \{\Levy(\aleph_0, \lambda) \such \lambda \mbox{ regular uncountable
cardinal}, \lambda < \chi \}$ of notions of forcing.
\label{levy-key}

\begin{definition}\label{3.3}%(Cf.\ \cite[Def.~1.1.6]{Sh:630})
A theory $\zfc^\ast \subseteq \zfc$ is called {\em
normal} if the following holds:
\[
\mbox{For every sufficiently large regular }\chi, \cH(\chi) \models \zfc^\ast.
\]
\end{definition}

We assume that the
forcing $\bP$ is defined by formulas $\phi_0(x)$ and $\phi_1(x,y)$
that describe $x \in \bP$ and $x\leq_\bP y$. The formulas are in a
countable language $\tau \subset H(\omega)$ and use a
parameter ${\mathfrak B}\subseteq H(\omega_1)$.
We let $\bar{\phi} = ( \phi_0,\phi_1)$ for the description of
$\bP$ and $\leq_\bP$. In the 
stronger form of nep that is called ``explicit nep''
we have  $\bar{\phi} = ( \phi_0,\phi_1,\phi_2)$ with
$\phi_0$ and $\phi_1$ in the same roles, whereas the
additional first order formula $\phi_2$ describes
predense sets. 

\begin{definition}\label{3.4}%(See \cite[Def.~1.1]{Sh:630})
We call $(\bar{\phi},{\mathfrak B}, \zfc^*)$ a {\em definition of a
forcing}. We call $M$ a {\em $(\bar{\phi},{\mathfrak B},
\zfc^\ast)$-candidate} if $M$ is a countable transitive $\zfc^\ast$
model and $\mathfrak B \in M$.
\end{definition}

This is a simplification, since we say transitive.
The evaluation of $\bP$ over countable transitive models shall give
relevant information about its forcing behaviour in $\bV$.
Hence it is natural to require $\bP = \bigcup_{M\text{ a candidate}} \bP^M$,
where $\bP^M=\phi_0^M$. 
\relax From the requirement that $\phi_0$ is upwards absolute
we get $\bP^M= \bP \cap M$.  Then only $\bP \subseteq H(\omega_1)$
can fulfil the natural requirement.
Fortunately many well-known useful forcings with conditions of size
$\omega$ have $\bP\subseteq H(\omega_1)$. 
However already iterations of small iterands 
(i.e., with names in $H(\omega_1)$) of lengths $\geq \omega_2$ 
are not $\subseteq H(\omega_1)$ anymore. 
As a technical means to handle this situation one can use ord-hc candidates.
$M \models \zfc^*$ instead of transitive models. 
We refer the reader to \cite{Sh:630} and \cite{Kr2009}, and we will work
here only with transitive models.

 In the next section, we show that
$(T,Y,\cS)$-preserving is an iterable property. So it is enough to
give a sufficient criterion for
 $(T,Y,\cS)$-preserving just for one nep iterand. Iterands usually
are small and we do not lose any of the creature forcings.

\begin{definition}\label{3.5}%(\cite[Def.~1.3 (1)(c)]{Sh:630})
If $M$ is a candidate then $G\subseteq \bP^M = \{ p \in M \such M
\models \phi_0(p) \}$ is {\em $(M,\bP)$-generic} if for all $A \in
M$, if $ M \models \mbox{``}A \subseteq \bP$ is a maximal
antichain'', then $|G \cap A| = 1$. (The incompatibility in $\bP$
might be not absolute, so $G \cap A \neq \emptyset$ is not enough.)
$q$ is called {\em $(M,\bP)$-generic} if for all $A \in M$ such that
$M \models\mbox{``} A$ is a maximal antichain, $q \Vdash |G \cap A|
=1$.
\end{definition}

\begin{definition}\label{3.6}%(\cite[Def.~1.3 (1)]{Sh:630})
Let $K$ be a class of forcings. $(\bar{\phi},{\mathfrak B},
\zfc^\ast)$ is called a $K$-{\em $(\bar{\phi},{\mathfrak B},
\zfc^\ast)$-definition of a nep forcing} if the following holds in
$\bV$ and in all extensions of $\bV$ by members of $K$:
\begin{myrules1}
\item[(a)]
$\phi_0$ defines the set  of elements of $\bP$  and $\phi_0$
is upwards absolute from candidates to $\bV$,
 in $\bV$ and in all $(\bar{\phi},{\mathfrak B},
 \zfc^\ast)$-candidates,
\item[(b)] $\phi_1$ defines the quasi ordering $\leq_\bP$ in $\bV$ and in
every $(\bar{\phi},{\mathfrak B},
 \zfc^\ast)$-candidate,
$\phi_1$ is upwards absolute from candidates to $\bV$,
 in $\bV$ and in all $(\bar{\phi},{\mathfrak B},
 \zfc^\ast)$-candidates,
\item[(c)]
if $M$ is a $(\bar{\phi}, {\mathfrak B}, \zfc^\ast)$-candidate and
$p \in \bP^M$ then there is an $(M,\bP)$-generic $q \geq p$.
\end{myrules1}
\end{definition}

We isolate a property:

\begin{myrules1}
\nothing{\item[$(\ast_1)$] If $M_i$, $i =1,2$, are
$(\bar{\phi},{\mathfrak B}, \zfc^\ast)$- candidates and $M_2 \models
\mbox{``}M_1$
 is a $(\bar{\phi},{\mathfrak B}, \zfc^\ast)$-candidate''
and $M_1 \models \phi_0(x)$ then $M_2 \models \phi_0(x)$ and the
same holds for $\leq$.}
\item[$(\heartsuit)$]
If $M_1$ is a  $(\bar{\phi},{\mathfrak B}, \zfc^\ast)$-candidate
and
$M_1 \models \mbox{``}M_0$
 is a $(\bar{\phi},{\mathfrak B}, \zfc^\ast)$-candidate and $p \in \bP^{M_0}$'' then
then  there is $q \in \bP^{M_1}$, $q\geq p$ such that
  $M_1 \models $``$q$ is $(M_0,\bP)$-generic'' and such that
in $\bV$, $q$ is $(M_0,\bP)$-generic.
\end{myrules1}

In the following we show that $\heartsuit$ follows from quite natural strengthenings of the notion of non-elementary properness. Many well-known forcings
are non-elementary proper in one of these strong variants.

\begin{definition}\label{3.7}%(See \cite[Def.~1.3 (2)]{Sh:630})
 We add the adverb ``explicitly'', so say
``$(\bar{\phi},{\mathfrak B}, \zfc^\ast)$ is called a
$(\bar{\phi},{\mathfrak B}, \zfc^\ast)$-definition of an
$K$-explicitly nep forcing'' if $\bar{\phi} = \langle \phi_0,
\phi_1,\phi_2\rangle$ and
$(\phi_0,\phi_1)$ are as in Definition~\ref{3.6} and additionally
\begin{myrules1}
\item[$(b)^+$] We assume $\phi_2$ is an $(\omega+1)$-place relation 
that is upward
absolute from $(\bar{\phi},, {\mathfrak B}, \zfc^\ast)$-candidates. 
$\phi_2(p_i \such i\leq \omega)$ says  $\{p_i \such i\leq \omega\}
\subseteq \bP$ and $\{p_i \such i<\omega\}$ is a predense antichain
above $p_\omega$ not just in $\bV$ but in every $(\bar{\phi},{\mathfrak B},
\zfc^\ast)$-candidate which satisfies $\phi_2( p_i \such i\leq
\omega)$. In this situation we say $\{p_i \such i< \omega \}$ is
explicitly predense above $p_\omega$.
\item[(c)$^+$]
We add to (c) in the definition of nep:
There is $q \geq p$ with the following property: If $N \models \cI$ is a
predense antichain above $p$, so $\cI \in N $then for some list
$\langle p_i \such i < \omega \rangle$ of $\cI \cap N$ we have
$\phi_2(\langle p_i \such i < \omega \rangle \concat \langle q
\rangle)$. We then say ``$q$ is explicitly $(N, \bP)$-generic above
$p$.''
\end{myrules1}
\end{definition}

In our proof $K$ contains also the
Levy  collapse, so not only mild forcings. So as soon as the
definition of the forcing $\bP$ is sensitive to cardinals, $K$-nep
becomes a strong requirement. Think for example again of the forcing
specialising a normal Aronszajn tree: After collapsing, the Aronszajn tree
is just a perfect tree $\subseteq\omega^{<\omega}$.

So finally to get $(\heartsuit)$ we need even more than explicitly
nep:

\begin{definition}\label{3.8}
%(See \cite[Def.~5.13 (1), (3)]{Sh:630}
\begin{myrules1}
\item[(1)] 
A $(\bar{\phi},{\mathfrak B}, \zfc^\ast)$-definition of a
forcing notion $\bP$ is called {\em straight} nep if is it
$K$-explicitly nep and in addition
 For $\ell< 3$ the formula $\phi_\ell$ is of the form
\[
(\exists t)[t \in \cH(\omega_1)
\wedge (\exists s)((s \in t \vee s=t ) \wedge
\psi_\ell^{ \bP}(\bar{x},s))],
\]
 where in the formula
$\psi^{\bP}_\ell$ the quantifiers are of the form ($\exists s' \in
s)$ and the atomic formulae are $x \in y$, ``$x$ is an ordinal'',
``$x<y$ are ordinals'' and those of ${\mathfrak B}$. 
\item[(2)]
We say {\em very straight} if it is straight and in addition
\begin{myrules1}
\item[(f)] for some Borel %(see \cite[Def.0.5]{Sh:630})
 functions $\bB_1$, $\bB_2$, if $N$ is a
candidate and  $\bar{a}$ lists $N$ and $p \in \bP^N$, then $q =
\bB_1(p,\bar{a}, N)$ is explicitly $(N,\bP)$-generic and
$\bB_2(p,\bar{a},N)$ is a witness, that is it witnesses $p \leq q$
and $\phi_2(\langle p_{\cI,n} \such n < \omega\rangle, q)$ for some
sequence $\langle p_{\cI,n} \such n < \omega\rangle$ of members of
$\cI$ for every predense antichain $\cI$ of $\bP^N$ in $N$.
\end{myrules1}
\end{myrules1}
\end{definition}

The property from Definition~\ref{3.8}(1) guarantees: If $p,q \in
M_1$ and $p \leq q$ in $\bV$, then  $M_1 \models p \leq q$.
Upwards absoluteness is included in the more basic canon of nep properties
Definition~\ref{3.6} (a), (b).

The following lemma shows that there are many examples of forcing notions 
that meet our version non-elementary properness.
Its proof is long and will not be repeated here.

\begin{lemma}\label{3.9}%  (\cite[Lemma~3.1]{Sh:630})
 We can use $\zfc^\ast ={\rm ZC}^-$ which is $K$-good
for $K$ from Page~\pageref{levy-key}
and normal.
\begin{myrules1}
\item[(1)]
Suppose that $\bP$ is a forcing of one of the following types:
\begin{myrules1}
\item[(a)]
$\bQ^{\rm tree}_e(K,\Sigma)$  for some finitary tree creating pair
$(K,\Sigma)$ that is t-omittory without a condition on the norm
\nothing{something was missing, otherwise it need not even be
proper} for $e=0$ and 2-big in the case of $e=1$ (this covers Sacks
forcing).
\item[(b)] 
\nothing{$\bQ^\ast_{{\rm s},\infty}$ for some finitary
creating pair $(K,\Sigma)$ which is growing \cite[Def.~2.1.1
(3)]{RoSh:470} condensed \cite[Def.~4.3.2.(3)]{RoSh:470} and of the
AB-type, or of the AB$^+$ type and satisfies $\oplus_0$, $\oplus_3$
of \cite[4.3.8]{RoSh:470}. This captures}
the Blass--Shelah  forcing notion.
\item[(c)] $\bQ^\ast_{{\rm w},\infty}(K,\Sigma)
$ for some finitary creating pair which
captures singletons, that is ($K,\Sigma)$ is forgetful and for every
$(t_0,\dots t_n)$ and for each $u \in \basis(t_0)$ and $v\in
\pos(u,t_0,\dots t_n)$ there is ($s_0,\dots, s_k)$ such that
$(t_0,\dots t_n) \leq (s_0,\dots s_k)$ and $m_{\rm dn}^{t_0} =m_{\rm
dn}^{s_0}$ and $m_{\rm up}^{t_n} = m_{\rm up}^ {s_k}$ and
$\pos(u,s_0,\dots, s_k)= \{v\}$. $(K,\Sigma)$ is forgetful if for
every $t\in K$ and $\la w,u \ra \in \val[t]$ and $ w' \in
\prod_{n<m_{\rm dn}^t} H(n)$ also $\la  w', w' \concat u \rest
[m_{\rm dn}^t, m_{\rm up}^t)\ra  \in \val (t)$.
%\cite[Def.~2.1.10]{RoSh:470}.
\nothing{\item[(d)] $\bQ^\ast_f(K,\Sigma)$ for some finitary 2-big creating
pair with the Halving Property which is either simple
\cite[2.1.7]{RoSh:470}  of gluing \cite[2.1.7]{RoSh:470} and an $H$-fast
function $f \colon \omega \times \omega \to \omega$.}
\end{myrules1}
Then $\bP$ is an explicitly nep very straight forcing notion with a
Souslin definition  (see \cite[Def.~1.9]{Sh:630}).
 \item[(2)]
Suppose that $\bP$ is a forcing of one of the following types:
\begin{myrules1}
\item[(a)]
$\bQ^{\rm tree}_e(K,\Sigma)$  for some countable tree creating pair
$(K,\Sigma)$ that is t-omittory without a condition on the norm
(see Def.~\ref{2.15})
\nothing{We get more examples by weakening t-omittory. 
This is o.k.
by the rework of the tree forcing \cite[Prop 2.3.6]{RoSh:470}
according to the remarks in Section \ref{S2}.} for $e=0$ and 2-big
in the case of $e=1$ (this Miller forcing and Laver forcing).
\item[(b)]
$\bQ^\ast_{\infty}(K,\Sigma)$ for some finitary growing
%\cite[Def.~2.1.1 3.]{RoSh:470}
pair $(K,\Sigma)$. This covers the Mathias forcing notion.
$(K,\Sigma)$ is called growing if for any sequence $(t_0,\dots, t_{n-1})$
with $m_{\rm dn}^{t_i} = m_{\rm up}^{t_{i-1}}$ for $i < n$ there is
$t\in \Sigma(t_0,\dots, t_{n-1})$ such that $\norm(t) \geq
\max_{i<n} \norm(t_i)$.
\end{myrules1}
Then $\bP$ is an explicitly nep very straight forcing notion.
\end{myrules1}
\end{lemma}

\nothing{
\begin{proposition}\label{4.15}%(\cite[Proposition~3.2]{RoSh:470})
All the proper forcing notions $\bP$ defined in \cite{RoSh:470,
RoSh:628} are simple\nothing{, correct (see \cite[Def.~1.3 (11)]{Sh:630})},
very straight and we can use $\zfc^\ast ={\rm ZC}^-$ which is $K$-good
and normal. Also the relations $p \perp q$ in $\bP$ is upwards
absolute from candidates (as well as $p \in \bQ$, $p \not\in \bQ$,
$p \leq q$, $p \not\perp q$).
\end{proposition}
}
For a proof see
(\cite[Proposition~3.2]{RoSh:470}.

\begin{definition}\label{3.10}
Suppose $\mathfrak B$ is a model with domain $\subseteq
\cH(\omega_1)$, $\bP$ is an $(\bar{\varphi}, {\mathfrak B},
\zfc^*)$ definition of a nep forcing. $\bP$ is called
\emph{$\omega$-Cohen
preserving for $(\bar{\varphi}, {\mathfrak B}, \zfc^*)$-candidates},
iff the following holds: If $N$ is a $(\bar{\varphi}, {\mathfrak B},
\zfc^*)$-candidate and for each $n$, $x_n \in {}^\omega \omega$ is a
Cohen real over $N$ and $p \in \bP^N$ then there is $q \in \bP$, $q
\geq p$ that is $(N,\bP)$-generic and $q \Vdash (\forall n) (x_n $
is a Cohen real over $N[\name{{\bG}_\bP}])$.
\end{definition}

For proper forcings that are given 
by a definition $(\bar{\varphi}, {\mathfrak B},\zfc^\ast)$
 this is a strengthening of Def.~\ref{3.1}, since
for countable $M \prec H(\chi)$, the transitive collapse is a candidate.

So finally we can state the main theorem in this section:

\begin{theorem}\label{3.11}
Suppose $\mathfrak B$ is a model with domain $|{\mathfrak B}|
\subseteq  H(\omega_1)$ and countable signature, the definition
$(\bar{\varphi}, {\mathfrak B}, \zfc^*)$ of $\bP$ is explicitly very
straightly nep, $\zfc^\ast$ is normal and $K$-good.
Suppose that $\bP$ is $\omega$-Cohen preserving for $(\bar{\varphi},
{\mathfrak B}, \zfc^*)$-candidates.
 Then $\bP$ is $(T,Y,\cS)$-preserving for
all triples $(T,Y,\cS)$.
\end{theorem}

\proof Let $\lambda\geq 2^{|H(\omega_1)|}$ be large enough such that
$(\forall \lambda'\geq \lambda )(\cH(\lambda') \models \zfc^*)$. Let
$ \lambda_1 = |\cH(\lambda)|$, and let
 $\chi > \lambda_1 \geq \lambda$.  Let $N' \prec
\cH(\chi)$ be countable such that $\{\lambda_1,p,\cS,T,
\bar{\varphi},{\mathfrak B}\}\subset N'$. 
Our aim is to show that $N'$ is as in Def.~\ref{1.16}.

Let $N$ be the Mostowski
collapse of $N'$ and say $\pi^{N'} \colon N' \to N$ is
the collapsing function. Let $\delta = N' \cap \omega_1$.
So $N\models \delta =\aleph_1$.
Assume that $N' \cap \omega_1 \in \cS$.

\medskip

We let $\bR = \Levy(\aleph_0,
\pi(\lambda_1))^N$.

Claim 1: In $\bV$, there is $\bg$ such that

\begin{myrules1}
\item[(a)] $\bg$ is $\bR$-generic over $N$,
\item[(b)]
if $t \in \pi^{N'}(Y(\delta))$ then $\pi^{N'}(T_{<t})$ is $\pi(T)$-generic
over $N[\bg]$.
\end{myrules1}

Proof: $\bR$
is a forcing notion in $N$ and hence in $\bV$.
Let $\bg_1$ be generic over $\bV$ not just over $N$. We first show
that $\bg_1$ would be as desired, were it in $\bV$. In $\bV[\bg_1]$,
let $t \in \pi^{N'}(Y(\delta))$. We show that $\pi^{N'}(T_{<t})$ is generic
over $N[\bg_1]$. Let $\cI \in N[\bg_1]$ be a subset of $\pi(T_{<\delta})$
such that $N[\bg_1] \models \cI \mbox{ is dense and open in  the
forcing } \pi((T,<_T))$. By the forcing theorem there is $p \in
\bg_1\cap \bR$ such that $N \models p \Vdash_\bR \name{\cI} \mbox{
is dense and open}$. Assume towards a contradiction that $\bV
\models [p \Vdash_\bR \pi(T_{<t}) \cap \name{\cI} = \emptyset].$

The set
\[\cI^* = \{q \in {\rm Cohen} \such
(\exists \nu \in \pi^{N'}(T_{<t}))( q \not\Vdash_\bR \nu \not\in \cI)\}
\]
is dense and open in the Cohen poset, since
$\bR$ is just Cohen forcing, and the iteration of two Cohen
forcings is equivalent to the iteration in the reversed order. So
there is $\nu \in \pi^{N'}(T_{<_T t})$  and there is $q \geq p$ such that
$ q \Vdash_{\bR} \nu \in \cI$ which is a contradiction.

Now $\bg_1$ is not in $\bV$. The requirements on $\bg_1$ have only
quantifiers bounded by sets and hence are absolute for
transitive models
\begin{multline*} (\forall D \in N)(D \mbox{ dense in
}\Levy(\aleph_0,\pi(\lambda_1))^N) \rightarrow D \cap \bg_1 \neq
\emptyset) \wedge\\
 (\forall t\in \pi(Y(\delta))) (\forall D \in N[\bg_1] \mbox{ that are dense in ${\pi}(T))$})
(\exists s \in \pi(T_{<t}) \cap D).
\end{multline*}
 So ``there is such a
$\bg$ with these properties'' is a $\Sigma_1$
sentence with parameters in $\bV$ 
that is true in $\bV[\bg_1]$. By absoluteness, it is true also in $\bV$.

Let $t \in Y\cap \pi^{N'}(T)$ be given.
Note that $\pi^{N'}(T) = T_{<\delta}$. By the assumption that
$T$ is $(Y,\cS)$-proper, and hence after the
Levy collapse, $\pi(T_{<_T t})$ is $(N,\pi^{N'}(T))$-generic.
Let $M = N[\bg]$ with a $\bg$ as in the claim. Then $\pi^{N'}(T_{<t})$ is
also $(M,\pi^{N'}(T))$-generic by the choice of $\bg$. $M$ is  a
candidate since $\bP$ is nep as in the condition of the
theorem and as $\bR  \in K$.

Now we use that $\bP$ is $\omega$-Cohen preserving for the
$(\bar{\phi},{\mathfrak B}, \zfc^\ast)$-candidate $M$.
 We choose a dense embedding $h \colon (\omega^{<\omega}, \triangleleft) \to
{\pi}^{N'}(T, <_T)$, $h \in N[\bg]$. Note the $N[\bg]$ thinks that 
${\pi}^{N'}(T)$
is countable, since $\bR$ collapses $\pi^{N'}(\omega_1)<\pi^{N'}(\lambda_1)$ from $N$
to $\omega$. So
in $N[\bg]$ the Cohen forcing
$(\omega^{<\omega}, \triangleleft)$ and ${\pi}^{N'}(T, <_T)$
are equivalent. We let
$\eta_t \in {}{^\omega }\omega$ be such that $n < \omega \rightarrow
h(\eta_t\rest n ) < t$. So $\eta_t$ is Cohen over $M$ iff
$\pi^{N'}(T_{<_T t})$ is $(M,\pi^{N'}(T))$ generic, and this holds also for
extensions of $M$ since there is the isomorphism $h$ in them. Since
$\pi^{N'}(T_{<_T t})$ is $(M,\pi^{N'}(T))$-generic, $\eta_t$ is Cohen generic
over $M$. Now we use that $\bP$ is Cohen preserving for the
candidate $M$. So there is $q \geq \pi(p)$ that is $(M,\bP)$-generic and
$q \Vdash (\forall t \in \pi^{N'}(Y(\delta)))( \eta_t$ is Cohen over
$M[\name{\bG}_{\bP}]$). So
\begin{multline}\label{abs1}
q\Vdash
(\forall t \in \pi(Y(\delta)))( \pi(T_{<_T t})\mbox{ is }
(M[\bG_\bP],\pi(T))\mbox{-generic})\\
\mbox{ and $q$ is $(M,\bP)$-generic}.
\end{multline}

Now we get from the latter
\begin{multline}\label{abs2}
(\exists q_3 \geq \pi(p))(
 q_3 \Vdash \mbox{``}(\forall t \in \pi(Y(\delta)))\\ \pi(T_{<_T t}) \mbox{ is
$(N[\name{\bG}_\bP],\pi(T))$-generic'' and $q_3$ is $(N,\bP)$-generic}).
\end{multline}

 Why?
We use nep again. We take $\chi_1$ such that $N \models (\chi_1$ is
sufficiently large such that ${\mathcal P}(\bP) \in H(\chi_1)$ and
$\chi_1$ is sufficiently small so that
$2^{\chi_1}$ exists). Let $N_1 = N \rest H(\chi_1)^N$. In $N$, $N_1$
is a candidate.

 By $(\heartsuit)$ there is
$q_1 \geq \pi(p)$, $q_1 \in N\subseteq M$, $N \models \mbox{``}\pi(p) \leq
q_1$ and $q_1$ is $(N_1, \bP)$-generic'' and $q_1 \geq \pi(p)$ also in $\bV$ and
$q_1$ is $(N_1,\bP)$-generic also in $\bV$.

We claim: $q_1$ is as required in the first half of
\eqref{abs2}, that is: For any $\cI \in M[\name{\bG}_\bP]$ that has a $\bP$-name in $N$ and is (forced by the weakest condition to be)
a dense set in $\pi(T,<_T)$
in the sense of $M[\name{\bG}_\bP]$,
$q_1
\Vdash (\forall t \in \pi(Y(\delta)))(\name{\cI} \cap \pi(T_{<_T t})
\neq \emptyset)$. All $\name{\cI} \in N$ that are $\bP$-names for
dense sets in $\pi(T,<_T)$ have names $\name{\cI}
 \in N_1$. Now we argue in $N$: For
any $q_2 \geq q_1$, $q_2 \in M$,
 there is $ q \geq q_2$, $q$ is
$(M,\bP)$-generic and $q \Vdash (\forall t \in
\pi(Y(\delta)))(\name{\cI} \cap \pi(T_{<_T t}) \neq \emptyset)$, as
we have shown above, in Equation~\eqref{abs1}, proved for
 $q_1$ instead of $\pi(p)$
and proved in $N \models \zfc^*$ instead of in $\cH(\chi)$ and $\bV$.
Then $q \Vdash \pi^{N'}(T_{<_T t}) \cap \name{\cI} \cap N_1\neq
\emptyset$ since $(\pi^{N'}(T_{<_T t}))^{N[\bG_{\bP}]} = (\pi(T_{<_T
t}))^{M[\bG_{\bP}]} \subseteq N_1$. So $q \Vdash \pi^{N'}(T_{<_T t}) \cap
\name{\cI} \cap N\neq \emptyset$ and hence $\pi^{N'}(T_{<_T t}) $ is
$(N[\bG_{\bP}], \bP)$-generic.
Since $N$ is elementary equivalent to $\cH(\chi)$, also in $\bV$ we have: 
For every $t \in \pi^{N'}(Y(\delta)$,
for every $\cI \in M[G]$ that has a $\bP$-name in $N$ and is a dense
set in $\pi^{N'}(T,<_T)$, $q_1 \Vdash (\forall t \in
\pi(Y(\delta)))(\name{\cI} \cap \pi(T_{<_T t}) \neq \emptyset)$. 
So $q_1$ has the property required of $q_3$ in \eqref{abs2}.

Now we use nep again and find $q_3\geq q_1$ that is
 $(N,\bP)$-generic.
So $q_3$ witnesses that \eqref{abs2} is proved.
 Now after taking the
reverse image of the Mostowski collapse (the nep forcing $\bP$
is moved from model to model by just taking its interpretation) we have
$q_4 \geq p$ such that
\begin{multline*} q_4 \Vdash (\forall t \in Y(\delta)) T_{<_T t} \mbox{ is
$(N'[\name{\bG}_\bP],T)$-generic'' and $q_4$ is $(N',\bP)$-generic}.
\end{multline*}
\proofend

\begin{remark}\label{3.12}
In the special case that the $Y \cap T_\delta$ is a singleton (or
empty) for all $\delta \in \cS$, we need only a weaker form of Cohen
preserving, with one Cohen generic $\eta$. In this special
case ``$T$ is $(Y,\cS)$-proper'' implies $T \rest \{\sup(a) \such a \in
\cS \wedge Y(\sup(a)) \neq \emptyset \}$ has no specialisation,
\nothing{, and in the reverse direction, if $T \rest \{\sup(a) \such
a \in \cS\}$ has no specialisation, then there is a $Y$ with
$\{\delta \such Y(\delta) \neq \emptyset\} \supseteq \{\sup(a) \such
a \in \cS\}$ such that $T$ is $(Y,\cS)$-proper,} see
\cite[Ch.~IX]{Sh:f}.

We may also consider the well-known stronger variant of
$(Y,\cS)$-properness for forcing with finite products of $T$: If $N
\prec \cH(\chi) $ and $\delta =  N \cap \omega_1$ and $N \cap
\omega_1 \in \cS$ and $t_0, \dots,
 t_{n-1}$ are pairwise distinct then $\{\bar{s} \in {}^n (T_{<\delta}) \such \bigwedge_{\ell<n}
s_\ell <_{T} t_\ell\}$ is  $(({}^n(T_{<\delta}),N)$-generic. Also
preserving this kind of $(Y,\cS)$-properness in a consequence of
$\omega$-Cohen-preserving and nep, by the same proof as above. By an analogue
of the results of the next section, this preservation property is
iterable.
\end{remark}

The following theorem is similar to Theorem~\ref{main}, and in the
case of nep forcing it strengthens Theorem~\ref{main} by adding the
intermediate step in the implication $\Pr^2(\bP) \rightarrow \bP$ is
$\omega$-Cohen preserving $\rightarrow \bP$ is
$(T,Y,\cS)$-preserving.

\begin{theorem}\label{3.13}%wie label main that appear at 1.17
\begin{myrules1}
\item[(1)]
If $\Pr^2(\bP)$, then $\bP$ is $\omega$-Cohen preserving.
\item[(2)]
Assume $\alpha(\ast)=\omega_1$ and $\cS \subseteq
[\omega_1]^{\omega}$ is stationary.  If $\Pr^2_\cS(\bP)$, then $\bP$
is $\omega$-Cohen preserving for Cohen reals over $N$ with $N \cap
\omega_1 \in \cS$.
\end{myrules1}
\end{theorem}
Proof (1):
Assume $N \prec H(\chi)$, $N \cap
\omega_1 \in \cS$,  $N \cap \omega_1 = \delta$, and 
$\bP \in N$, $p \in N\cap \bP$, and assume
 for every $i <
\omega$, $x_i$ be Cohen generic over $N$. Let $x = {\bf st}$ for a winning strategy ${\bf st}$ for
player COM in $\Game^2(N,\bP,p)$. We show that there is a $q$ as required.

\smallskip

 Let $\{ \name{\cI}_k \such k \in \omega \}$ list all
$\bP$-names in $N$ of comeager sets  and let $\{\cJ_n
\such n \in \omega\}$ list all the dense sets in $\bP$ in $N$.  Now
take a play $\langle
\langle \bar{p_n},\bar{q_n}\such n \in \omega
\rangle$ in which COM plays
according to ${\bf st}$. 
By Lemma~\ref{1.4} COM can strengthen his moves and still wins. 
COM plays in every round $n$ in every
part $p_{n,\ell}$, $\ell <\ell_n$, so strong that
$p_{n,\ell} \in \bigcap_{r<n}\cJ_r$ such that for every $i < n$
\[
p_{n,\ell} \Vdash_{\bP} x_i \in \bigcap_{k< n}\name{\cI}_{k}.
\]

Such $p_{n,\ell}$  exist for the following reason:
Since $ \bigcap_{k< n}\name{\cI}_{k}$ is (forced by
the weakest condition to be) comeagre,  for every $n \in
\omega$, the set $\cJ_{n} = \{s \in {\mathbb C} \such \{ q \in \bP
\such q\not\Vdash_{\bP} [s] \cap \bigcap_{k<
n}\name{\cI}_{k}=\emptyset\}$ is dense and open in $\bP\}$ is
open and dense in the Cohen forcing $\mathbb C$ in the ground model.
The Cohen real $x_i$ fulfils $x_i \in \cJ_{n}$. So for every $i$, $\{ q \in \bP
\such q\not\Vdash_{\bP} x_i \not\in \bigcap_{k<
n}\name{\cI}_{k}\}$ is dense in $\bP$.

COM wins the play because he played according to the strategy. So
for every $u$, in particular for $u = \omega$,  there is $q_u \geq
p$ such that
\begin{equation}\label{1d}
q_u \Vdash (\exists^\infty n \in u )(\exists \ell < \ell_n)
(p_{n,\ell} \in \name{\bG}_\bP).
\end{equation}
Let $k \in \omega$ and $q'\geq q_u$ be given.
Then there is $q{''}\geq q'$ and $n \geq k$ such that $q{''} \Vdash
n\in u$. 
So there is $i <\ell_n$, $q{''}\Vdash q_{n,i}\in \name{\bG}_\bP$ and hence
\begin{equation}\label{2neu}
q{''} \Vdash_{\bP}
x_{k} \in \bigcap_{k'< n}\name{\cI}_{k'}.
\end{equation}
Now we unfreeze $k$ and combine the equations
\eqref{1d} and \eqref{2neu} and thus get
\[
q_u \Vdash (\forall k < \omega)(x_k \mbox{
is } (N[\name{\bG}_\bP],T)\mbox{-generic.})
\]
\relax From $q_{n,i}\in \bigcap_{r<n} \cJ_r$ we also get
that $q_u$ is $(N,\bP)$-generic. \proofend

\begin{remark}\label{3.14}
$\omega$-Cohen preserving is not a necessary condition for preserving Souslin trees:
By \cite[Lemma~3.1]{Sh:630}, Blass--Shelah forcing is nep in the
strong form that is used in Theorem~\ref{3.11}. Blass--Shelah forcing
is not Cohen preserving. This follows from the fact that the generic
real is not split by any real in the ground model (see
\cite{BsSh:242}). Hence the ground model is meagre after Blass
Shelah forcing. So Blass-Shelah forcing is not positivity preserving
for the meagre ideal in the sense of \cite[Def.~3.1]{KrSh:828}. So
by \cite[Lemma~5.6]{KrSh:828} is it not true positivity preserving
\cite[Def.~5.5]{KrSh:828}. Now by \cite[Lemma~5.8]{KrSh:828}
Blass-Shelah does not preserve the Cohen genericity.
 So it is a nep forcing not covered by
Theorem~\ref{3.11}. Nevertheless Blass--Shelah preserves Souslin trees
by Theorem~\ref{2.32}. 
\end{remark}

\section{Preserving the Souslinity of an $\omega_1$-tree}
\label{S4}

The topic of the section is the preservation of 
properties of notions of forcing in countable support iterations.
We return to the  $\omega_1$-trees and the
properties of $(T,Y,\cS)$ from the first section. In this section we
give a self-contained proof of the following theorem:

\begin{theorem}\label{4.1}
%(\cite[Ch.~XVIII, Conclusion 3.9.F]{Sh:f})
 Let $\langle \bP_\alpha,
\name{\bQ}_\beta \such \alpha \leq \gamma, \beta < \gamma \rangle$
be a countable support iteration of proper forcings. Suppose that
$T$ is a Souslin tree in $\bV$ and that for every $\alpha < \gamma$,
in $\bV^{\bP_\alpha}$,  $\Vdash_{\bQ_\alpha} \mbox{``}T \mbox{ is
Souslin''}$. Then $T$ is Souslin in $\bV^{\bP_\gamma}$.
\end{theorem}

By Lemma~\ref{1.15} $T$ is Souslin in $\bV^{\bP_\gamma}$ iff
$T$ is $(Y,\cS)$-proper for a stationary $\cS\subseteq[\omega_1]^\omega$
and $Y= \bigcup\{T_{\sup(a)}\such a\in \cS\}$. 
So Theorem~\ref{4.1} is a special
case 
\nothing{for $\cS$ being stationary in $[\omega_1]^{\omega}$, $Y =
\{t_n^\delta, n < \omega, \delta \in W\}$, $Y(\delta) = \{t_n^\delta
\such n < \omega\}=T_\delta$ for $\delta \in W$ and $W \subseteq
\{\sup(a) \such a \in \cS\}$ stationary, of the proof} 
of the
following theorem
\nothing{ that is not more complex than the proof of
Theorem~\ref{4.1}, which is}
 (\cite[Ch.~XVIII, Conclusion 3.9 F]{Sh:f}):

\begin{theorem}\label{4.2}
%(\cite[Ch.~XVIII, Conclusion 3.9 F]{Sh:f})
 Let $\langle \bP_\alpha,
\name{\bQ}_\beta \such \alpha \leq \gamma, \beta < \gamma \rangle$
be an countable support iteration of proper forcings. Suppose that
$T$ is an $\omega_1$-tree that is $(Y,\cS)$-proper and that for
every $\alpha < \gamma$, in $\bV^{\bP_\alpha}$, $\name{\bQ}_\alpha$
is $(T,Y,\cS)$-preserving. Then $\bP_\gamma$ is
$(T,Y,\cS)$-preserving and  $T$ is $(Y,\cS)$-proper in
$\bV^{\bP_\gamma}$.
\end{theorem}

The proof of Theorem~\ref{4.2} is not more complex than the proof of 
Theorem~\ref{4.1}.
It involves preserving unbounded families in certain relations.
A simple 
similar example is to preserve a $\leq^*$-unbounded family.
The relations are binary relations on spaces ${}^a a$ for 
countable sets $a$. For the proof of Theorem~\ref{4.2} we need
that the union of all considered $a$ covers $\omega_1$.
For each fixed $a$, ${}^a a$ is just a copy of the
Baire space $\omega^\omega$, which means $a$ has the discrete topology and
${}^a a$ carries the product topology.
 We consider $\aleph_1$ different sets $a$, and on each fixed 
$a$ we work with countably many relations 
$R_{\alpha,a}$, $\alpha \in a$.

Let $\cS \subseteq [\omega_1]^\omega$.
Let for $a \in \cS$, ${\bf g}_a \in {}^a a$ and for $\alpha \in a$,
$R_{\alpha,a} \subseteq {}^a a  \times {}^a a $ be a relation.

We assume that 
 a fixed ${\bf g}_a$, the sets 
\begin{equation}
\label{null}
\{f \in {}^a a \such  f R_{\alpha,a} {\bf g}_a\}
\end{equation}
are closed in ${}^a a$. 
We will see that open relations do not harm since they are
the union of countably many closed relations.
Let $N \prec (H(\chi),\in,<)$.
Now suppose that we have a collection $\cG$ of ${\bf g}_a$'s such that
\begin{equation}\label{eins}
(\forall f \in N \cap {}^a a)\bigvee_{\alpha \in a}\bigvee_{{\bf g}_a \in \cG} f R_{\alpha,a} {\bf g}_a.
\end{equation}

Is there a $(N,\bP,p)$-generic condition such that
$q$ forces that $q$ forces \eqref{eins} holds for
$f\in N[G]$?
Suppose the answer is positive for each iterand, 
what can be said about the countable support limit?

The proof of the iteration theorem in its basic form  uses only
that 
\[\{f \in {}^a a \such f R_{\alpha,a} \bg_a\}\]
is closed for each $R_{\alpha, a}$, $\bg_a$.
There is an example of relations
$R_{\alpha,a}$, $\bg_a$, such that 
$(\bar{R}, \cS, \bar{\bg})$ preserving coincides with 
$(T,\cS,Y)$-preserving.

The reader can jump ahead to Definition~\ref{4.9}
to see what particular $a$, $g_a$  and $R_{\alpha,a}$
we use for the proof of Theorem~\ref{4.2}.
We let
\[\alpha(\ast)= \bigcup \cS=\omega_1\]

We point to the sources
in \cite{Sh:f}:
Our presentation belongs to Case (b) from
\cite[Ch.~XVIII, Context 3.1]{Sh:f}.
Within this Case (b) we focus onto the
Possibilities (also sometimes called ``Cases'' there) A and C in
\cite[Ch.~XVIII, Def.~3.3, 3.4]{Sh:f}.
\relax From Def.~\ref{4.9} one reads off that in the triple
$(\bar{R},\cS, \bar{\bg})$ describing $(T,Y,\cS)$-preservation has closed
and open relations $\bar{R}$. 
So for the relevant relations $\bar{R}$,
preservation for Possibility A and for
Possibility C are equivalent, see  Lemma~\ref{4.11}.

 In our presentation 
in contrast to
\cite{Goldstern93, BJ}, $\bar{R}$ is of size $\omega_1$ and
$\alpha(*) = \omega_1$. The widely known presentations of iteration 
theorems for the relations on the real numbers
\cite{Goldstern93, BJ} have usually countably many relations and
the equivalent to $\alpha(\ast)$  is $\omega$.
The relations we preserve are still on the
Baire space and its topology matters for all the considered
possibilities. However, there is for each $a \in \cS$ an incarnation
${}^a a$ of the Baire space. Moreover, since $a = N \not\in N$ 
names $f \in N$ for functions $f \colon \dom(f) \to \omega_1$
such that $a \subseteq\dom(f)$, $f(x) \in a$ for $x\in a$ now necessarily
are names for functions with larger domains. This does not cause problems,
since the evaluation will be always on $a$.

\medskip

Note that we change one Definition, namely \cite[Ch.~XVIII,
Def.~3.4]{Sh:f}. So our version of ``$(\bar{R},\cS,\bg)$-preserving
for Possibility A'' has not been named in the definitions in the
book nor anywhere else. 
Namely items (iv) and (vi) of Def.~\ref{4.5} are new.

However, this Def.~\ref{4.5} is the one used in the proof of
the preservation theorem for the limit case in \cite[Ch.~XVIII,
Theorem~3.6]{Sh:f}. Our Possibility A here and the proofs here
(which are the ones from the book with some additional explanations)
do not need the distinction whether reals are added or not. The
original definition of Possibility A in the book works as well,
however, the proofs are longer.
There are two proofs based on the old definition:
For forcings that add reals the technique is much shorter
(\cite{Goldstern93})  than for the 
general case that was proved later by Goldstern and
Kellner \cite{GoldsternKellner}. Our proof given here is short and works
in the general case.

\smallskip  

The letter $\bQ$ now stands for an iterand. We let $(2^{|\bQ|})^+ <
\chi$, $\cS \subseteq [\omega_1]^{<\omega_1}$ be stationary, usually
the elements of $\cS$ are of the form
 $N \cap \omega_1$ for a countable $N \prec
\cH(\chi)$. In the language of \cite{Sh:f}, we have for $a \in \cS$,
$d[a], c[a] =  a \not \in a$ and we are in Case (b) of
\cite[Ch.~XVIII, Context 3.1]{Sh:f}, $d[a] \not\in a$, and
$d'[a]=c'[a]=\omega_1$. We will not mention the functions $c$, $d$,
$c'$, $d'$ henceforth since they are fixed. We stay with our special case
of $\cS$ and $\alpha(\ast)= \omega_1$. So we cut down
a lot in comparison to the rich Section~3 of Chapter XVIII of
\cite{Sh:f}. On the other hand, we add numerous proofs to claims
that are written there without a proof.

\begin{definition}\label{4.3}%(Part of \cite[Ch.~XVIII, Def.~3.2]{Sh:f})
\begin{myrules1}
\item[(0)]
$N$ is $(\bar{R},\cS, \bg)$-good means: $a:= N
\cap \bigcup\cS \in \cS$, and for every $f \in N$, $f \colon \bigcup
\cS  \supset\!\to \bigcup \cS$ with $a \subseteq \dom(f)$ for some
$\beta \in a \cap \alpha(\ast)$ we have $f\rest a R_{\beta,a}
\bg_a$.
\item[(1)]
We say {\em $(\bar{R},\cS,\bg)$ covers in $\bV$} iff for
sufficiently large $\chi$ for every $x \in H(\chi)^\bV$ there is a
countable $N \prec \cH(\chi)$ to which $(\bar{R}, \cS, \bg)$ and $x$
belong such that $N$ is $(\bar{R},\cS, \bg)$-good.
\item[(2)]
Let $\cS$ be stationary.
We say {\em $(\bar{R},\cS, \bg)$ fully covers in $\bV$} iff for some
$x \in H(\chi)$,  for every countable $N \prec \cH(\chi)$ to which
$(\bar{R},\cS,\bg)$ and $x$ belong and which fulfils $N \cap \bigcup
\cS \in \cS$ we have that $N$ is $(\bar{R},\cS,\bg)$-good.
\end{myrules1}
\end{definition}

\begin{definition}\label{4.4}%(Part of  \cite[Ch.~XVIII, Def.~3.3]{Sh:f})
\begin{myrules1}
\item[(1)]
We say {\em $(\bar{R}, \cS, \bg)$ strongly covers in Case A} iff it
covers in $\bV$ and each $R_{\alpha,a}$ \nothing{(strictly speaking
$R_{\alpha,a}$ but we will have $\forall \alpha \forall a,b ( \alpha
\in a \cap b \rightarrow (R_{\alpha,a} \leftrightarrow
R_{\alpha,b})))$} is a closed or an open binary relation on ${}^a a$.
We assume from now on that for $\alpha \in a$, $a \in \cS$, 
\[\{f \colon a \to a \such f R_{\alpha,a} \bg_a\}\]
is closed. This is sufficient.

\item[(2)]
We say {\em $(\bar{R}, \cS, \bg)$ strongly covers in Case C} iff it
covers in $\bV$ and in addition for each $a \in \cS$ in any forcing
extension (or at least for any forcing extension in by a forcing in
a family of forcings we are interested in) of $\bV$ player II has a
winning strategy in the following game: In the $n$-th move player I
chooses $N_n$, $H_n$ such that
\begin{myrules1}
\item[(a)] $(N_n,\in \rest N_n)$ is a countable not
necessarily transitive model of ${\rm ZFC}^-$, $N_n \cap \bigcup \cS
= a \in \cS$, $a$, $\cS$, $\bg$, $\bar{R} \in N_n$, $\ell< n
\rightarrow N_\ell \subset N_n$ and $N_n \models (\bar{R}, \cS,
\bg)$ covers, and $f \in N_n \rightarrow (f \rest a) R_a \bg_a$.
\item[(b)]
$H_n \subseteq \{\langle f_0, \dots, f_{n-1} \rangle \such$ for some
finite $d \subseteq \omega_1$, $(\forall \ell <n) (f_\ell \in {}^d
\omega_1)\}$ and $H_n \in N_n$ is not empty,
\item[(c)]
if $\langle f_0, \dots, f_{n-1} \rangle \in H_n$ and $d \subseteq
\dom(f_0)$ is finite then $\langle f_0\rest d, \dots ,f_{n-1} \rest
d \rangle \in H_n$,
\item[(d)] if $\langle f_0, \dots,  f_{n-1}\rangle \in H_n$ and
$\dom(f_0) \subseteq d$, $d$ finite, $d \subseteq\omega_1$ then for
some $\langle f'_0, \dots, f'_{n-1}\rangle \in H_n$ we have
$\dom(f'_\ell) = d$ and $f_\ell \subseteq f'_\ell$,
\item[(e)] $m<n$ and $\langle f_0, \dots, f_{n-1}\rangle \in H_n
\rightarrow \langle f_0, \dots, f_{m-1}\rangle \in H_m^*$ (see
below).
\end{myrules1}
Player II chooses $\langle f^n_0,\dots, f^n_{n-1}\rangle \in H_n
\cap N_n$, $f^n_\ell \supseteq f^m_\ell$ for $\ell \leq m <n$.
Finally, the definition 
$H^*_n = \{\langle f_0, \dots, f_{n-1}\rangle$ for each $\ell$
the functions $f_\ell$, $f^n_\ell$ are compatible$\}$
completes the induction step.

In the end player II wins if for every $m < \omega$, $\bigcup_{n
\geq m} f^n_m$ is a function with domain $a$ and 
$\bigvee_{\alpha\in a} \bigcup_{m\geq n } f^n_m R_{\alpha,a} \bg_a$
\end{myrules1}

\end{definition}

If $R_{\alpha,a}$ is open then we can write $R_{\alpha,a}
= \bigcup_{n\in \omega} R_{\alpha,n,a}$ where each $R_{\alpha,n,a}$
is closed and use $\omega_1=\alpha(\ast) = \omega \alpha(\ast)$,
$R'_{\omega\alpha + n,a} =R_{\alpha,n,a}$ and work with the closed
relations $R'_{\beta,n}$, $\beta < \alpha(\ast)$.

As we already mentioned, we changed the following definition in
Possibility A in comparison to the definition \cite[Ch.~XVIII, Def.~3.4]{Sh:f}
 in the book, so that
it is fitting to the proof in the book.
The items (iv) and (vi) are changed.

\begin{definition}\label{4.5}
%(Compare to \cite[Ch.~XVIII, Def.~3.4]{Sh:f})
We say $\bQ$ is {\em
$(\bar{R},\cS,\bg)$-preserving for Possibility A} iff the following
holds for any $\chi$, $\chi_1$, $N$, $p\in \bQ\cap N$, $k<\omega$: Assume
\begin{myrules1}
\item[($\ast$)]
\begin{myrules1}
\item[(i)] $\chi_1$ is large enough and $\chi > 2^{\chi_1}$,
\item[(ii)] $N \prec \cH(\chi)$ is countable, $N \cap \bigcup \cS = a \in
\cS$, and $\bQ,\cS, \bg, \chi_1 \in N$,
\item[(iii)] $N$ is $(\bar{R}, \cS, \bg)$-good and $ p \in \bQ \cap
N$,
\item[(iv)] $k \in \omega$ and for $\ell < k$ we have a $\bQ$-name
for a function $\name{f}_\ell \in N$, and $\Vdash_{\bQ}
\dom(\name{f}_\ell) \supseteq a$,
\item[(v)] for $\ell <k$, $m<\omega$, $f^*_{m,\ell}$ is a function
from $a$ to $a$ in $N$,
\item[(vi)] for $n<\omega$, $p \leq p_n \leq p_{n+1}$,
\item[(vii)] for $x \in \dom(f^*_{m,\ell})$, $\ell< k$,
for every $m$ there is $n_0$ such that
for $n \geq n_0$, $p_n \Vdash \name{f}_\ell(x)=
f^*_{m,\ell}(x)$,
\item[(viii)] for $\ell < k$, $m<\omega$, $f^*_{m,\ell} R_{\beta^m_\ell,a} \bg_a$ for some
$\beta^m_\ell \in a$, $\beta^{m+1}_\ell \leq \beta^m_\ell$, and
$\beta^*_\ell = \lim_{m \to \omega} \beta^m_\ell$,
\item[(ix)] if $\cI$ is a dense open set of $\bQ$ and $\cI \in N$, then for some
$n$, $p_n \in \cI$.
\end{myrules1}\end{myrules1}
{\em Then} there is $q \geq p$, $q \in \bQ$ that is
$(N,\bQ)$-generic and
\begin{myrules1}
\item[(a)]
for $\ell < k$, $q \Vdash_\bQ (\exists \gamma_\ell \in  a,
\gamma_\ell \leq \beta^\ast_\ell) (\name{f}_\ell\rest a
R_{\gamma_\ell,a} \bg_a)$,
\item[(b)]
$q \Vdash N[\name{\bG_{\bQ}}]$ is $(\bar{R},\cS,\bg)$-good.
\end{myrules1}
\end{definition}

Note that conclusion (a) expresses a sort of directedness:
${\bf g}_a$ is the same for any $\name{f}_\ell$, $\ell <k$.
We will use the possibility to work with
unboundedly many $k$  in the proof of the preservation of
``$(\bar{R},\cS,\bg)$-preserving for Possibility A''
for iterations when the cofinality of the iteration length is
countable.

\begin{definition}\label{4.6}
We say {\em $\bQ$ is $(\bar{R},\cS,\bg)$-preserving for Possibility
C} iff the following holds: Assume
\begin{myrules1}
\item[(i)] $\chi_1$ is large enough and $\chi > 2^{\chi_1}$,
\item[(ii)] $N \prec \cH(\chi)$ is countable, $N \cap \bigcup \cS = a \in
\cS$, and $\bQ,\cS, \bg, \chi_1 \in N$,
\item[(iii)] $N$ is $(\bar{R}, \cS, \bg)$-good and $p \in \bQ \cap
N$.
\end{myrules1}
{\em Then} there is $q \geq p$, $q \in \bQ$ that is
$(N,\bQ)$-generic and $q \Vdash N[\bG_{\bQ}]$ is
$(\bar{R},\cS,\bg)$-good.
\end{definition}

\begin{lemma}\label{4.7} %(\cite[Ch.~XVIII, Claim 3.5]{Sh:f})
%Assume $(\bar{R}, \cS,\bg)$ covers for possibility.
\begin{myrules1}
\item[1)] If $(\bar{R}, \cS,\bg)$ covers in $\bV$ and $\bQ$ is an
$(\bar{R},\cS,\bg)$-preserving forcing notion (for any Possibility)
then in $\bV^{\bQ}$, $(\bar{R},\cS,\bg)$ still covers.

\item[2)] The property ``$(\bar{R},\cS,\bg)$-preserving for Possibility
A (respectively C)'' is preserved by composition of forcing notions.
\end{myrules1}
\end{lemma}

\proof (1) Let $\bG$ be $\bP$-generic over $\bV$.  $N[\bG] \prec
\cH(\chi)^{\bV[\bG]}$ for $N$ being $(\bar{R},\cS,\bg)$-good  in
$\bV$ is a witness for covering in $\bV[\bG]$. For Possibility A, we
can take $k=0$, so $(\ast)$ is vacuously true. Conclusion (b)
suffices.

\smallskip

(2) The proof for Possibility A  fits the old
(\cite[Ch.~XVIII, Def.~3.4]{Sh:f}) and the new
definition of Possibility A (\ref{4.5}). We fix $\bQ = \bQ_0 \ast
\name{\bQ_1}$, $\chi$, $\chi_1$, $N$, $a$, $k$, $\name{f_\ell}$,
$\beta_\ell$, $f^*_{m,\ell}$ for $\ell < k$, $m < \omega$,
$p=(q^0_0,\name{q^0_1})$, $p^n =p_n = (q^n_0,\name{q^n_1})$ as in
($\ast$) of Definition~\ref{4.5} Possibility A. We take $p^0 = p$.
By condition (vi) of ($\ast$) for each $n < m < \omega$, $q^m_0
\Vdash_{\bQ_0} \name{q_1^0} \leq_{\bQ_1} q^n_1 \leq_{\bQ_1} q^m_1$
hence without loss of generality by clause (ix) of ($\ast$) by
taking different names $\name{q^n_1}$ that are above $q^n_0$ the
same,
\begin{myrules1}
\item[$(\ast)_1$]
$\Vdash_{\bQ_0}\name{q_1^0} \leq_{\bQ_1} \name{q^n_1} \leq_{\bQ_1}
\name{q^m_1}$, and
\item[$(\ast)_2$]
for every $x \in a$ or every sufficiently large $n < \omega$,
$(\emptyset,q^n_1)$ forces $\name{f_\ell}(x)$ to be equal to some
specific $\bQ_0$-name $\name{f'_{n,\ell}(x)} \in N$ for each
$\ell<k$.
\end{myrules1}
\nothing{Now we define $\name{f'_\ell}$ a $\bQ_0$-name of a member
of ${}^a a$ such that $\Vdash_{\bQ_0}$ ``for each $x \in a$ for
every sufficiently large $n$, $q^n_1 \Vdash_{\bQ_1}
\name{f'_\ell}(x) = \name{f_\ell}(x)$''. }

 Since $\bQ_0$ is
$(\bar{R},\cS,\bg)$ preserving there is $q_0 \in \bQ_0$ which is
$(N,\bQ_0)$-generic and is above $q^0_0$ in $\bQ_0$ and forces
$N[\name{\bG_{\bQ_0}}]$ to be $(\bar{R},\cS,\bg)$-good and for
some ${\gamma'}^\ast_\ell \leq \beta^\ast_\ell$, %$\gamma'_\ell \in N$
we have $q_0 \Vdash_{\bQ_0} \bigwedge_{\ell< k} \name{f_\ell}
R_{\gamma'_\ell,a} \bg_a$.

Let $\bG_0 \subseteq \bQ_0$ be generic over $\bV$ and $q_0 \in
\bG_0$. We want to apply Definition~\ref{4.5} with $N[\bG_0]$,
$\name{q^0_1}[\bG_0]$, $\langle \name{q^n_1}[\bG_0] \such n < \omega
\rangle$, $\langle \name{f_\ell}[\bG_0] \such \ell<k\rangle$,
$\langle \name{f'_{n,\ell}}[\bG_0]\such \ell<k, n<\omega\rangle$,
$\langle {\gamma'}^\ast_\ell \such \ell < k \rangle$,
$\name{\bQ_1}[G_0]$ there in $(\ast)$  and check that all the items
are fulfilled.

Clause (i) follows from clause (i) for $\bQ_0 \ast \name{\bQ_1}$,

clause (ii): as $q_0$ is $(N,\bQ_0)$-generic we have $N[\bG_0] \cap
\bigcup \cS = N \cap \bigcup \cS \in \cS$,

clause(iii) holds by the choice of $q_0$ and by conclusion (b) in
Definition~\ref{4.5} for $\bQ_0$,

clause (iv) follows from clause (iv) for $\bQ_0 \ast \name{\bQ_1}$,

clause (v): if $x \in a$ then there are $\ell$ and a $\bQ_0$-name
$\name{\tau} \in N$ such that $\Vdash_{\bQ_0} [q_\ell^1
\Vdash_{\name{\bQ_1}} \name{f'_{m,\ell}}(x) = \name{\tau} \in a]$,
as the set of $(r_0,\name{r_1}) \in \bQ_0 \ast \name{\bQ_1}$ such
that $r_0\Vdash_{\bQ_0} \name{r}_1 \Vdash_{\name{\bQ_1}}
\name{f'_{m,\ell}}(x) = \name{\tau}$ for some $\bQ_0$-name
$\name{\tau}$ is a dense open subset of $\bQ_0 \ast\name{\bQ_1}$
some $(q_\ell^0,q_\ell^1)$ is in it and there is such a
$\name{\tau}$, by properness w.l.o.g. $\name{\tau} \in N$. So
$\name{f'_{m,\ell}}[\bG_0] = \name{\tau}[\bG_0] \in a$.

clause (vi) was ensured by our choice $(\ast)_1$,

%clause (vii) by $(\ast)_2$,

clause (vii) by the choice of $\name{f'_{m,\ell}}$ and $\langle
q^n_1 \such n < \omega \rangle$,

clause (viii) by the choice of $q_0$ and $\gamma_\ell'$,

clause (ix) follows from clause(ix) for  $\bQ_0 \ast \name{\bQ_1}$
and a density argument as in (v). In details: If $N[\bG_0] \models
\cI \subseteq\name{\bQ_1}$ is dense and open, then since $\cI \in
N[\bG_0]$  for some $\name{\cI'} \in N$ we have $\Vdash_{\bQ_0}
\name{\cI'}$ is a dense open subset of $\name{\bQ_1}$ and
$\name{\cI'}[\bG_0] = \cI$. Let $\cJ= \{(r_0,\name{r_1}) \in \bQ_0
\ast \name{\bQ_1} \such \Vdash \name{r_1} \in \name{\cI'}\}$. $\cJ
\in N$ is a dense open subset of $\bQ_0\ast \name{\bQ_1}$. Hence for
every sufficiently large $\ell$, $(q^0_\ell,\name{q^1_\ell}) \in
\cJ$  and so $\name{q^1_\ell}[\bG_0] \in \name{\cI'}[\bG_0] = \cI$
and we finish.

\medskip

The proof for Possibility C is particularly easy: We read the
definition of $(\bar{R},\cS,\bg)$-preserving in this case and see
that given $N \prec \cH(\chi)$, $N \cap \omega_1 \in \cS$, $p  \in
\bQ_0 \ast \name{\bQ_1}\cap N$, $p = (q^0_0,q^0_1)$ there is $q \geq
p$, $q = (q^1_0, q^1_1)$, that is $\bQ_0 \ast \name{\bQ_1}$-generic
and
$$(q^1_0,\name{q^1_1})
\Vdash_{\bQ_0 \ast \name{\bQ_1}} N[\bG_{\bQ_0} \ast
\bG_{\name{\bQ_1}[\bG_{\bQ_0}]}] \mbox{ is
$(\bar{R},\cS,\bg)$-good.}$$ First we take  by the hypothesis on
$\bQ_0$ a $(N,\bQ_0)$-generic condition $q^1_0 \geq q^0_0$ such that
$$q^1_0 \Vdash_{\bQ_0} N[\bG_{\bQ_0]} \mbox{ is
$(\bar{R},\cS,\bg)$-good}.$$ Then we take $\bG_{\bQ_0}$,
$\bQ_0$-generic over $\bV$ such that $q^1_0 \in \bG_{\bQ_0}$. Now in
$\bV[\bG_{\bQ_0}]$, $N[\bG_{\bQ_0} ] \cap \omega_1 \in \cS$ and
hence there is $q^1_1 \geq \name{q^1_0}[\bG_{\bQ_0}]$ such that ,
$$q^1_1 \Vdash_{\name{\bQ_1}[\bG_0]}
N[\bG_{\bQ_0} \ast \bG_{\name{\bQ_1}[\bG_{\bQ_0}]}]
\mbox{ is $(\bar{R},\cS,\bg)$-good.}$$ \proofend

The following theorem is central.

\begin{theorem}\label{4.8}%(Compare with \cite[Ch.~XVIII, Theorem~3.6]{Sh:f})
Suppose that $(\bar{R},\cS,\bg)$ strongly covers in $ \bV$ for
Possibility A (resp.\ C), and that $\bP=\langle \bP_i,\bQ_j \such i
\leq \alpha, j< \alpha \rangle$ is a countable support iteration of
proper $(\bar{R},\cS, \bg)$-preserving forcing notions for
Possibility A (resp.\ C). Then $\bP_\alpha$ is a
$(\bar{R},\cS,\bg)$-preserving forcing notion for Possibility A
(resp.\ C) and $(\bar{R},\cS,\bg)$ strongly covers in $\bV^{\bP}$
for the respective Possibility.
\end{theorem}

\proof We prove by induction of $\zeta\leq \alpha$ that for every
$\xi \leq \zeta$, $\bP_\zeta/\bP_\xi$ is
$(\bar{R},\cS,\bg)$-preserving for Possibility A  (resp.\  C) in
$\bV^{\bP_\xi}$, moreover in Definition~\ref{4.5} we can get
$\dom(q)  \setminus \xi = \zeta \cap N$. For $\zeta=0$ there is
nothing to prove, for $\zeta$ successor we use the previous lemma. So
let $\zeta$ be a limit. We first consider $\cf(\zeta) = \omega$.
We fix a strictly increasing sequence $\la \zeta_\ell \such \ell <\omega\ra$
with $\zeta_0= \xi$ and $\sup \zeta_\ell = \zeta$.

\smallskip

First we consider Possibility A. We let $\{\tau_j \such j\in \omega
\}$ list the $\bP_\zeta$-names of ordinals which belong to $N$. Let
$N$ be $(\bar{R},\cS,\bg)$-good. 
In the following we use the convention that the first index indicates
that we deal with a $\bP_{\zeta_\ell}$-name 
$\tau$ or $\name{f}$ (for a $\bP_\zeta/\bP_{\zeta_\ell}$-name)
and the second index is for the enumeration of the particular subset of $N$.

We choose by induction on $j$, 
$k_j<\omega$ such that
\begin{myrules1}
\item[(A)] $k_j < k_{j+1}$,
\item[(B)] there is a sequence $\langle \name{\tau_{\ell,j} }\such
\ell<j\rangle$ such that $\name{\tau_{\ell,j}}$ is a
$\bP_{\zeta_\ell}$-name and
\begin{myrules1}
\item[($\alpha$)]  $p_{k_j} \rest[\zeta_j,\zeta) \Vdash_{\bP_\zeta}
\name{\tau_j} = \name{\tau_{j,j}}$,
\item[($\beta$)] for $\ell<j$ we have $p_{k_j}
\rest[\zeta_\ell,\zeta_{\ell+1}) \Vdash_{\bP_{\zeta_{\ell+1}}}
\name{\tau_{\ell+1,j}} = \name{\tau_{\ell,j}}$,
\end{myrules1}
\item[(C)] if $j= i+1$, $\ell<i$ then $\Vdash_{\bP_{\zeta_{\ell+1}}} p_{k_i}
\rest [\zeta_\ell,\zeta_{\ell+1}) \leq p_{k_j}
\rest[\zeta_\ell,\zeta_{\ell+1})$,
\item[(D)] if $j=i+1$ then
$\Vdash_{\bP_\zeta} p_{k_i} \rest[\zeta_i,\zeta) \leq p_{k_j}
\rest[\zeta_i,\zeta)$.
\end{myrules1}

Given $k_i$, $\la \tau_{\ell,i} \such \ell<i\ra$ we by induction
hypothesis the $p$ that fulfil the requirement for
$p_{k_{i+1}}$ are dense in $\bQ\cap N$, hence by $(ix)$
there is a $k_{i+1}$ such that $p_{k_{i+1}}$ is in that dense set.

Now let
 $\name{f_\ell}$, $\ell < k$, be given as in ($\ast$) of Def.~\ref{4.5}.
Let  $\{\name{f_j} \such \ell< j<\omega\}$ list the $\bP_\zeta$-names
of members on $N$ that are extensions of functions from $a$ to $a$.
For $\ell<k$ let them be the $f^\ast_{m,\ell}$ as given in ($\ast$) of
Definition~\ref{4.5}.  Since $N$ is $(\bar{R},\cS,\bg)$-good, $p_n$
from above can serve as $p_n$ in ($\ast$). We will now show how to
choose $\name{f^\ast_{m,j}} \in N$, $m<\omega$, $j < \omega$.

Let $h(j,x) <\omega$ be such that $\name{\tau_{h(j,x)}} =
\name{f_j}(x)$.  We can now define for $n< \omega$, $j<\omega$,
$\name{f^\ast_{n,j}}$ a $\bP_{\zeta_n}$-name of a function from $a$
to $a$. Let $\name{f^\ast_{n,j}}(x)= \tau_{n,h(j,x)}$ if $h(j,x)
\geq n$ and $\tau_{h(j,x),h(j,x)}$ if $h(j,x)<n$. So 
$\name{f^\ast_{0,j}}(x)= f_j(x)$ for $j < k$.
So also for the names $\name{f^\ast_{m,j}}(x)$ we have (viii) of the
hypothesis ($\ast$), since (viii) holds objects $f_{m,j}$ from there
and the sets
\[ \{f \in {}^a a \such f R_{\beta ,a} \bg_a\} 
\]
are closed for $\beta \in a$.
We choose by induction on $n$, $q_n$,  $\name{\alpha^n_\ell}$ for
$\ell < k+n$ such that
\begin{myrules1}
\item[(a)] $q_n \in \bP_{\zeta_n}$, $\dom(q_n) \setminus \xi = N \cap
\zeta_n$, $q_{n+1}\rest \zeta_n = q_n$,
\item[(b)]
$q_n$ is $(N,\bP_{\zeta_n})$-generic,
\item[(c)]
$q_n \Vdash_{\bP_{\zeta_n}} N[\bG_{\bP_{\zeta_n}}]$ is
$(\bar{R},\cS,\bg)$-good,
\item[(d)] $p_0 \rest \zeta_0 \leq q_0$ in $\bP_{\zeta_0}$,
\item[(e)] $q_{n+1} \rest \zeta_n \Vdash_{\bP_{\zeta_n}} p_{n}
\rest[\zeta_n,\zeta_{n+1}) \leq p_{n+1} \rest[\zeta_n,\zeta_{n+1})$
(in $\bP_{\zeta_{n+1}}/\bP_{\zeta_n}$),
\item[(f)] for $\ell<k+n$, $\name{\alpha^n_\ell}$ is a $\bP_{\zeta_n}$-name
of an ordinal in $a$, $q_{n+1} \Vdash \name{\alpha^{n+1}_\ell} \leq
\name{\alpha^n_\ell}$, $\alpha^0_\ell \leq \beta^\ast_\ell$, for
$\ell<k$,
\item[(g)] for $\ell<k+n$, $q_n \Vdash_{\bP_{\zeta_n}}
\name{f^\ast_{n,\ell}} R_{\name{\alpha^n_\ell},a } \bg_a $.
\end{myrules1}

The induction step is by the induction hypothesis and by
Definition~\ref{4.5} Possibility A with $k+n$ in the role of
$k$. In the end we let $q =
\bigcup_{n<\omega}q_n$.

We show that $q$ is $(N,\bP_\zeta)$-generic and that is satisfies
conclusions (a) and (b) of Def.~\ref{4.5}. Let $q \in
\bG_{\bP_\zeta} \subseteq \bP_\zeta$, $\bG_{\bP_\zeta}$ be
$\bP_\zeta$-generic over $\bV$. $\bG_{\bP_\xi} = \bG_{\bP_\zeta}
\cap \bP_\xi$ for $\xi < \zeta$ and $\bG_{\bP_{\zeta_n}} =
\bG_{\bP_\zeta} \cap \bP_{\zeta_n}$. Now for each $\bP_\zeta$-name
$\name{\tau}$ for an ordinal there is some $j$ such that
$\name{\tau}=\name{\tau_j}$. $q_j$ forces $\name{\tau_{j,j}} \in N$
and $p_j\rest [\zeta_j,\zeta)$ forces $\name{\tau_j} =
\name{\tau_{j,j}}$. $p_j \rest[\zeta_j,\zeta) \leq q$ by (d) and
(e). So $q$ forces $\name{\tau_j} = \name{\tau_{j,j}}$ and $q \Vdash
\name{\tau_j}\in N \cap \Ord$, so $q$ is $(N,\bP_\zeta)$-generic.

For each $\ell$, $\langle \alpha^n_\ell \such \ell \leq n <
\omega\rangle$ is not increasing by (f) and hence eventually
constant, say with value $\alpha^*_\ell$. If $x \in a$, $j<\omega$,
then for $n > h(j,x)$, $p_n \Vdash \name{f_j}(x)
=\name{f^\ast_{n,j}}(x)$. So for every finite $b \subseteq a$,
$\langle(\name{f^\ast_{n,j}}\rest b) [\bG_{\bP_{\zeta_n}}] \such n <
\omega \rangle$ is eventually constant, equal to $ (\name{f_j}\rest
b)[\bG_{\bP_\zeta}]$.  By (g), for sufficiently large  $n$, 
\begin{myrules}
\item[(1)]
$q
\Vdash_{\bP_\zeta} (\name{f_j}\rest b) [\name{\bG_{\bP_\zeta}}] =
(\name{f^\ast_{n,j}} \rest b )[\name{\bG_{\bP_{\zeta_n}}}]$ and
\item[(2)]
$q_n\Vdash \name{f^\ast_{n,j}}[\name{\bG_{\bP_{\zeta_n}}}]
R_{\alpha^n_j ,a} \bg_a$ and 
\item[(3)] $\alpha^n_j = \alpha^\ast_j$.
\end{myrules}
 Since $R_{\alpha^\ast_j,a}$ is closed,
and $\la \name{f^\ast_j} \rest b [\bG_{\bP_\zeta}] \such b \subseteq a, b 
\mbox{ finite} \ra$
converges to $\name{f_j}$, we get 
$q\Vdash_{\bP_\zeta} \name{f_j}%[\name{\bG_{\bP_\zeta}}]
R_{\alpha^\ast_j, a} \bg_a$. 
This finishes the proof of (b), that $q
\Vdash_{\bP_\zeta} N[\name{\bG_{\bP_\zeta}}]$ is
$(\bar{R},\cS,\bg)$-good. Now for (a) note that for there is $n$
such that  for $\ell< k$,
$p_n \Vdash \alpha^\ast_\ell \leq \alpha_\ell^n[\bG_{\bP_\zeta}] \leq
\beta^\ast_\ell$. Thus we finished the proof for the limit of
countable cofinality for Possibility A.

\medskip

Again the proof for possibility C is short. Let
$\langle \name{f_\ell} \such \ell<\omega\rangle$ enumerate the
$\bP_\zeta$-names $\name{f} \colon \omega_1 \to \omega_1$ with
$\name{f} \in N$. Let $\langle\name{\tau_n} \such n < \omega
\rangle$ list the $\bP_\zeta$-names of ordinals which belong to $N$.
We choose by induction on $n$, $\name{p_n}$, $q_n$,  $\name{H_n}$,
$\langle \name{f_\ell^n} \such \ell \leq n \rangle$ such that
\begin{myrules1}
\item[(a)] $q_n \in \bP_{\zeta_n}$, $\dom(q_n) \setminus \xi \subseteq N \cap
\zeta_n$, $q_{n+1} \rest \zeta_n = q_n$,
\item[(b)] $q_n$ is $(N[\name{\bG_{\bP_{\zeta_n}}}],
\bP_{\zeta_n})$-generic,
\item[(c)]
$q_n \Vdash N[\name{\bG_{\bP_{\zeta_n}}}]$ is $(\bar{R},\cS,
\bg)$-good,
\item[(d)] $\name{p_n}$ is a $\bP_{\zeta_n}$-name of a member of  $\bP_\zeta
\cap N$, $q_n \Vdash_{\bP_{\zeta_n}} \name{p_n} \rest \zeta_n \in
\name{\bG_{\bP_{\zeta_n}}}$,
 \item[(e)] $\name{H_n}$ is a
$\bP_{\zeta_n}$-name, $\name{H_n} = \{ \langle g_0,\dots ,g_{n-1}
\rangle \such d \subseteq a$ is finite and $\name{p_n}
\not\Vdash_{\bP_\zeta/\name{\bG_{\zeta_n}}} \langle \name{f_0}\rest
d, \dots,\name{f_{n-1}}\rest d\rangle \neq \langle
g_0,\dots,g_{n-1}\rangle\}$,
\item[(f)]
$\name{f_\ell^n}$ is a $\bP_{\zeta_n}$-name such that
\begin{equation*}
\begin{split}
q_n \Vdash_{\bP_{\zeta_n}} &
\langle \name{f^n_\ell} \such
\ell<n\rangle \in \name{H_n} \mbox{ and for every $m \leq n$ we
have}\\
& \name{p_{n+1}} \not\Vdash_{\bP_\zeta/\bP_{\zeta_n}} \neg
\bigwedge_{\ell<m} \name{f_\ell} \supseteq \name{f_\ell^m},
\end{split}
\end{equation*}
\item[(g)]
$q_n\Vdash_{\bP_{\zeta_n}} \name{p_{n+1}}$ forces a value to
$\name{\tau_n}$.
\end{myrules1}

There is no problem to carry out the definition and we still have
the freedom to choose $\langle \name{f_\ell^n} \such
\ell<\omega\rangle$. For this we use the winning strategy from
possibility C of Definition~\ref{4.4} choosing there the $n$-th move
of player I as $N_n= N[\bG_{\bP_{\zeta_n}}]$ and
\begin{equation*}
\begin{split}
\name{H_n}[\bG_{\bP_{\zeta_n}}] = &\{\langle
g_0,\dots,g_{n-1}\rangle \such \mbox{ for some finite $d \subseteq
a$ we have }\\
&g_\ell \in {}^d \omega_1\mbox{  for $\ell<n$
and}\\
&\name{p_n}[\bG_{\bP_{\zeta_n}}]
\not\Vdash_{\bP_\zeta/\bG_{\bP_{\zeta_n}}} \langle \name{f_0} \rest
d,\dots, \name{f_{n-1}}\rest d\rangle \neq \langle
g_0,\dots,g_{n-1}\rangle \} \end{split}
\end{equation*}
 so the
$n$-th move is defined in $\bV^{\bP_{\zeta_n}}$ according to the
winning strategy for $\bP_{\zeta_{n+1}}$. We can work in $\bV^{{\rm
Levy}(\aleph_0,(2^{|\bP_\alpha|})^+)}$. Now of course while playing
the universe changes but as the winning strategy is absolute there
is no problem.

\smallskip

We let $q = \bigcup q_n$. In the end player II wins, that means for
every $m < \omega$ $\bigcup_{n \geq n} f_m^n$ is a function which
has domain $a$ and $\bigcup_{n \geq m} f^n_m R_a \bg_a$. By the
choice of players I's moves for every $m < \omega$,
$q\Vdash_{\bP_\zeta} \name{f_m} = \bigcup_{n \geq m} f^n_m$. So $N$
is $(\bar{R},\cS,\bg)$-good. Moreover, also for $\bP_\alpha$ player
II has a winning strategy in the game: It is just the winning
strategy sketched above.
If I plays a proper, non-empty subset of the $\name{H_n}$
from Definition~\ref{4.4} then II takes 
an $n$-tuple from the subset according to 
his winning strategy in $\bV^{\bP_{\zeta_n}}$. 
 Thus we finish the limit step of countable
cofinality.

\medskip

Now we continue to look at Possibility C:

If $\cf(\zeta) > \aleph_0$, suppose an initial part of the game
$\langle \name{N_n}, \name{H_n}, \langle \name{f_\ell^n} \such \ell
\leq n \rangle \such n \leq m \rangle$ for $\bP_\zeta$ is played.
\nothing{ then let $\{\name{f_\ell} \such \ell < \omega\}$ be again
an enumeration of all functions in $N$ and let $\{\name{\tau_\ell}
\such \ell < \omega\}$  be an enumeration of all names in $N$ for
ordinals.
 For each $k < \omega$} We choose $\xi<\zeta$ such these finitely many names
 for countable objects are $\bP_\xi$-names. One can collapse all
 objects and thus get them into $H(\omega_1)$ and hence they are hereditarily
 countable and then \cite[Lemma 5.13]{BsSh:242} is applicable
 and such a
 $\xi$ exists. Now player II can play according to the strategy for
 $\bP_\xi$. Then player I moves and we choose a new $\xi$ for
 catching the longer initial segment of the game. In the end II
 wins, since $\bP_\xi \lessdot \bP_\zeta$ for all the $\xi$'s on the
 way and each $f \in N[\bG_{\bP_\zeta}]$, $\colon a \to a$
 also appears at some stage $\xi<\zeta$, so we know
 then $N[\bG_{\bP_\zeta}]$ is $(\bar{R},\cS,\bg)$-good.

\smallskip

In the case of Possibility A, all the $\bP_\zeta$-names $\name{f}
\in N$ for functions in ${}^a a$ and the $\name{f_\ell}$, $\ell<k$,
are $\bP_\xi$-names for a $\xi<\zeta$. \proofend

\begin{definition}\label{4.9}%(\cite[Ch.~XVIII, Def.~3.9 A]{Sh:f})
Now fix $(T,Y,\cS)$ as in Section~\ref{S1}. Assume that
$\cS\subseteq [\omega_1]^{\omega}$ stationary and $\delta = \sup(a)$
for $a \in \cS$.
$Y \subseteq T$ is given, and we fix an enumeration of it as follows:
$Y(\delta) = \{ t_n^\delta \such n < \gamma_\delta\} \subseteq T_\delta$
and $\gamma_\delta$ is finite or $\omega$.
Now we choose $R_{\alpha,a}$ and ${\bf g}_a$
by defining $\{ f \colon a
\to a \such f R_{\alpha,a} \bg_a\}$ for $\sup(a) = \delta$ as
\begin{myrules1}
\item[($\alpha$)] $\alpha=0$ and  $f(0) \in \gamma_\delta$ 
and $f^{-1}[\{1\}] \cap \{
s\in T_{<\delta}\cap a  \such 0 <_T s <_T t_{f(0)}^\delta \} \neq
\emptyset$, or 
\item[($\beta$)] $0 < \alpha< \delta$ and  $f(0) \in \gamma_\delta$ 
and  $f^{-1}[\{1\}] \cap \{ s\in
T_{<\delta}\cap a  \such t^\delta_{f(0)} \rest \alpha \leq s \} =
\emptyset$
\item[($\gamma$)] $f(0) \not \in \gamma_\delta$. 
\end{myrules1}
$R_{\alpha, a}$ is a countable union of closed relations, so
Possibility A applies.
\end{definition}

\begin{lemma}\label{4.10} %(\cite[Ch.~XVIII, Claim 3.9 2) and 3)]{Sh:f})
 1) Iff $T$ is
$(Y,\cS)$-proper, then $(\bar{R},\cS,\bg)$ fully covers.

2) If $(\bar{R},\cS,\bg)$ covers then $(\bar{R},\cS,\bg)$ strongly
covers for Possibility A.
\end{lemma}

\proof 1) We read the meaning of $\{ f \in {}^a a \such f R_{\beta,a}
\bg_a\}$ from Definition~\ref{4.9}.
$f(0) \in a \cap \omega_1$.
The part (a) of the disjunction means $f^{-1}[\{1\}]$ is a subset of
the forcing $T_{<\delta}$, and 
$\{s \such s <_t t^\delta_{f(0)}\}$ meets $f^{-1}[\{1\}]$.
This is an open relation.
Note that it is not written that $f^{-1}[\{1\}]$ be dense.
The disjunction (b) means $f^{-1}[\{1\}]$ is not dense in $(T\cap N, <_T)$
since above $t^\delta_{f(0)}\rest \alpha$ there is no element.
This is a closed relation.
The disjunction (c) means 
$t_{f(0)}^\delta$ need not be considered as a generic filter in 
$T_{<\delta}$ as it is not
in $Y(\delta$). This is a clopen relation, since it speaks only about
$f(0)$.
Now $(\bar{R},\cS,\bg)$ covers $N$ iff $T$ is $(Y,\cS)$ proper. 
2) The $R_{\beta,a}$ are open or closed. \proofend

\begin{lemma}\label{4.11}%(\cite[Ch.~XVIII Claim 3.9 C]{Sh:f}) Let
$(\bar{R},\cS,\bg)$ be as in Definition~\ref{4.9}. A forcing notion
$\bQ$ is $(\bar{R},\cS,\bg)$-preserving for Possibility C iff it is
so for Possibility A.
\end{lemma}

\proof ``$\leftarrow$'': A winning strategy  for player II is: In
the $n$-th round he uses that I played $N_n \supseteq N_\ell$ and
$H_n$ for $n \geq \ell$ such that $N_n$ is as in (a) of
Definition~\ref{4.4}(2). In particular $N_n \models
(\bar{R},\cS,\bg)$ covers. So there is $N' \in N_n$ that is
$(\bar{R},\cS,\bg)$-good and $N' \models (\bar{R},\cS,\bg)$ covers
and is $(\bar{R},\cS,\bg)$-good.

Player II plays
 $\langle f_0^n,\dots,f^n_{n-1}\rangle\in H_n\cap N_n$ and $d_n\supseteq
d_\ell$ for $\ell <n$ and $\beta^n_{k} \leq \beta^\ell_{k}$ for $k
\leq n$ and $\ell< n$ such that $\dom(f_\ell^n) = d_n$, for $\ell<n$
and such that  $(\forall \ell< n)(\exists m)(f^\ast_{m,\ell} \supset
f^n_\ell \rest d_n \wedge f^\ast_{m,\ell} \colon a \to a \wedge
f^\ast_{m,\ell} R_{\beta^m_{\ell},a} {\bf g}_a)$, and such that $\bigcup
d_n = a$. Since for $\ell < \omega$, the $\beta^m_{\ell}$, $m <
\omega$, become eventually constant to $\beta^*_\ell$, and since
$R_{\beta^*_\ell,a}$ is closed, he thus ensures $(\forall \ell \in
\omega) (\bigcup_{n \geq \ell} f^n_\ell R_{\beta^*_\ell,a} \bg_a)$.

\smallskip

``$\rightarrow$'': We take the conditions from Possibility A: Let
$N\prec \cH(\chi)$ be countable, $(\bar{R},\cS,\bg)\in N$ and $p$,
$\langle p_n \such n<\omega\rangle$, $\langle \name{f_\ell} \such
\ell<k \rangle$, $\langle f^\ast_{m,\ell} \such \ell<k, m<\omega
\rangle$, $\langle \beta^m_\ell, \beta^*_\ell\such \ell<k, m<\omega
\rangle$ be as in ($\ast$) of Definition~\ref{4.5}.

Let $\delta = N \cap \omega_1$. Let $w_m = \{ \ell< k \such
f^\ast_{m,\ell}(0) < \gamma_\delta \wedge \beta^m_\ell \neq 0 \}$.
For some $m_0$ all of the finitely many possible $w_m$ appeared. Fix
such an $m_0$ and fix the finitely many witnesses $x_{m,\ell}$, $m <
m_0$, $\ell < k$. For $\ell \in k \setminus w_m$ choose
$t_{f^\ast_{m,\ell}(0)}^\delta\rest \beta^m_\ell \leq_T
x_{m,\ell}\in T\cap N$ such that $x_{m,\ell} <_T
t^\delta_{f^*_{m,\ell}(0)}$ and $\bigvee p_n \Vdash
\mbox{``}\name{f_\ell}(x_{m,\ell})=1 \vee \name{f_\ell}(0) \geq
\gamma_\ell$.'' So for some $n(\ast)$,
$$
p_{n(\ast)} \Vdash \bigwedge_{\ell \in k \setminus w}
\name{f_\ell}(x_{m,\ell})= 1 \vee \name{f_\ell}(0) \geq
\gamma_\ell.$$
Let \begin{equation*} \begin{split} \cI = \{ q \in
\bQ \such&
 \mbox{
for each }\ell \in w, q \mbox{ forces a value to }
\name{f_\ell},\mbox{ say $m_\ell$ and it forces}\\
& \mbox{a truth value to } (\exists x)(t^\delta_{m_\ell}\rest
\beta_\ell <_T x \wedge \name{f_\ell}(x)=1) \}. \end{split}
\end{equation*}
 So for some $n
> n(\ast)$ we have $p_n \in \cI$, and hence all truth values it forces
are false, since for $f^\ast_{m,\ell}$ they are false, since
$f^\ast_\ell R_{\beta^\ast_\ell,a} \bg_a$ and since for any $x$
there is $p_{n'} \geq p_n$ and $m$ forcing $\name{f_\ell}(x)=
f^\ast_{m,\ell}(x)$ and we let $\beta^{\ast\ast}_\ell =
\max\{\beta^m_\ell \such m < m_0\}$. So any $(N,\bQ)$-generic $q
\geq p_{n'}$ we have
\[
q\Vdash \bigwedge_{m \in
\gamma_\delta} \{s \in \delta \such s <_t t_m^\delta\} \mbox{ is
$(N[\name{\bG_\bQ}],T)$ generic},
\]
 or, in other words, $q \Vdash
N[\name{\bG_{\bQ}}]$ is $(\bar{R},\cS, \bg)$-good. The existence of
such a $q$ follows from Possibility C. For $\name{f_\ell}$, $\ell <
k$ moreover $q \Vdash \name{f_\ell} R_{\beta^{\ast\ast}_\ell, a}
\bg_a$ and $\beta^{\ast\ast}_\ell \leq \gamma_\ell$. So the
conclusions (a) and (b) of Def.~ \ref{4.5} of preserving in Case A
are shown.

\proofend

\begin{lemma}\label{4.12}%(\cite[Ch.~XVIII Fact 3.9 E]{Sh:f}) Let
$(\bar{R},\cS,\bg)$ be as in Definition~\ref{4.9}. A forcing notion
$\bQ$ is $(\bar{R},\cS,\bg)$-preserving for Possibility A iff it is
$(T,\cS,Y)$-preserving.
\end{lemma}

\proof For the forward direction we read the meaning of
$(\bar{R},\cS,\bg)$. It is easy to fulfil ($\ast$) of
Definition~\ref{4.5} for $k=0$. Conclusion (b) in
Definition~\ref{4.5} for Possibility A ensures that $\bQ$ is
$(T,\cS,Y)$-preserving.

\smallskip

For the backward direction: We look at Possibility A. All the $\{f
\such f R_{\alpha,a} \bg_a\}$ are open or closed, and by
reorganising, we can assume that all of them are closed. Fix $k \in
\omega$ and $p$, $\name{f_\ell}$, $\ell< k$, $f^\ast_{m,\ell}$,
$\ell , m < \omega$, $\beta^m_\ell$, $m< \omega$, $\langle p_n \such
n < \omega \rangle$ as in ($\ast$) of Def.~\ref{4.5}. We have to
find a particular generic $q \geq p$ that also  satisfies
conclusions (a) and (b) of Def.~\ref{4.5}.

Now $$ \cI=\{ q \in \bQ \cap N \such  q \not\Vdash (\exists \ell <
k) (\neg \name{f_\ell} R_{\beta^\ast_\ell,a} \bg_a)\}$$ is dense
above some $p_n$ in $\bP$ since $\bigcup_{\ell< k} \{ f \in {}^a a
\such \neg f R_{\beta^\ast_\ell,a} \bg_a \}$ is open in the Baire
space ${}^a a$. So first we take $p_1 \geq p$, $p_1 \in \cI \cap N$,
and then we take, according to preserving in possibility C,  $q \geq
p_1$ such that $q$ is $(N,\bP)$-generic and $q \Vdash
N[\name{\bG_\bQ}]$ is $(\bar{R}, \cS,\bg)$-good.
 \proofend

So now for the proof of Theorem~\ref{4.2} we can use Lemma~\ref{4.6}
and Theorem~\ref{4.8} for Possibility A or for Possibility C, and we
are done.

\def\germ{\frak} \def\scr{\cal} \ifx\documentclass\undefinedcs
  \def\bf{\fam\bffam\tenbf}\def\rm{\fam0\tenrm}\fi % f**k-amstex!
  \def\defaultdefine#1#2{\expandafter\ifx\csname#1\endcsname\relax
  \expandafter\def\csname#1\endcsname{#2}\fi} \defaultdefine{Bbb}{\bf}
  \defaultdefine{frak}{\bf} \defaultdefine{=}{\B} % doublef**k-amstex!!
  \defaultdefine{mathfrak}{\frak} \defaultdefine{mathbb}{\bf}
  \defaultdefine{mathcal}{\cal}
  \defaultdefine{beth}{BETH}\defaultdefine{cal}{\bf} \def\bbfI{{\Bbb I}}
  \def\mbox{\hbox} \def\text{\hbox} \def\om{\omega} \def\Cal#1{{\bf #1}}
  \def\pcf{pcf} \defaultdefine{cf}{cf} \defaultdefine{reals}{{\Bbb R}}
  \defaultdefine{real}{{\Bbb R}} \def\restriction{{|}} \def\club{CLUB}
  \def\w{\omega} \def\exist{\exists} \def\se{{\germ se}} \def\bb{{\bf b}}
  \def\equivalence{\equiv} \let\lt< \let\gt>

%\bibliographystyle{plain}
%\bibliography{../freiburg/sh/lit,../freiburg/sh/listb}

\end{document}